\documentclass[10pt]{article}
\usepackage[margin=2.5cm]{geometry}
\usepackage[utf8]{inputenc}
\usepackage{bussproofs}
\usepackage{calc}
\usepackage{amsmath}
\usepackage{amsthm}
\usepackage{amssymb}
\usepackage{stmaryrd}
\usepackage{setspace}
\usepackage{MnSymbol}
\usepackage{cmll}
\usepackage{longtable}
\theoremstyle{definition}
\newtheorem{thm}{Theorem}
\newtheorem{cnv}[thm]{Convention}
\newtheorem{defn}[thm]{Definition}
\newtheorem{crl}[thm]{Corollary}

\newtheorem{prp}[thm]{Proposition}

\usepackage{setspace}
\title{Calculi of epistemic grounding based on Prawitz's theory of grounds}
\author{Antonio Piccolomini d'Aragona\\\begin{small}Aix Marseille Univ, CNRS, Centre Gilles Gaston Granger, Aix-en-Provence, France\end{small}\\\texttt{antonio.piccolomini-d-aragona@univ-amu.fr}}
\date{}

\begin{document}

\maketitle

\begin{abstract}
We define a class of formal systems inspired by Prawitz's theory of grounds. The latter is a semantics that aims at accounting for epistemic grounding, namely, at explaining why and how deductively valid inferences have the power to epistemically compel to accept the conclusion. Validity is defined in terms of typed objects, called grounds, that reify evidence for given judgments. An inference is valid when a function exists from grounds for the premises to grounds for the conclusion. Grounds are described by formal terms, either directly when the terms are in canonical form, or indirectly when they are in non-canonical form. Non-canonical terms must reduce to canonical form, and two terms may be said to be equal when they converge towards equivalent grounds. In our systems these properties can be proved through rules distinguished according to whether they concern types or logic. Type rules involve type introduction and elimination, equality for application of operational symbols, and re-writing equations for non-canonical terms. The logic amounts to a sort of intuitionistic system in a Gentzen format. To conclude, we show that each system of our class enjoys a normalization property.
\end{abstract}

\medskip


\section{Introduction}

Prawitz's theory of grounds (ToG) is a constructive semantics belonging to the field of proof-theoretic semantics [the expression is due to Schroeder-Heister 1991; for a general introduction see Francez 2015, Schroeder-Heister 2018]. It aims at providing a notion of valid inference through which the epistemic power of deduction could be adequately explained. The basic idea is that one is justified in asserting something when in possession of a ground for this assertion. Inferences are applications of operations on grounds for the premises, and an inference is valid if the operation applied in it actually yields a ground for the conclusion.

Prawitz has devoted many papers to a \emph{philosophical} articulation of ToG [Prawitz 2009, 2012, 2013, 2015, 2018]. The \emph{formal} part, however, is only at an embryonic stage. From Prawitz's examples and hints, it seems that ToG should at least include a theory of \emph{formal languages of grounding}, whose terms denote grounds and operations on grounds. Here, we build upon this idea and expand it towards a theory of \emph{formal systems of grounding}, through which we can prove that terms of languages of grounding enjoy some relevant properties.

The languages of grounding envisaged by Prawitz seem to contain terms, but not formulas. Our first step will thus be that of showing how term-languages can be expanded by adding to them two binary predicates and the first-order logical constants. The predicates express, respectively, the fact that a term denotes a ground for a certain assertion, and the fact that two terms denote the same ground for a certain assertion. An expansion of this kind is exemplified in Section 3, over a ground-theoretic version of Gentzen's natural deduction system for first-order intuitionistic logic. In Section 4 we introduce a system of grounding for the enriched Gentzen-language, and prove some theorems within it, as well as a meta-theorem about it. In Section 5 we generalize the Gentzen-based example, and outline a hierarchy of systems of grounding where, \emph{mutatis mutandis}, we can prove all the results already provable with respect to the example of Section 4. Finally, in Section 6 we define reduction and permutation functions for non-normal derivations, and a double measure for maximal formulas. We then prove a normalization result that can be applied to each system in the hierarchy.

Before all this, however, we must mention some issues from the background of our approach, that is, from Prawitz's own ground-theoretic ideas, and from Prawitz's proposal of developing languages of grounding whose terms denote grounds and operations on grounds.

\section{Background: grounds and languages of grounding}

ToG aims at explaining why and how valid inferences have the epistemic power to compel us to accept their conclusion if we have accepted their premises as justified. Since a proof can be understood as a chain of valid inferences, ToG also aims at accounting for the epistemic power of deduction. Prawitz's previous semantics of valid arguments and proofs [mainly Prawitz 1971, 1973, 1977] proves unable to accomplish this task, mainly because of an inversion in the explanatory order according to which proofs are chains of valid inferences. In Prawitz's earlier approach, indeed, an inference is valid when it yields valid arguments or proofs when attached to valid arguments or proofs, so that valid arguments and proofs are prior and cannot be defined through valid inferences [an analysis of this point, as well as of the other reasons that led Prawitz to the adoption of ToG, can be found in d'Aragona 2019a, 2021b, Usberti 2015].

Thus, in ToG Prawitz introduces a non-inferential notion of ground that provides a non-circular definition of inferential validity [for this point see mainly Prawitz 2019]. Informally, a ground is an abstract object that reifies what one is in possession of when one has evidence for a given assertion. This naturally leads to consider a ground as an object typed on the sentence for the assertion of which it is a ground. Since evidence states are achieved by performing operations of some kind, grounds must be conceived also as epistemic objects, namely, as objects built up of the operations performed to achieve the evidence they stand for. Also the operations are typed, although their types involve a domain - the input values of the operations - and a co-domain - the output grounds they produce.

Of course, evidence for assertions must somehow depend on the meaning of the asserted sentences. And since we want to explain possession of evidence as possession of grounds, a connection must hold between grounds and meaning. Prawitz's line of thought seems to be inspired here by the BHK explanation of meaning in terms of proof-conditions for formulas of different logical kinds [Heyting 1956], as well as by Gentzen's claim that introduction rules in his natural deduction calculus define the meaning of the logical constant they introduce [Gentzen 1934 - 1935]. In ToG, this becomes the idea that the meaning of a sentence is given by laying down what counts as a ground for asserting it. In the atomic case, these ground-conditions vary depending on the content of the assertion. Atomic grounds are thus given by so-called atomic bases which, as known, also occurred in Prawitz's previous semantic approaches [mainly Prawitz 1973]. As we shall see, an atomic base consists of the sets of the individual, functional and relational symbols of a given language (usually, a first-order one) plus an atomic system on this language. Atomic systems are normally understood as sets of production rules whose premises and conclusion are atomic formulas of the underlying language. The general idea is therefore that a production rule involving an individual, functional or relational symbol $s$ provides the meaning of $s$ by showing its deductive behavior. So, atomic bases play in Prawitz's semantics the same role as models in model-theory. The main difference is that, while in the latter the meaning of the individual, functional and relational symbols is given through interpretation functions that associate these symbols to set-theoretic objects, in a Prawitzian framework this is done in deductive terms, in compliance with the epistemic concerns of the overall constructivist setup.

In the logically complex cases, Prawitz provides clauses - over atomic bases $\mathfrak{B}$ - that proceed by induction on logical complexity, by attributing to each logical constant $k$ a primitive and meaning-constitutive operation $k I$. The simplest clauses are those for $\wedge$, $\vee$ and $\exists$. Given $\alpha, \beta$ ($i = 1, 2$) and $t$ closed - i.e. $t$ does not involve individual variables -

\begin{itemize}
    \item[($\wedge$)] $g$ is a ground for $\alpha \wedge \beta$ over $\mathfrak{B}$ iff $g$ is of the form $\wedge I(g_1, g_2)$, where $g_1$ is a ground over $\mathfrak{B}$ for $\alpha$ and $g_2$ is a ground over $\mathfrak{B}$ for $\beta$;
    \item[($\vee$)] $g$ is a ground over $\mathfrak{B}$ for $\alpha \vee \beta$ iff $g$ is of the form $\vee I[\alpha \rhd \alpha \vee \beta](h)$, where $h$ is a ground over $\mathfrak{B}$ for $\alpha$, or of the form $\vee I[\beta \rhd \alpha \vee \beta](h)$, where $h$ is a ground over $\mathfrak{B}$ for $\beta$;
    \item[($\exists$)] $g$ is a ground over $\mathfrak{B}$ for $\exists x \alpha(x)$ iff $g$ is of the form $\exists I[\alpha(t) \rhd \exists x \alpha(x)](h)$, where $h$ is a ground over $\mathfrak{B}$ for $\alpha(t)$.
\end{itemize}
Observe that in ($\vee$) and ($\exists$) we had to specify the type of the operation through an expression of the form $x \rhd y$ where $x$ is the domain and $y$ the co-domain of the operation. Operational types become central when we consider the clauses for $\rightarrow$ and $\forall$, since the latter involve an unrestricted notion of operation on grounds. Clearly, not all operations on grounds are primitive and meaning-constitutive. Others are non-primitive and, Prawitz says, they must be specified by indicating their domain, their co-domain and - since we expect them to be constructive - an equation that shows how their input values can be transformed into grounds of an expected type. In general, given $\alpha_1, ..., \alpha_m, \beta$ open or closed with free individual variables $x_1, ..., x_n$ ($n, m \geq 0$), we can say that a first-level operation on grounds over $\mathfrak{B}$ of operational type

\begin{center}
    $\alpha_1, ..., \alpha_m \rhd \beta$
\end{center}
- here and below, the symbol $\rhd$ is meant to separate the domain $\alpha_1, ..., \alpha_m$ from the co-domain $\beta$ - is a constructive function that can be indicated by employing variables $x$ for individuals and typed-variables $\xi^\alpha$ for grounds for $\alpha$, i.e.

\begin{center}
    $f(x_1, ..., x_n, \xi^{\alpha_1}, ..., \xi^{\alpha_m})$,
\end{center}
and such that, for every $k_1, ..., k_n$ individuals in $\mathfrak{B}$, for every ground $g_i$ over $\mathfrak{B}$ for $\alpha_i[k_1, ..., k_n/x_1, ..., x_n]$ ($1 \leq i \leq m$),

\begin{center}
    $f(k_1, ..., k_n, g_1, ..., g_m)$
\end{center}
is a ground over $\mathfrak{B}$ for $\beta[k_1, ..., k_n/x_1, ..., x_n]$. Some of these operations may bind individual variables, with usual restrictions on so-called proper variables. An example of non-primitive first-level operation on grounds - over every base - is an operation of conjunction elimination $f_{\wedge, i}(\xi^{\alpha_1 \wedge \alpha_2})$. Its operational type is $\alpha_1 \wedge \alpha_2 \rhd \alpha_i$ ($i = 1, 2$), since it is specified by the equation

\begin{center}
    $f_{\wedge, i}(\wedge I(g_1, g_2)) = g_i$.
\end{center}

Given $\Gamma_1, ... \Gamma_n, \Delta$ sets of open or closed formulas with $\Delta \subseteq \bigcup_{i \leq n} \Gamma_i$, and given $\alpha_1, ..., \alpha_n, \beta$ open or closed, we can now say that a second-level operation on grounds over $\mathfrak{B}$ of operational type

\begin{center}
    $(\Gamma_1 \rhd \alpha_1, ..., \Gamma_n \rhd \alpha_n) \rhd (\Delta \rhd \beta)$
\end{center}
is a constructive function

\begin{center}
    $f(\xi^{\Gamma_1 \rhd \alpha_1}, ..., \xi^{\Gamma_n \rhd \alpha_n})$
\end{center}
such that, for every $\phi_i$ ground over $\mathfrak{B}$ for $\alpha_i$ if $\Gamma_i = \emptyset$, or first-level operation on grounds over $\mathfrak{B}$ of operational type $\Gamma_i \rhd \alpha_i$ if $\Gamma_i \neq \emptyset$ ($1 \leq i \leq n$),

\begin{center}
    $f(\phi_1, ..., \phi_n)$
\end{center}
is a ground over $\mathfrak{B}$ for $\beta$ if $\Delta = \emptyset$, or a first-level operation on grounds over $\mathfrak{B}$ of operational type $\Delta \rhd \beta$ if $\Delta \neq \emptyset$. Some of these operations may bind typed variables. An example of non-primitive second-level operation on grounds - over every base - is an operation of existential elimination $f_\exists(\xi^{\exists x \alpha(x)}, \xi^{\alpha(x) \rhd \beta})$ binding $x$ and $\xi^{\alpha(x)}$ on the second entry. Its operational type is $[\exists x \alpha(x), (\alpha(x) \rhd \beta) \rhd \beta]$, since it is specified by the equation

\begin{center}
    $f_\exists(\exists I[\alpha(t) \rhd \exists x \alpha(x)](g), f(x, \xi^{\alpha(x)})) = f(t, g)$.\footnote{One may want that an operation that binds one or more individual variables is not of first, but of second level. This would be a way of classifying the ‘‘complexity" of the operations according to whether they bind or not individual or ground-variables. Here, we have pursued a different approach, since we have classified the ‘‘complexity" of the operations according to whether expressions of the form $\Gamma \rhd \alpha$ occur or not in the domain - namely, if the domain contains only formulas, then the operation is of first-level, while if it contains at least an expression of the form $\Gamma \rhd \alpha$, then the operation is of second level. These two approaches essentially boil down to the same, with respect to our concerns here. For further details see [d'Aragona 2021a].}
\end{center}

We can now introduce the clauses for $\rightarrow$ and $\forall$. They involve primitive and meaning-constitutive operations that bind some variables. We indicate the binding by writing the bound variable after the symbol of the operation. So, given $\alpha, \beta$ and $\forall x \alpha(x)$ closed,

\begin{itemize}
    \item[($\rightarrow$)] $g$ is a ground over $\mathfrak{B}$ for $\alpha \rightarrow \beta$ iff $g$ is of the form $\rightarrow I \xi^\alpha(f(\xi^\alpha))$, where $f(\xi^\alpha)$ is an operation on grounds over $\mathfrak{B}$ of operational type $\alpha \rhd \beta$;
    \item[($\forall$)] $g$ is a ground over $\mathfrak{B}$ for $\forall x \alpha(x)$ iff $g$ is of the form $\forall I x(f(x))$, where $f(x)$ is an operation on grounds over $\mathfrak{B}$ of operational type $\alpha(x)$.
\end{itemize}
Observe that the clauses proceed by simultaneous recursion, via the compositional character of the primitive operations. Observe also that the notion of first-level operation on grounds presupposes the notion of ground for sentences, whereas the notion of second-level operation on grounds presupposes the notions of ground for sentences and of first-level operation on grounds [for a more detailed reconstruction of the notions of ground and operation on grounds the reader may refer to d'Aragona 2021a].

An inference from a set of premises $\Gamma$ to a conclusion $\alpha$ can be now described as the application of an operation on grounds for the elements of $\Gamma$. The inference is valid if the operation has operational type $\Gamma \rhd \alpha$. A proof is defined as a finite chain of valid inferences. ToG reverses the order that Prawitz's previous semantics adopted for the intertwinement of the concepts of valid inference and proof. As already said, in these semantics inferential validity was defined by saying that an inference is valid when it yields a valid argument or proof when appended to valid arguments or proofs. The priority given to valid arguments and proofs, and the fact that inferential validity is defined in terms of these notions, impedes to account for the epistemic power of correct deduction. Correct deduction compels to accept its conclusion \emph{because} it is made of valid - and hence compelling - inferences, so that the notion of valid inference must be explanatorily prior to those of valid argument or proof. This is achieved in ToG by defining proofs in terms of valid inferences, and valid inferences in terms of non-inferential (operations on) grounds. 

In spite of this, though, it is doubtful whether ToG fulfills the task of accounting for epistemic compulsion. In fact, ToG shares with Prawitz's earlier approach a problem of recognizability [Prawitz 2015, but see also d'Aragona 2019a, 2021b, Usberti 2015, 2019] and a problem of vacuous validity [Prawitz 2021, but see also d'Aragona 2021b]. In light of these problems, the simple performing of a valid inference may not be enough for the inferential agent to be justified in asserting the conclusion. The mere fact that the agent is in possession of an operation that transforms grounds for the premises into grounds for the conclusion - i.e. that the agent \emph{knows how} to obtain a ground for the conclusion given grounds for the premises - may not imply that the agent also recognizes that the operation has that behavior - i.e. that the agent \emph{knows that} a ground for the conclusion is actually delivered by the applied operation. It is for example clear that this recognition is impossible if we understand the required recognizability as a decision procedure. Because of G\"{o}del's incompleteness theorems, the class of operations on grounds cannot be recursive. Already in the case of a base for first-order arithmetic, no recursive set of operations can account for all the grounds on this base. But even when we limit ourselves to pure first-order logic, from the definitions above it is easy to see that the empty function is an operation on grounds over $\mathfrak{B}$ of operational type $\alpha \rhd \beta$ for any $\beta$, if no ground for $\alpha$ exists on $\mathfrak{B}$. To decide that empty functions of this kind are operations on grounds of the required type, we should have a decision procedure for the problem whether $\alpha$ has or not a ground on $\mathfrak{B}$, and we know that this may be impossible.

A discussion of these problems is however beyond the scope of this paper. Rather, we focus on Prawitz's suggestion of a formal development of ToG in terms of languages of grounding whose terms denote grounds and operations on grounds, and encode a deductive activity through which these objects can be built. These terms may be thus considered as ‘‘programs" that contain instructions for obtaining a ground or an operation on grounds $g$, and whose execution yields $g$. Prawitz specifies that, because of G\"{o}del incompleteness theorems, languages of grounding must be thought of as ‘‘open" since, as seen, no ‘‘closed" language for e.g. first-order arithmetic can contain grounds for undecidable first-order arithmetic formulas [Prawitz 2015].

A ‘‘closed" language can be here understood as a language with a recursive set of operational symbols that denote an equally recursive set of operations on grounds. Thus, an ‘‘open" language can be understood as a ‘‘closed" language to which we can indefinitely add new symbols for new operations. Although the overall language of ToG is hence expected to be non-recursive, the idea of indefinitely extending ‘‘closed" languages suggests that ‘‘open" languages can be dealt with by outlining a hierarchy of ‘‘closed" languages, ordered through a (partial) extension relation. This intuition may be explored in two ways. Either languages of grounding come with denotation functions that associate terms, as well as the operational symbols they are built up of, to grounds or operations on grounds - so, the extension of a ‘‘closed" language will also involve the extension of the denotation function defined over it [this idea is developed in d'Aragona 2021a]. Otherwise, one can define formal systems over these languages, which requires an enrichment of the term-languages to languages with formulas. The hierarchy of languages thereby becomes a hierarchy of systems, and the extension relation among languages and denotation functions can be understood as an extension relation among systems defined over the languages. This latter strategy is the one we pursue in this paper.

\section{Gentzen-languages}

In this section we exemplify a language of grounding with terms and formulas. In the next section, we build a system of grounding over it. The language and the system are particularly suitable for a sharp illustration of the general structure of \emph{any} language and of \emph{any} system of grounding. Later on, this will permit us to generalize the examples to generic languages and generic systems.

The language of grounding we define in this section is a sort ground-theoretic translation of an intuitionistic first-order system in a Gentzen natural deduction style, although we enrich it through formulas that express denotation and equivalent denotation. But first, we have to introduce the preliminary notions of background language and atomic base.

\subsection{Background language and atomic base}

Following Prawitz's suggestions [mainly Prawitz 2015], as well as Prawitz's standard approach in his previous semantics [mainly Prawitz 1973], the languages of grounding must be relative to a background language over an atomic base.

Background languages play a twofold role. On the one hand, they provide types for grounds and operations on grounds, as well as for the alphabet and the terms of the languages of grounding. Through this typing, on the other hand, they indicate for what formula a given object is a ground, as well as for what formula a given term denotes a ground. So, an object of type $\alpha_1, ..., \alpha_n \rhd \beta$ is to be a first-level operation on grounds that transforms grounds for $\alpha_i$ into grounds for $\beta$. The operation is to be denoted in some suitable language of grounding by some suitable term of type $\beta$ that involves free ground-variables $\xi^{\alpha_1}, ..., \xi^{\alpha_n}$ - the example can be easily generalised to simple grounds, or to other kinds of operations on grounds, but for further details see [d'Aragona 2021a].

A \emph{background language} is here a first-order logical language $L$ defined as usual - with an atomic constant $\bot$ for the absurd. We assume as defined in a standard way also:

\begin{itemize}
    \item free variables in terms - noted $FV(t)$ - and free and bound variables in formulas - noted $FV(\alpha)$ and $BV(\alpha)$, with $FV(\Gamma) = \{FV(\alpha) \ | \ \alpha \in \Gamma\}$;
    \item substitution functions of individual variables with terms in terms - noted $t[u/x]$ - and formulas - noted $\alpha[t/x]$;
    \item substitution functions of atomic formulas with formulas in formulas - noted $\alpha[\beta/\gamma]$;
    \item a term being free for a variable in a formula;
    \item a formula being free for an atomic formula in a formula.
\end{itemize}

\noindent We indicate with $\texttt{TERM}_L$ the set of the terms of $L$, with $\texttt{ATOM}_L$ the set of the atomic formulas of $L$, and with $\texttt{FORM}_L$ the set of the formulas of $L$. An expansion of $L$ is, of course, a first-order logical language $L^*$ whose alphabet contains the alphabet of $L$ - we write $L \subseteq L^*$.

Atomic bases play a twofold role too. They serve as inductive basis in definitions about grounds, operations on grounds and languages of grounding, and they provide meaning for the atomic formulas of $L$ through an atomic system over it. The idea is that the rules of an atomic system over $L$ that concern a symbol $s$ in $L$ fix the meaning of $s$ by showing how $s$ behaves in deduction. Both the concept of ‘‘rule concerning a symbol" and that of ‘‘determination of meaning through deductive behavior" may be assumed as primitive - although they are seemingly very complex [a semantics inspired by these concepts is developed in Cozzo 1994].

Following Prawitz, an \emph{atomic system} $\texttt{S}$ \emph{over} $L$ can be defined as a Post-system over $L$, i.e. a pair $\langle L^*, \mathfrak{R} \rangle$ where $L^* \subseteq L$ and $\mathfrak{R}$ is a recursively axiomatisable set of production rules

\begin{prooftree}
\AxiomC{$\alpha_1, ..., \alpha_n$}
\UnaryInfC{$\beta$}
\end{prooftree}
such that:

\begin{itemize}
    \item $\alpha_i \in \texttt{ATOM}_{L^*}$ and $\alpha_i \neq \bot$ ($1 \leq i \leq n$);
    \item  $\beta \in \texttt{ATOM}_{L^*}$ or $\beta = \bot$ and, if $x \in FV(\beta)$, then there is $1 \leq i \leq n$ such that $x \in FV(\alpha_i)$.
\end{itemize}
We indicate with $\texttt{DER}_{\texttt{S}}$ the set of the derivations of $\texttt{S}$ - defined in an usual inductive way. An expansion of $\texttt{S}$ is an atomic system $\texttt{S}^* = \langle L^{**}, \mathfrak{R}^* \rangle$ with $L^{*} = L^{**}$ and $\mathfrak{R} \subseteq \mathfrak{R}^{*}$, or $L^{*} \subset L^{**}$ and $\mathfrak{R} \subset \mathfrak{R}^{*}$ - we write $\texttt{S} \subseteq \texttt{S}^*$. An \emph{atomic base} $\mathfrak{B}$ \emph{over} $L$ is a triple $\langle \texttt{C}, \texttt{R}, \texttt{S} \rangle$ where $\texttt{C}$ is the set of the individual constants of $L$, $\texttt{R}$ is the set of the relational constants of $L$, and $\texttt{S}$ is an atomic system over $L$. An expansion of $\mathfrak{B}$ is an atomic base $\mathfrak{B}^*$ over $L^{*}$ with atomic system $\texttt{S}^*$ such that $L \subseteq L^{*}$ and $\texttt{S} \subseteq \texttt{S}^*$ - we write $\mathfrak{B} \subseteq \mathfrak{B}^*$.

\subsection{Basic Gentzen-language of grounding}

\begin{defn}
Given a background language $L$ and an atomic base $\mathfrak{B}$ over it, a \emph{Gentzen-language of grounding} over $\mathfrak{B}$ - indicated with $\texttt{Gen}$ - is determined by an \emph{alphabet} $\texttt{Al}_{\texttt{Gen}}$ and a set of \emph{typed terms} $\texttt{TERM}_{\texttt{Gen}}$. $\texttt{Al}_{\texttt{Gen}}$ contains the following elements:

\begin{itemize}
    \item an individual constant $\delta_i$ naming the $i$-th $\Delta \in \texttt{DER}_{\texttt{S}}$ with no undischarged assumptions and unbound variables ($i \in \mathbb{N}$)\footnote{If binding and dischargements are not allowed in the rules of $\texttt{S}$, then $\Delta$ has to be a derivation whose leaves are axioms, and where no individual variable occurs. In what follows, however, we shall abstract from whether the rules of $\texttt{S}$ may or not bind variables or discharge assumptions, although this question becomes central when dealing with other aspects of Prawitz's semantics - e.g. completeness [Piecha, de Campos Sanz \& Schroeder-Heister 2015, Piecha \& Schroeder-Heister 2018].};
    \item typed-variables $\xi^\alpha_i$ ($\alpha \in \texttt{FORM}_L$, $i \in \mathbb{N}$);
    \item operational symbols of the form $F[\omega]$ - where $\omega$ is meant to be the type of the operational symbol - possibly binding individual and typed-variables. We call primitive the following symbols:
    \begin{itemize}
        \item $\wedge I[\alpha, \beta \rhd \alpha \wedge \beta]$ ($\alpha, \beta \in \texttt{FORM}_L$)
        \item $\vee I[\alpha_i \rhd \alpha_1 \vee \alpha_2]$ ($\alpha_i \in \texttt{FORM}_L$, $i = 1, 2$);
        \item $\rightarrow I[\beta \rhd \alpha \rightarrow \beta]$ binding $\xi^\alpha_i$ in the term of type $\beta$ to which the symbol is applied ($\alpha, \beta \in \texttt{FORM}_L$, $i \in \mathbb{N}$);
        \item $\forall I[\alpha(x_i) \rhd \forall y_j(\alpha(y_j/x_i)]$ binding $x_i$ ($i, j \in \mathbb{N}$, $\alpha \in \texttt{FORM}_L$);
        \item $\exists I[\alpha(t) \rhd \exists x_i \alpha(x_i)]$ ($t \in \texttt{TERM}_L$, $\alpha \in \texttt{FORM}_L$, $i \in \mathbb{N}$)
    \end{itemize}
    We call non-primitive the following symbols:
    \begin{itemize}
        \item $\wedge_{E, i}[\alpha_1 \wedge \alpha_2 \rhd \alpha_i]$ ($\alpha_i \in \texttt{FORM}_L$, $i = 1, 2$);
        \item $\vee E[\alpha_1 \vee \alpha_2, \gamma, \gamma \rhd \gamma]$ binding $\xi^{\alpha_i}$ in the $i + 1$-th entry ($\alpha_i, \gamma \in \texttt{FORM}_L$, $i = 1, 2$);
        \item $\rightarrow E[\alpha \rightarrow \beta, \alpha \rhd \beta]$ ($\alpha, \beta \in \texttt{FORM}_L$);
        \item $\forall E[\forall x_i \alpha(x_i) \rhd \alpha(t/x_i)]$ ($\alpha \in \texttt{FORM}_L$, $t \in \texttt{TERM}_L$, $i \in \mathbb{N}$);
        \item $\exists E[\exists x_i \alpha(x_i), \beta \rhd \beta]$ binding $y_j$ and $\xi^{\alpha(y_j)}_k$ on the second entry in the term of type $\beta$ to which the symbol is applied ($\alpha, \beta \in \texttt{FORM}_L$, $i, j, k \in \mathbb{N}$);
        \item $\bot_\alpha[\bot \rhd \alpha]$ ($\alpha \in \texttt{FORM}_L$).
    \end{itemize}
\end{itemize}
We omit operational types, subscripts and superscripts whenever possible. $\texttt{TERM}_{\texttt{Gen}}$ is the smallest set $X$ such that - we indicate typing by colon, and binding of variables by writing the bound variables immediately after the operational symbol:

\begin{itemize}
    \item $\delta : \alpha \in X$ where $\alpha$ is the conclusion of the $\Delta \in \texttt{DER}_S$ named by $\delta$;
    \item $\xi^\alpha_i : \alpha \in X$;
    \item $T : \alpha, U : \beta \Rightarrow \wedge I(T, U) : \alpha \wedge \beta \in X$;
    \item $T : \alpha_1 \wedge \alpha_2 \in X \Rightarrow \wedge_{E, i}(T): \alpha_i \in X$;
    \item $T : \alpha_i \in X \Rightarrow \vee I[\alpha_i \rhd \alpha_1 \vee \alpha_2](T) : \alpha_1 \vee \alpha_2 \in X$;
    \item $T : \alpha \vee \beta, U : \gamma, Z : \gamma \in X \Rightarrow \vee E \ \xi^\alpha \ \xi^\beta(T, U, Z) : \gamma \in X$;
    \item $T : \beta \in X \Rightarrow \ \rightarrow I \xi^\alpha(T) : \alpha \rightarrow \beta \in X$;
    \item $T : \alpha \rightarrow \beta, U : \alpha \in X \Rightarrow \ \rightarrow E(T, U) : \beta \in X$;
    \item $U : \alpha(x) \in X \Rightarrow \forall I[\alpha(x) \rhd \forall y \alpha(y)] x (U) : \forall y \alpha(y) \in X$;
    \item $T : \forall x \alpha(x) \in X \Rightarrow \forall E[\forall x \alpha(x) \rhd \alpha(t)](T) : \alpha(t) \in X$;
    \item $T : \alpha(t) \in X \Rightarrow \exists I[\alpha(t) \rhd \exists x \alpha(x)](T) : \exists x \alpha(x) \in X$;
    \item $T : \exists x \alpha(x), U : \beta \in X \Rightarrow \exists E y \ \xi^{\alpha(y)}(T, U) : \beta \in X$;
    \item $T : \bot \in X \Rightarrow \bot_\alpha(T) : \alpha \in X$.
\end{itemize}
In the clause for $\forall I$, it must hold that $x \notin FV(\beta)$ for $\xi^\beta \in FV^T(U)$. In the clause for $\exists E$, it must hold that $y \notin FV(\gamma)$ for $\gamma \neq \alpha(y)$ and $\xi^\gamma \in FV^T(U)$. In the clause for $\forall E$ and in that for $\exists I$, $t$ must be free for $x$ in $\alpha(x)$. For the notation $FV^T$, see below.
\end{defn}

\noindent We say that $T$ is \emph{canonical} iff it is an individual constant, or it is a typed-variable or its outermost operational symbol is primitive. $T$ is \emph{non-canonical} otherwise. Observe that, for every logical constant $k$, $k I$ corresponds to the Gentzen introduction rule for $k$, while $k E$ corresponds to the Gentzen elimination rule for $k$.

The following notions can be defined in a standard way: the set $S(T)$ of the \emph{sub-terms} of $T$; the sets $FV^I(T)$ and $BV^I(T)$ of, respectively, the \emph{free} and \emph{bound individual variables} of $T$; the sets $FV^T(U)$ and $BV^T(U)$ of, respectively, the \emph{free} and \emph{bound typed-variables} of $U$. $U$ is \emph{closed} iff $FV^I(U) = FV^T(U) = \emptyset$, \emph{open} otherwise.

A \emph{substitution of} $x$ \emph{with} $t$ \emph{in} $T$ and a \emph{substitution of} $\xi^\alpha$ \emph{with} $W : \alpha$ \emph{in} $T$ are functions $\texttt{TERM}_{\texttt{Gen}} \to \texttt{TERM}_{\texttt{Gen}}$ defined in a usual inductive way. The inductive basis is respectively

\begin{center}
    $\delta[t/x] = \delta$ and $\xi^\alpha[t/x] = \xi^{\alpha[t/x]}$
\end{center}

\begin{center}
    $\delta[W/\xi^\alpha] = \delta$ and $\xi^\beta[W/\xi^\alpha] = \begin{cases} \xi^\beta & \text{if} \ \alpha \neq \beta\\ W & \text{if} \ \alpha = \beta \end{cases}$
\end{center}
but one has also to deal with possible clashes of bound variables, and of terms required to be free for individual variables in formulas. This can be done following the principles that individual variables are replaced in the types of bound typed-variables, that bound variables are never replaced, that no variable becomes bound after substitution, and finally that substitution respects the restrictions on terms free for individual variables in formulas. So for example

\begin{center}
    $\rightarrow I \xi^\alpha(U)[t/x] = \ \rightarrow I \xi^{\alpha[t/x]}(U[t/x])$
\end{center}
\medskip

\begin{center}
    $\forall I y(U)[t/x] = \begin{cases} \forall I y(U) & \text{if} \ x = y \\ \forall I y(U[t/x]) & \text{if} \ x \neq y, y \notin FV(t) \\ \forall I z((U[z/y])[t/x]) & \text{if} \ x \neq y, y \in FV(t), \ \text{for} \ z \notin FV(t), z \notin FV^I(U) \end{cases}$
\end{center}
\medskip

\begin{center}
    $\exists I(U)[t/x] = \begin{cases} \exists I(U) & \text{if} \ U : \alpha(s/y), \exists I(U) : \exists y \alpha(y), x \in FV(s), t \ \text{not free for} \ y \ \text{in} \ \alpha(y) \\ \exists I(U[t/x]) & \text{otherwise} \end{cases}$
\end{center}
\medskip

\begin{center}
    $\vee E \ \xi^\beta_i \ \xi^\gamma_j (U, Z_1, Z_2)[W/\xi^\alpha_h] =$
\end{center}

\begin{itemize}
    \item $\vee E \ \xi^\beta_i \ \xi^\gamma_j(U[W/\xi^\alpha_h], Z_1, Z_2)$\\if $\xi^\alpha_h = \xi^\beta_i = \xi^\gamma_j$;
    \item $\vee E \ \xi^\beta_i \ \xi^\gamma_j(U[W/\xi^\alpha_h], Z_1[W/\xi^\alpha_h], Z_2)$\\if $\xi^\alpha_h \neq \xi^\beta_i$, $\xi^\beta_i \notin FV^T(W)$, $\xi^\alpha_h = \xi^\gamma_j$;
    \item $\vee E \ \xi^\beta_s \ \xi^\gamma_j(U[W/\xi^\alpha_h], (Z_1[\xi^\beta_s/\xi^\beta_i])[W/\xi^\alpha_h], Z_2)$\\if $\xi^\alpha_h \neq \xi^\beta_i$, $\xi^\beta_i \in FV^T(W)$, $\xi^\alpha_h = \xi^\gamma_j$, for $\xi^\beta_s \notin FV^T(Z_1)$, $\xi^\beta_s \notin FV^T(W)$;
    \item $\vee E \ \xi^\beta_i \ \xi^\gamma_j(U[W/\xi^\alpha_h], Z_1, Z_2[W/\xi^\alpha_h])$\\if $\xi^\alpha_h = \xi^\beta_i$, $\xi^\alpha_h \neq \xi^\gamma_j$, $\xi^\gamma_j \notin FV^T(W)$;
    \item $\vee E \ \xi^\beta_i \ \xi^\gamma_s(U[W/\xi^\alpha_h], Z_1, (Z_2[\xi^\gamma_s/\xi^\gamma_j])[W/\xi^\alpha_h])$\\if $\xi^\alpha_h = \xi^\beta_i$, $\xi^\alpha_h \neq \xi^\gamma_j$, $\xi^\gamma_j \in FV^T(W)$, for $\xi^\gamma_s \notin FV^T(Z_2)$, $\xi^\gamma_s \notin FV^T(W)$;
    \item $\vee E \ \xi^\beta_i \ \xi^\gamma_j(U[W/\xi^\alpha_h], Z_1[W/\xi^\alpha_h], Z_2[W/\xi^\alpha_h])$\\if $\xi^\alpha_h \neq \xi^\beta_i$, $\xi^\beta_i \notin FV^T(W)$, $\xi^\alpha_h \neq \xi^\gamma_j$, $\xi^\gamma_j \notin FV^T(W)$;
    \item $\vee E \ \xi^\beta_s \ \xi^\gamma_j(U[W/\xi^\alpha_h], (Z_1[\xi^\beta_s/\xi^\beta_i])[W/\xi^\alpha_h], Z_2[W/\xi^\alpha_h])$\\if $\xi^\alpha_h \neq \xi^\beta_i$, $\xi^\beta_i \in FV^T(W)$, $\xi^\alpha_h \neq \xi^\gamma_j$, $\xi^\gamma_j \notin FV^T(W)$, for $\xi^\beta_s \notin FV^T(Z_1)$, $\xi^\beta_s \notin FV^T(W)$;
    \item $\vee E \ \xi^\beta_i \ \xi^\gamma_s(U[W/\xi^\alpha_h], Z_1[W/\xi^\alpha_h], (Z_2[\xi^\gamma_s/\xi^\gamma_j])[W/\xi^\alpha_h])$\\if $\xi^\alpha_h \neq \xi^\beta_i$, $\xi^\beta_i \notin FV^T(W)$, $\xi^\alpha_h \neq \xi^\gamma_j$, $\xi^\gamma_j \in FV^T(W)$, for $\xi^\gamma_s \notin FV^T(Z_2)$, $\xi^\gamma_s \notin FV^T(W)$;
    \item $\vee E \ \xi^\beta_s \ \xi^\gamma_t(U[W/\xi^\alpha_h], (Z_1[\xi^\beta_s/\xi^\beta_i])[W/\xi^\alpha_h], (Z_2[\xi^\gamma_t/\xi^\gamma_j])[W/\xi^\alpha_h])$\\if $\xi^\alpha_h \neq \xi^\beta_i$, $\xi^\beta_i \in FV^T(W)$, $\xi^\alpha_h \neq \xi^\gamma_j$, $\xi^\gamma_j \in FV^T(W)$, for $\xi^\beta_s \notin FV^T(Z_1)$, $\xi^\beta_s \notin FV^T(W)$, $\xi^\gamma_t \notin FV^T(Z_2)$, $\xi^\gamma_t \notin FV^T(W)$.
\end{itemize}

\subsection{An enriched Gentzen-language}

The terms of $\texttt{Gen}$ are suitable for denoting grounds and operations on grounds. For example, every constant $\delta : \alpha$ may be understood as standing for the atomic derivation of which it is a name, and thus for a ground for $\alpha$. Likewise, the operation $\wedge I[\alpha, \beta \rhd \alpha \wedge \beta]$ and the term $\wedge I (\xi^\alpha, \xi^\beta)$ can be understood as standing for an operation that takes grounds for $\alpha$ and grounds for $\beta$ returning a ground for $\alpha \wedge \beta$, and thus for an operation on grounds of operational type $\alpha, \beta \rhd \alpha \wedge \beta$. Observe that, in the latter case, we have an operational symbol and a term that stand for the same object, but this does not generally occur; as soon as the operational symbol binds some variables - e.g. $\rightarrow I [\beta \rhd \alpha \rightarrow \beta]$ - the object it stands for \emph{cannot} be the same as the one that some term stands for - say, $\rightarrow I \xi^\alpha (T)$, where $T : \beta$ stands for an operation on grounds of operational type $\alpha \rhd \beta$. This is because we have not allowed typed-variables of the form $\xi^{\Gamma \rhd \alpha}$, so while the operational symbols may stand for operations on grounds that take \emph{other operations on grounds} as values, a term will stand for either a ground, or an operation on grounds that takes simple grounds as values. If one stays within a denotational approach [like the one partially anticipated in d'Aragona 2018 in terms of the notion of \emph{operational symbol fitting with an operational type}, and developed in d'Aragona 2021a], this difficulty can be overcome by first distinguishing carefully the denotation of the operational symbols from that of the terms, and then by defining compositionally the latter in terms of the former. In the systems we are going to develop, the limitation is instead overcome by enriching the set of the terms through typed functional variables ranging over terms that are to stand for operations on grounds of a specific operational type, and then by quantifying over these variables.

However, the fact that $\texttt{Gen}$ only contains terms may be considered as a shortcoming. Without formulas, we cannot express \emph{within it} properties of terms and of other components of the alphabet of $\texttt{Gen}$ like those indicated above. We are obliged to move to a meta-level, and the proposed formalisation of ToG would come out partial, if not defective. It is for this reason that we now enrich $\texttt{Gen}$ towards a more powerful language, with formulas that concern two aspects: the fact that a given term denotes a ground for a given formula, or an operation on grounds with a given domain and a given co-domain, and the fact that two terms denote (extensionally) equal (operations on) grounds. This is achieved by adding to $\texttt{Gen}$ two binary predicates - one for speaking of denotation, and another for speaking of equality - and the first-order logical constants.

To avoid confusion between formulas of $L$ and formulas of the enriched $\texttt{Gen}$, we shall indicate the latter with the first capital letters of the Latin alphabet $A, B, C, ...$. To avoid confusion between the logical constants of $L$ and the logical constants of the enriched $\texttt{Gen}$, we shall indicate the latter as follows: conjunction is $\times$, disjunction is $+$, implication is $\supset$, universal quantification is $\Pi$ and existential quantification is $\mathfrak{E}$. Finally, to distinguish the absurd constant of $L$ from that of the enriched $\texttt{Gen}$, we indicate the latter with $\bot^G$.

\begin{defn}
Given $\texttt{Gen}$ over $\mathfrak{B}$, an \emph{enriched Gentzen-language over} $\mathfrak{B}$ - indicated with $\texttt{Gen}^\varepsilon$ - is determined by an \emph{alphabet} $\texttt{Al}_{\texttt{Gen}^\varepsilon}$, a set of \emph{typed terms} $\texttt{TERM}_{\texttt{Gen}^\varepsilon}$ and a set of \emph{formulas} $\texttt{FORM}_{\texttt{Gen}^\varepsilon}$ - whose set of \emph{atomic formulas} is indicated with $\texttt{ATOM}_{\texttt{Gen}^\varepsilon}$. $\texttt{Al}_{\texttt{Gen}^\varepsilon}$ is such that $\texttt{Al}_{\texttt{Gen}} \subset \texttt{Al}_{\texttt{Gen}^\varepsilon}$ and moreover it contains

\begin{itemize}
    \item functional variables $\texttt{h}^\alpha_{x, i}$, $\texttt{f}^\alpha_j$ ($\alpha \in \texttt{FORM}_L$, $i, j \in \mathbb{N}$) - whose role and ‘‘meaning" is explained below;
    \item binary relation symbols $Gr$ and $\equiv$;
    \item logical constants $\times$, $+$, $\supset$, $\Pi$, $\mathfrak{E}$ and $\bot^G$.
\end{itemize}
Whenever possible, we will omit subscripts and superscripts. $\texttt{TERM}_{\texttt{Gen}^\varepsilon}$ is the smallest set $X$ such that $\texttt{TERM}_{\texttt{Gen}} \subset X$ and

\begin{itemize}
    \item $\texttt{h}^\alpha_x(t) : \alpha(t/x) \in X$ ($t \in \texttt{TERM}_L$);
    \item $\texttt{f}^\alpha(T) : \alpha \in X$ ($T \in \texttt{TERM}_{\texttt{Gen}}$).
\end{itemize}
plus recursive clauses for forming complex terms through the operational symbols of $\texttt{Gen}$ as indicated in definition 1. $\texttt{FORM}_{\texttt{Gen}^\varepsilon}$ is the smallest set $X$ such that:

\begin{itemize}
    \item $T : \alpha \in \texttt{TERM}_{\texttt{Gen}^\varepsilon} \Rightarrow Gr(T, \alpha) \in X$;
    \item $T : \alpha \in \texttt{TERM}_{\texttt{Gen}^\varepsilon}, U : \alpha \in \texttt{TERM}_{\texttt{Gen}^\varepsilon} \Rightarrow T \equiv U \in X$;
    \item $\bot^G \in X$;
    \item $A, B \in X \Rightarrow A \square B \in X$ (for $\square = \times, +, \supset$)
    \item $A \in X \Rightarrow \square \ \nu \ A \in X$ (for $\square = \Pi, \mathfrak{E}$ and $\nu = x_i, \xi^\alpha_i$ for $i \in \mathbb{N}$ and $\alpha \in \texttt{FORM}_L$).
\end{itemize}
\end{defn}
\noindent We put as usual $\neg^G A \stackrel{def}{=} A \supset \bot^G$. Observe that for an atomic formula $Gr(T, \alpha)$ to be well-formed, we require that $T$ has type $\alpha$ in $\texttt{Gen}^\varepsilon$, and that for an atomic formula $T \equiv U$ to be well-formed we require that $T$ and $U$ have the same type in $\texttt{Gen}^\varepsilon$. These restrictions are due to the fact that, if we did not adopt them, we would have some unwanted consequences in the formal system of grounding of the next section.

Formulas of the kind $Gr(T, \alpha)$ are meant to express that $T$ is a ground for $\alpha$, or an operation on grounds with co-domain $\alpha$ when $T$ contains individual or typed unbound variables. In order not to have an excessively heavy notation, we shall from now on write atomic formulas of this kind simply as $T : \alpha$, but we remark that this notation must not be confused with the \emph{meta-linguistic} expression of equal form that has been used so far to indicate that $T$ is a term of type $\alpha$ in $\texttt{Gen}$. Formulas of the kind $T \equiv U$ are meant to express that $T$ and $U$ denote equivalent grounds or operations on grounds - basically, equivalence amounts to \emph{extensional} identity.

It may now be worth saying something about the multi-layered linguistic structure that we are dealing with here. First of all, observe that this structure is made of three levels: (a) a background language $L$, which provides types for the terms of (b) the language $\texttt{Gen}$, which contains terms whose properties are expressed by the formulas of (c) the language $\texttt{Gen}^\varepsilon$. These three levels are meant to refer to two domains: (1) the domain of the objects that the terms of $L$ stand for and (2) the domain of grounds and operations on grounds, named by terms of $\texttt{Gen}$, that justify assertions about (1). The individual variables of $\texttt{Gen}$ (and hence of $\texttt{Gen}^\varepsilon$) range over the elements of (1), while the typed-variables of $\texttt{Gen}$ (and hence of $\texttt{Gen}^\varepsilon$) range over some of the (terms standing for) elements of (2). Additionally, in $\texttt{Gen}^\varepsilon$ we have typed functional variables. A functional variable of the form $\texttt{h}^\alpha_x$ is meant to range over terms that are to stand for operations on grounds of operational type $\Gamma \rhd \alpha$ involving $x$, for $x$ not occurring free in any element of $\Gamma$ - so $\texttt{h}^\alpha_x(t)$ will be any such operation applied to the individual that $t$ stands for. Instead, a functional variable of the form $\texttt{f}^\alpha$ is meant to range over terms that are to stand for operations on grounds of operational type $\Gamma \rhd \alpha$ - so $\texttt{f}^\alpha(Z)$ will be any such operation applied to the ground or operation on grounds that $Z$ stands for. Typed functional variables are thus meant to account for (properties of) operations that bind individual and of typed-variables. For example, if we had not functional variables, we could not express the fact that a ground for $\alpha \rightarrow \beta$ is obtained by applying $\rightarrow I$ to an operation that, when applied \emph{to any} ground for $\alpha$, yields a ground for $\beta$.

Observe that individual variables occur in a term $T$ of $\texttt{Gen}$ only when they occur in the type of the typed-variables of $T$, or in the type of typed-variables used in $T$ to indicate bindings - like e.g. $\rightarrow I \xi^{\alpha(x)} (U)$. In fact, this is not surprising, for a term $T$ of $\texttt{Gen}$ does not \emph{directly} concern the individuals that the terms of $L$ stand for. Rather, $T$ is to stand for a ground that justifies \emph{an assertion} $\alpha$ about these individuals. So, it may concern these individuals only \emph{indirectly}, to the extent that $\alpha$ contains variables ranging over them. Much more importantly, we remark that the quantifiers of $\texttt{Gen}^\varepsilon$ can bind indifferently \emph{all kinds} of variables, be they individual, typed- or functional. This is because, when we express properties of (terms standing for) grounds and operations on grounds, we may need to quantify both over individuals denoted by terms of $L$, and over (terms standing for) grounds and operations on grounds denoted by terms of $\texttt{Gen}$ and $\texttt{Gen}^\varepsilon$. To give just three examples, we must be able to express the fact that:

\begin{itemize}
    \item an operation on grounds of operational type $\alpha(x)$, \emph{whenever} applied to an individual $k$, produces a ground for $\alpha(k)$ (and that, given a term that denotes this operation, \emph{whenever} we replace $x$ with a term denoting $k$, we obtain a term that denotes a ground for $\alpha(k)$);
    \item an operation on grounds of operational type $\alpha \rhd \beta$, \emph{whenever} applied to a ground $g$ for $\alpha$, produces a ground for $\beta$ (and that, given a term that denotes this operation, \emph{whenever} we replace a typed-variable $\xi^\alpha$ with a term denoting a ground for $\alpha$, we obtain a term that denotes a ground for $\beta$);
    \item an operation on grounds of operational type $(\alpha \rhd \beta) \rhd \gamma$, \emph{whenever} applied to an operation on grounds of operational type $\alpha \rhd \beta$, produces a ground for $\gamma$ (and that, given a term that denotes this operation, \emph{whenever} we replace a functional variable $\texttt{f}^\beta$ applied to a typed-variable $\xi^\alpha$ - i.e. $\texttt{f}^\beta(\xi^\alpha)$ - with a term denoting an operation on grounds of operational type $\alpha \rhd \beta$, we obtain a term that denotes a ground for $\gamma$).
\end{itemize}

With this established, we can provide some examples of reading of formulas of $\texttt{Gen}^\varepsilon$ - incidentally, all the formulas we exemplify here are provable in the system we develop in the next section. The formula

\begin{center}
    $\Pi \xi^{\alpha \vee \beta} (\xi^{\alpha \vee \beta} : \alpha \vee \beta \Leftrightarrow \mathfrak{E} \xi^\alpha \mathfrak{\xi^\beta}((\xi^{\alpha \vee \beta} \equiv \vee I (\xi^\alpha) \times \xi^\alpha : \alpha) + (\xi^{\alpha \vee \beta} \equiv \vee I (\xi^\beta) \times \xi^\beta : \beta)))$
\end{center}
is to be read as follows - $\Leftrightarrow$ is the usual abbreviation of $(... \supset ---) \times (--- \supset ...)$ for bi-implication:

\begin{quote}
   for every term $\xi^{\alpha \vee \beta}$ (i.e. for every term of type $\alpha \vee \beta$), $\xi^{\alpha \vee \beta}$ denotes a ground for $\alpha \vee \beta$ iff for some term $\xi^\alpha$ (i.e. for some term of type $\alpha$), for some term $\xi^\beta$ (i.e. for some term of type $\beta$), $\xi^{\alpha \vee \beta}$ is equivalent to $\vee I(\xi^\alpha)$ and $\xi^\alpha$ denotes a ground for $\alpha$, or $\xi^{\alpha \vee \beta}$ is equivalent to $\vee I(\xi^\beta)$ and $\xi^\beta$ denotes a ground for $\beta$. 
\end{quote}
Thus, the formula says that every (term denoting a) ground for $\alpha \vee \beta$ (is equivalent to a term that) begins with the primitive operation $\vee I$ applied either to a (term denoting a) ground for $\alpha$ or to a (term denoting a) ground for $\beta$. The formula

\begin{center}
    $\Pi \xi^{\forall x \alpha(x)} (\xi^{\forall x \alpha(x)} : \forall x \alpha(x) \Leftrightarrow \mathfrak{E} \texttt{h}^{\alpha(x)} (\xi^{\forall x \alpha(x)} \equiv \forall I x (\texttt{h}^{\alpha(x)}(x)) \times \Pi x (\texttt{h}^{\alpha(x)}(x) : \alpha (x))))$
\end{center}
is to be read as follows - $\Leftrightarrow$ is again the usual abbreviation of $(... \supset ---) \times (--- \supset ...)$ for bi-implication:

\begin{quote}
    for every $\xi^{\forall x \alpha(x)}$ (i.e. for every term of type $\forall x \alpha(x)$), $\xi^{\forall x \alpha(x)}$ denotes a ground for $\forall x \alpha(x)$ iff there is term $\texttt{h}^{\alpha(x)}$ where $x$ does not occur free in free typed-variables (i.e. there is term of type $\alpha(x)$ where $x$ does not occur free in free typed-variables) such that $\xi^{\forall x \alpha(x)}$ is equivalent to $\forall I x (\texttt{h}^{\alpha(x)}(x))$ and, for every $x$, $\texttt{h}^{\alpha(x)}(x)$ denotes a ground for $\alpha(x)$.
\end{quote}
Thus, the formula says that every (term denoting a) ground for $\forall x \alpha(x)$ (is equivalent to a term that) begins with the primitive operation $\forall I$ applied to a (term denoting an) operation on grounds of operational type $\alpha(x)$ - with $x$ bound. Finally, the formula

\begin{center}
    $\Pi \xi^{\exists x \alpha(x)} \Pi \texttt{f}^\beta (\xi^{\exists x \alpha(x)} : \exists x \alpha(x) \times \Pi x \Pi \xi^{\alpha(x)} (\xi^{\alpha(x)} : \alpha(x) \supset \texttt{f}^\beta(\xi^{\alpha(x)}) : \beta) \supset \exists E \ x \ \xi^{\alpha(x)} (\xi^{\exists x \alpha(x)}, \texttt{f}^\beta(\xi^{\alpha(x)})) : \beta)$
\end{center}
is to be read as follows:

\begin{quote}
    for every term $\xi^{\exists x \alpha(x)}$ (i.e. for every term of type $\exists x \alpha(x)$), for every term $\texttt{f}^\beta$ possibly involving free typed-variables (i.e. for every term of type $\beta$ possibly involving free typed-variables), if $\xi^{\exists x \alpha(x)}$ denotes a ground for $\exists x \alpha(x)$ and if, for every $x$, for every term $\xi^{\alpha(x)}$ (i.e. for every term of type $\alpha(x)$), if $\xi^{\alpha(x)}$ denotes a ground for $\alpha(x)$ then $\texttt{f}^\beta(\xi^{\alpha(x)})$ denotes a ground for $\beta$, then $\exists E \ x \ \xi^{\alpha(x)} (\xi^{\exists x \alpha(x)}, \texttt{f}^\beta(\xi^{\alpha(x)}))$ denotes a ground for $\beta$.
\end{quote}
Thus the formula says that the operation $\exists E$ yields a (term denoting a) ground for $\beta$ whenever applied to a (term denoting a) ground for $\exists x \alpha(x)$ and to a (term denoting an) operation on grounds of operational type $\alpha(x) \rhd \beta$ - binding $x$ and $\xi^{\alpha(x)}$ on the second entry.

We conclude this section with some standard notions concerning free and bound variables in terms and formulas, as well as substitution of variables with terms in formulas. The set $S(T)$ of the \emph{sub-terms} of $T$ is defined in a usual inductive way in the ‘‘traditional" cases. For functional variables we have

\begin{center}
    $S(\texttt{h}^\alpha(t)) = \{\texttt{h}^\alpha(t)\}$ and $S(\texttt{f}^\alpha(T)) = S(T) \cup \{\texttt{f}^\alpha(T)\}$.
\end{center}
The sets $FV^I(U)$ and $BV^I(U)$ of respectively the \emph{free} and \emph{bound individual variables} of $U$ are defined in a usual inductive way in the ‘‘traditional" cases. For the functional variables we have

\begin{center}
    $\star^I(\texttt{h}^\alpha_x(t)) = \star(t) \cup \star(\alpha(t/x))$ and $\star^I(\texttt{f}^\alpha(T)) = \star^I(T) \cup \star(\alpha)$
\end{center}
where $\star = FV, BV$. The sets $FV^T(U)$ and $BV^T(U)$ of respectively the \emph{free} and \emph{bound typed-variables} of $U$ are defined as usual in an inductive way in the ‘‘traditional" cases. For the functional variables we have

\begin{center}
    $\star^T(\texttt{h}^\alpha(t)) = \emptyset$ and $\star^T(\texttt{f}^\alpha(U)) = \star^T(U)$
\end{center}
where $\star = FV, BV$. Finally, the sets $FV^F(U)$ and $BV^F(U)$ of respectively the \emph{free} and \emph{bound functional variables} of $U$ are defined inductively starting from the base

\begin{center}
    $\star^F(\texttt{h}^\alpha(t)) = \{\texttt{h}^\alpha\}$ and $\star^F(\texttt{f}^\alpha(T)) = \{\texttt{f}^\alpha\}$
\end{center}
where $\star = FV, BV$. A term $U$ is \emph{closed} iff $FV^I(U) = FV^T(U) = FV^F(U) = \emptyset$, \emph{open} otherwise.

A \emph{substitution of} $x$ \emph{with} $u$ in $T$ is a function $\texttt{TERM}_{\texttt{Gen}^\varepsilon} \to \texttt{TERM}_{\texttt{Gen}^\varepsilon}$ defined in a usual inductive way in the ‘‘traditional cases". For the functional variables we have

\begin{center}
    $\texttt{h}^\alpha(t)[u/x] = \texttt{h}^{\alpha[u/x]}(t[u/x])$ and $\texttt{f}^\alpha(T)[u/x] = \texttt{f}^{\alpha[u/x]}(T[u/x])$.
\end{center}
A \emph{substitution of} $\xi^\alpha$ \emph{with} $U : \alpha$ in $T$ is a function $\texttt{TERM}_{\texttt{Gen}^\varepsilon} \to \texttt{TERM}_{\texttt{Gen}^\varepsilon}$ defined in a usual inductive way for the ‘‘traditional cases". For the functional variables we have

\begin{center}
    $\texttt{h}^\alpha(t)[U/\xi^\alpha] = \texttt{h}^\alpha(t)$ and $\texttt{f}^\alpha(Z)[U/\xi^\alpha] = \texttt{f}^\alpha(Z[U/\xi^\alpha])$.
\end{center}
A \emph{substitution of} $h^\beta_x$ \emph{with} $U : \beta$ in $Z$ is a function $\texttt{TERM}_{\texttt{Gen}^\varepsilon} \to \texttt{TERM}_{\texttt{Gen}^\varepsilon}$ that is defined only when $x \notin FV(\alpha)$ for $\xi^\alpha \in FV^T(U)$. It runs in a usual inductive way in the ‘‘traditional cases", while for the functional variables we have

\begin{center}
    $\texttt{h}^\beta_x(t)[U/\texttt{h}^\beta_x] = U[t/x]$ and $\texttt{f}^\alpha(T)[U/h^\beta_x] = \texttt{f}^\alpha(T)$
\end{center}
A \emph{substitution of} $\texttt{f}^\alpha$ \emph{with} $U : \alpha$ in $T$ is an always defined function $\texttt{TERM}_{\texttt{Gen}^\varepsilon} \to \texttt{TERM}_{\texttt{Gen}^\varepsilon}$. It runs in a usual inductive way in the ‘‘traditional cases", while for the functional variables we have

\begin{center}
    $\texttt{h}^\alpha(t)[U/\texttt{f}^\alpha] = \texttt{h}^\alpha(t)$ and $\texttt{f}^\alpha(Z)[U/\texttt{f}^\alpha] = U$
\end{center}
In the last three cases, we have the following restriction: the substitution is defined only when $U \in \texttt{TERM}_{\texttt{Gen}}$. This is to avoid cases like

\begin{center}
    $\texttt{f}^\alpha(\xi^\alpha)[\texttt{f}^\alpha(\xi^\alpha)/\xi^\alpha] = \texttt{f}^\alpha(\xi^\alpha[\texttt{f}^\alpha(\xi^\alpha)/\xi^\alpha]) = \texttt{f}^\alpha(\texttt{f}^\alpha(\xi^\alpha))$
\end{center}
which would produce ill-formed terms. 

The sets $FV^I(A)$ and $BV^I(A)$ of respectively the \emph{free} and \emph{bound individual variables} of $A$ can be defined in a usual inductive way, starting from

\begin{center}
    $\star^I(T : \alpha) = FV^I(T)$ and $\star^I(T \equiv U) = \star^I(T) \cup \star^I(U)$
\end{center}
where $\star = FV, BV$ - observe that, if $T : \alpha$, then $\star(\alpha) \subseteq \star^I(T)$ with $\star = FV, BV$. Bound individual variables in quantified formulas give

\begin{center}
    $BV^I(\square \ \nu \ A) = \begin{cases} BV^I(A) \cup \{\nu\} & \text{if} \ \nu \ \text{is an individual variable}\\BV^I(A) & \text{otherwise}\end{cases}$
\end{center}
The sets $FV^T(A)$ and $BV^T(A)$ of respectively the \emph{free} and \emph{bound typed-variables} of $A$, and the sets $FV^F(A)$ and $BV^F(A)$ of respectively the \emph{free} and \emph{bound functional variables} are defined analogously. We say that $A$ is \emph{closed} iff $FV^I(A) = FV^T(A) = FV^F(A) = \emptyset$, and that it is \emph{open} otherwise.

A \emph{substitution of} $x$ \emph{with} $t$ \emph{in} $A$ is a function $\texttt{FORM}_{\texttt{Gen}^\varepsilon} \to \texttt{FORM}_{\texttt{Gen}^\varepsilon}$ defined in a usual inductive way. We only specify that, in the case of quantified formulas, we have

\begin{center}
    $(\square \ \nu \ A)[t/x] = \begin{cases} \square \ \nu[t/x] \ A[t/x] & \text{if} \ \nu \ \text{is a typed- or functional variable}\\ \square \ \nu \ A[t/x] & \text{if} \ \nu \ \text{is an individual variable and} \ \nu \neq x\\\square \ \nu \ A & \text{if} \ \nu \ \text{is an individual variable and} \ \nu = x\end{cases}$
\end{center}
A \emph{substitution of} $\xi^\alpha$ \emph{with} $T : \alpha$ \emph{in} $A$ is a function $\texttt{FORM}_{\texttt{Gen}^\varepsilon} \to \texttt{FORM}_{\texttt{Gen}^\varepsilon}$ defined in a usual inductive way. We only specify that, in the case of quantified formulas, we have

\begin{center}
    $(\square \ \nu \ A)[T/\xi^\alpha_i] = \begin{cases} \square \ \nu \ A[T/\xi^\alpha_i] & \text{if} \ \nu \neq \xi^\alpha_i \\\square \ \nu \ A & \text{if} \ \nu = \xi^\alpha_i \end{cases}$
\end{center}
A \emph{substitution of} $h^\beta_x$ \emph{with} $U : \beta$ in $A$ is a function $\texttt{FORM}_{\texttt{Gen}^\varepsilon} \to \texttt{FORM}_{\texttt{Gen}^\varepsilon}$ that is defined only when $x \notin FV(\alpha)$ for $\xi^\alpha \in FV^T(U)$. It is defined in a usual inductive way. In the case of quantified formulas, it behaves like the substitutions of typed-variables in formulas. The same holds for a \emph{substitution of} $\texttt{f}^\alpha$ \emph{with} $T : \alpha$ in $A$, except that this substitution is defined for every $T : \alpha$.

The notion of $t$ \emph{free for} $x$ \emph{in} $A$ can be specified in a standard inductive way. As for quantified formulas $\square \ \nu \ A$, we must have that either $\nu = x$ or that $\nu \neq x$, $\nu \notin FV(t)$ and $t$ free for $x$ in $A$. Similarly, the notion of $U : \alpha$ \emph{free for} $\xi^\alpha$ or $\texttt{h}^\alpha$ or $\texttt{f}^\alpha$ \emph{in} $A$ is defined in a usual inductive way; for quantified formulas $\square \ \nu \ A$, either $\nu = \xi^\alpha, \texttt{h}^\alpha, \texttt{f}^\alpha$ or $\nu \neq \xi^\alpha, \texttt{h}^\alpha, \texttt{f}^\alpha$, $\nu \notin FV^I(U)$, $\nu \notin FV^T(U)$, $\nu \notin FV^F(U)$ and $U$ free for $\xi^\alpha, \texttt{h}^\alpha, \texttt{f}^\alpha$ in $A$.

\section{A formal system of grounding over $\texttt{Gen}^\varepsilon$}

We now introduce a formal system of grounding over $\texttt{Gen}^\varepsilon$. The system contains three kinds of rules:

\begin{itemize}
    \item type introductions and type eliminations;
    \item equivalence rules, i.e.
    \begin{itemize}
        \item standard rules for reflexivity, simmetry and transitivity, plus preservation of denotation;
        \item rules saying that operational symbols applied to equivalent terms yield equivalent terms;
        \item equations for computing non-canonical terms whose outermost operations are applied to (relevant) canonical terms;
    \end{itemize}
    \item first-order logical rules.
\end{itemize}
We show that three kinds of results are provable \emph{within} the system:

\begin{itemize}
    \item theorems translating in $\texttt{Gen}^\varepsilon$ the ground-clauses of Section 2;
    \item theorems saying that the non-primitive operational symbols of $\texttt{Gen}^\varepsilon$ stand for functions of a certain operational type thanks to the equations associated to them;
    \item theorems showing that if we add to $\texttt{Gen}^\varepsilon$ a non-primitive operational symbol $\phi$, thereby obtaining an expansion $\texttt{Gen}^\varepsilon_1$ of $\texttt{Gen}^\varepsilon$, and if we expand our system $\Sigma$ by adding to it an equation associated to $\phi$ that allows to prove that $\phi$ stands for a function of type $\Gamma \rhd \alpha$ (plus obviously other appropriate equivalence rules), thereby obtaining an expansion $\Sigma_1$ of $\Sigma$, and finally if $\Gamma \vdash \alpha$ is derivable in first-order intuitionistic logic, then for every term $T$ of $\texttt{Gen}^\varepsilon_1$ we can find a term $U$ of $\texttt{Gen}^\varepsilon$ such that $T \equiv U$ is provable in $\Sigma_1$.
\end{itemize}
Finally, we prove a theorem \emph{about} the system, namely, that if $U$ is a term of $\texttt{Gen}$ of type $\beta$ and $FV^T(U) = \{\xi^{\alpha_1}, ..., \xi^{\alpha_n}\}$, then $U : \beta$ is provable in the system under the assumptions $\xi^{\alpha_i} : \alpha_i$ ($1 \leq i \leq n$) only. This meta-theorem can be considered as a kind of \emph{correctness} result for our system - recall the meaning of $U : \alpha$.

\subsection{The system}

The rules of the system come in a Gentzen natural deduction format.

\subsubsection{Typing rules}

Rules for typing split into introductions and eliminations. Type eliminations come in the form of generalized elimination rules [see Schroeder-Heister 1984a, 1984b]. They are in line with Dummett's so-called \emph{fundamental assumption} - i.e. every closed proof of $\alpha$ must reduce to a closed \emph{canonical} proof of $\alpha$ [see Dummett 1991].

\paragraph{Type introductions}
\medskip

\begin{prooftree}
\AxiomC{}
\RightLabel{$\texttt{C}$, for every individual constant $\delta$ naming a $\Delta \in \texttt{DER}_{\texttt{S}}$ with conclusion $\alpha$}
\UnaryInfC{$\delta : \alpha$}
\end{prooftree}

\begin{prooftree}
 \AxiomC{$T : \alpha$}
 \AxiomC{$U : \beta$}
 \RightLabel{$\wedge I$}
 \BinaryInfC{$\wedge I (T, U) : \alpha \wedge \beta$}
 \AxiomC{$T : \alpha_i$}
 \RightLabel{$\vee I$, ($i = 1, 2$)}
 \UnaryInfC{$\vee I[\alpha_i \rhd \alpha_1 \vee \alpha_2](T) : \alpha_1 \vee \alpha_2$}
 \AxiomC{$[\xi^\alpha_i : \alpha]$}
 \noLine
 \UnaryInfC{$\vdots$}
 \noLine
 \UnaryInfC{$T(\xi^\alpha_i) : \beta$}
 \RightLabel{$\rightarrow I$}
 \UnaryInfC{$\rightarrow I \xi^\alpha_i (T(\xi^\alpha_i)) : \alpha \rightarrow \beta$}
 \noLine
 \TrinaryInfC{}
\end{prooftree}
\medskip

\begin{prooftree}
 \AxiomC{$T(x_i) : \alpha(x_i)$}
 \RightLabel{$\forall I$}
 \UnaryInfC{$\forall I y (T(y/x_i)) : \forall y \alpha(y/x_i)$}
 \AxiomC{$T : \alpha(t/x)$}
 \RightLabel{$\exists I$}
 \UnaryInfC{$\exists I[\alpha(t/x) \rhd \exists x \alpha(x)](T) : \exists x \alpha(x)$}
 \noLine
 \BinaryInfC{}
\end{prooftree}
\medskip

\noindent We have the restrictions: in $\rightarrow I$, $\xi^\alpha_i$ must not occur free in undischarged assumption other than $\xi^\alpha_i : \alpha$; in $\forall I$, $x_i$ must not occur free in any undischarged assumption.\footnote{Although seemingly non-standard, the restriction on $\rightarrow I$ actually turns out to be necessary. Without it, we could prove in the system such results as $\mathfrak{E} \texttt{f}^\beta (\rightarrow I \xi^\alpha (\wedge I (\xi^\alpha, \texttt{f}^\beta(\xi^\alpha))) : \alpha \rightarrow \alpha \wedge \beta)$ from the assumption $\mathfrak{E} \texttt{f}^\beta \mathfrak{E} \xi^\alpha (\texttt{f}^\beta(\xi^\alpha) : \beta)$ - i.e., morally, from the fact that there is an operation $\texttt{f}^\beta$ that \emph{for some} ground $\xi^\alpha$ for $\alpha$ yields a ground $\texttt{f}^\beta(\xi^\alpha)$ for $\beta$, we can conclude that there is an operation $\texttt{f}^\beta$ such that $\rightarrow I \xi^\alpha(\wedge I (\xi^\alpha, \texttt{f}^\beta(\xi^\alpha)))$ is a ground for $\alpha \rightarrow \alpha \wedge \beta$, which intuitively does not hold because we have to grant that $\texttt{f}^\beta$ yields a ground for $\beta$ \emph{for every} ground $\xi^\alpha$ for $\alpha$. Also, note the difference between the restriction on \emph{term formation} for $\forall I$, and the restriction on \emph{type introduction} for $\forall I$. In particular, we can have $\forall I x (\texttt{h}^{\alpha(x)}(x))$ \emph{as a term}, but we cannot \emph{prove} $\forall I x (\texttt{h}^{\alpha(x)}(x)) : \forall x \alpha(x)$ from the assumption $\texttt{h}^{\alpha(x)}(x) : \alpha(x)$. Without this restriction, e.g., $\Pi x (\texttt{h}^{\alpha(x)}(x) : \alpha(x))$ would be provable from $\mathfrak{E} x (\texttt{h}^{\alpha(x)}(x) : \alpha(x))$ - i.e. from the fact that, \emph{for some} $x$, $\texttt{h}^{\alpha(x)}(x)$ denotes a ground for $\alpha(x)$ we derive that this holds \emph{for all} $x$. In general, that a term enjoying given properties denotes a ground for a given formula is something that the system \emph{has to prove}, although what one proves is exactly the intended meaning of the term - in this case, an arbitrary term with no free individual variable $x$ occurring free in typed-variables that occur free within it \emph{is required} to stand for an operation on grounds of operational type $\Gamma \rhd \alpha$ for $x$ not free in any element of $\Gamma$, but this requirement must be \emph{proved} by the system. The individual variable $x$ is free in $\texttt{h}^{\alpha(x)}(x) : \alpha(x)$, even if $\texttt{h}^{\alpha(x)}(x)$ should stand for operations on grounds - applied to $x$ - of operational type $\Gamma \rhd \alpha$ for $x$ not free in any element of $\Gamma$. There is thus a difference between an individual variable being free in a term or in a formula, and an individual variable being suitable for quantification through $\forall I$. The same difference applies, in Prawitz's semantics of valid arguments, to closed valid arguments for $\alpha(x)$, to which we can apply licitly Gentzen's introduction for $\forall$, but where $x$ may occur free - and eventually be replaced by other individual terms before quantification. This observations also apply to $\Pi$-introduction displayed below - labelled with ($\Pi_I$) - which undergoes the same restrictions as those required for the $\forall I$-rule.}

\paragraph{Type eliminations}
\medskip

\begin{prooftree}
 \AxiomC{$T : \alpha \wedge \beta$}
 \AxiomC{$[T \equiv \wedge I(\xi^\alpha, \xi^\beta)]$}
 \AxiomC{$[\xi^\alpha : \alpha]$}
 \AxiomC{$[\xi^\beta : \beta]$}
 \noLine
 \TrinaryInfC{$\vdots$}
 \noLine
 \UnaryInfC{$A$}
 \RightLabel{$D_\wedge$}
 \BinaryInfC{$A$}
\end{prooftree}
\medskip

\begin{prooftree}
 \AxiomC{$T : \alpha \vee \beta$}
 \AxiomC{$[T \equiv \vee I(\xi^\alpha)]$}
 \AxiomC{$[\xi^\alpha : \alpha]$}
 \noLine
 \BinaryInfC{$\vdots$}
 \noLine
 \UnaryInfC{$A$}
 \AxiomC{$[T \equiv \vee I (\xi^\beta)]$}
 \AxiomC{$[\xi^\beta : \beta]$}
 \noLine
 \BinaryInfC{$\vdots$}
 \noLine
 \UnaryInfC{$A$}
 \RightLabel{$D_\vee$}
 \TrinaryInfC{$A$}
\end{prooftree}
\medskip

\noindent The two minor premises $A$ are the conclusions of two \emph{distinct} sub-derivations: the leftmost (possibly) depends on the assumptions $T \equiv \vee I(\xi^\alpha)$ and $\xi^\alpha : \alpha$ - discharged \emph{only in this sub-derivation} - while the rightmost (possibly) depends on the assumptions $T \equiv \vee I(\xi^\beta)$ and $\xi^\beta : \beta$ - discharged \emph{only in this sub-derivation}.

\begin{prooftree}
 \AxiomC{$T : \alpha \rightarrow \beta$}
 \AxiomC{$[T \equiv \ \rightarrow I \xi^\alpha (\texttt{f}^\beta(\xi^\alpha))]$}
 \AxiomC{$[\Pi \xi^\alpha(\xi^\alpha : \alpha \supset \texttt{f}^\beta(\xi^\alpha)) : \beta))]$}
 \noLine
 \BinaryInfC{$\vdots$}
 \noLine
 \UnaryInfC{$A$}
 \RightLabel{$D_\rightarrow$}
 \BinaryInfC{$A$}
\end{prooftree}
\medskip

\begin{prooftree}
 \AxiomC{$T : \forall x \alpha(x)$}
 \AxiomC{$[T \equiv \forall I x (\texttt{h}^{\alpha(x)}(x))]$}
 \AxiomC{$[\Pi x(\texttt{h}^{\alpha(x)}(x) : \alpha(x))]$}
 \noLine
 \BinaryInfC{$\vdots$}
 \noLine
 \UnaryInfC{$A$}
 \RightLabel{$D_\forall$}
 \BinaryInfC{$A$}
\end{prooftree}
\medskip

\begin{prooftree}
 \AxiomC{$T : \exists x \alpha(x)$}
 \AxiomC{$[T \equiv \exists I(\xi^{\alpha(x)})]$}
 \AxiomC{$[\xi^{\alpha(x)} : \alpha(x)]$}
 \noLine
 \BinaryInfC{$\vdots$}
 \noLine
 \UnaryInfC{$A$}
 \RightLabel{$D_\exists$}
 \BinaryInfC{$A$}
\end{prooftree}
\medskip

\noindent We have restrictions on \emph{proper variables}, e.g.: in $D_\rightarrow$, $\texttt{f}^\beta$ must not occur free in any undischarged assumption on which $A$ depends, other than those discharged by the rule; in $D_\forall$, for every $t \in \texttt{TERM}_L$, $\texttt{h}^{\alpha(t)}$ must not occur free in any undischarged assumption on which $A$ depends (in what we obtain by replacing $x$ with $t$ in the derivation of the minor premise), other than those discharged by the rule (with substitution of $x$ with $t$); in $D_\exists$, $x$ and $\xi^{\alpha(x)}$ must not occur free in any undischarged assumption on which $A$ depends other than those discharged by the rule. We also have a rule for eliminating the type $\bot$:

\begin{prooftree}
 \AxiomC{$T : \bot$}
 \RightLabel{$\bot$}
 \UnaryInfC{$\bot^G$}
\end{prooftree}

As a final remark, observe that the term $T$ mentioned in the major premise of the type elimination rules \emph{need not} begin with a primitive operation - it may be, for example, a mere typed-variable. Nonetheless, the leftmost assumption of the minor premises equates $T$ with a term that begins with a primitive operation. This is because, as said, the type elimination rules are inspired by Dummett's fundamental assumption. The idea behind the rules is thus the following. If $T$ denotes a ground for a formula $\alpha$ with main logical sign $s$ - what the major premise expresses - then $T$ must be equivalent to a term that \emph{directly} denotes a ground for $\alpha$, i.e. a term where the primitive operation $sI$ is applied to grounds for the sub-formula(s) of $\alpha$; therefore, the consequences of $T$ denoting a ground for $\alpha$ include those that may be derived from the equivalence of $T$ with a term that directly denotes a ground for $\alpha$ - what the assumptions that the minor premise depends on express.

\subsubsection{Equivalence rules}

Equivalence rules split in: standard rules for equivalence, plus preservation of denotation; rules for application of operations to equivalent terms; equations ruling the computation of non-primitive operational symbols.

\paragraph{Reflexivity, symmetry, transitivity, preservation of denotation}
\medskip

\begin{prooftree}
 \AxiomC{}
 \RightLabel{$\equiv_R$}
 \UnaryInfC{$T \equiv T$}
 \AxiomC{$T \equiv U$}
 \RightLabel{$\equiv_S$}
 \UnaryInfC{$U \equiv T$}
 \AxiomC{$U \equiv Z$}
 \AxiomC{$Z \equiv W$}
 \RightLabel{$\equiv_T$}
 \BinaryInfC{$U \equiv W$}
 \AxiomC{$T \equiv U$}
 \AxiomC{$U : \alpha$}
 \RightLabel{$\equiv_P$}
 \BinaryInfC{$T : \alpha$}
 \noLine
 \QuaternaryInfC{}
\end{prooftree}
\medskip

\paragraph{Application of operational symbols to equivalent terms (and vice versa in the canonical case)}
\medskip

\begin{prooftree}
 \AxiomC{$T \equiv U$}
 \AxiomC{$V \equiv W$}
 \RightLabel{$\equiv^\wedge_1$}
 \BinaryInfC{$\wedge I(T, V) \equiv \wedge I (U, W)$}
 \AxiomC{$\wedge I (T_1, T_2) \equiv \wedge I (U_1, U_2)$}
 \RightLabel{$\equiv^\wedge_{2, i}$, $i = 1, 2$}
 \UnaryInfC{$T_i \equiv U_i$}
 \AxiomC{$T \equiv U$}
 \RightLabel{$\equiv^\wedge_{3, i}$, $i = 1, 2$}
 \UnaryInfC{$\wedge_{E, i}(T) \equiv \wedge_{E, i}(U)$}
 \noLine
 \TrinaryInfC{}
\end{prooftree}
\medskip

\begin{footnotesize}
\begin{prooftree}
 \AxiomC{$T \equiv U$}
 \RightLabel{$\equiv^\vee_1$}
 \UnaryInfC{$\vee I(T) \equiv \vee I(U)$}
 \AxiomC{$\vee I[\alpha_i \rhd \alpha_1 \vee \alpha_2](T) \equiv \vee I[\alpha_i \rhd \alpha_1 \vee \alpha_2](U)$}
 \RightLabel{$\equiv^\vee_2$}
 \UnaryInfC{$T \equiv U$}
 \AxiomC{$T \equiv U$}
 \AxiomC{$Z_1 \equiv Z_2$}
 \AxiomC{$W_1 \equiv W_2$}
 \RightLabel{$\equiv^\vee_3$}
 \TrinaryInfC{$\vee E \xi^\alpha_i \xi^\beta_i(T, Z_1, W_1) \equiv \vee E \xi^\alpha_i \xi^\beta_i(U, Z_2, W_2)$}
 \noLine
 \TrinaryInfC{}
\end{prooftree}
\end{footnotesize}
\medskip

\begin{prooftree}
 \AxiomC{$T \equiv U$}
 \RightLabel{$\equiv^\rightarrow_1$}
 \UnaryInfC{$\rightarrow I \xi^\alpha_i(T) \equiv \ \rightarrow I \xi^\alpha_i(U)$}
 \AxiomC{$\rightarrow I \xi^\alpha(T) \equiv \ \rightarrow I \xi^\alpha(U)$}
 \RightLabel{$\equiv^\rightarrow_2$}
 \UnaryInfC{$T \equiv U$}
 \AxiomC{$T_1 \equiv T_2$}
 \AxiomC{$U_1 \equiv U_2$}
 \RightLabel{$\equiv^\rightarrow_3$}
 \BinaryInfC{$\rightarrow E(T_1, U_1) \equiv \ \rightarrow E(T_2, U_2)$}
 \noLine
 \TrinaryInfC{}
\end{prooftree}
\medskip

\begin{prooftree}
 \AxiomC{$T \equiv U$}
 \RightLabel{$\equiv^\forall_1$}
 \UnaryInfC{$\forall I x_i(T) \equiv \forall I x_i(U)$}
 \AxiomC{$\forall I x (T) \equiv \forall I x(U)$}
 \RightLabel{$\equiv^\forall_2$}
 \UnaryInfC{$T \equiv U$}
 \AxiomC{$T \equiv U$}
 \RightLabel{$\equiv^\forall_3$}
 \UnaryInfC{$\forall E(T) \equiv \forall E(U)$}
 \noLine
 \TrinaryInfC{}
\end{prooftree}
\medskip

\begin{prooftree}
 \AxiomC{$T \equiv U$}
 \RightLabel{$\equiv^\exists_1$}
 \UnaryInfC{$\exists I(T) \equiv \exists I(U)$}
 \AxiomC{$\exists I(T) \equiv \exists I(U)$}
 \RightLabel{$\equiv^\exists_2$}
 \UnaryInfC{$T \equiv U$}
 \AxiomC{$T \equiv U$}
 \AxiomC{$Z \equiv V$}
 \RightLabel{$\equiv^\exists_3$}
 \BinaryInfC{$\exists E \ x_i \ \xi^{\alpha(x_i)}_j (T, Z) \equiv \exists E \ x_i \ \xi^{\alpha(x_i)}_j (U, V)$}
 \noLine
 \TrinaryInfC{}
\end{prooftree}
\medskip

\noindent The rules concerning symbols that bind variables could be more fine-grained. One could require that if the indexes of the typed variables are such as to yield vacuous bindings, then the two terms can be said to be equivalent. So for example, we do not want $\rightarrow \xi^\alpha_1(\xi^\alpha_1) \equiv \ \rightarrow \xi^\alpha_2(\xi^\alpha_1)$, but we may want $\rightarrow \xi^\alpha_2(\xi^\alpha_1) \equiv \ \rightarrow \xi^\alpha_3(\xi^\alpha_1)$. However, we shall not deal with these details here.

\paragraph{Equations for computing non-canonical terms}
\medskip

\begin{prooftree}
 \AxiomC{}
 \RightLabel{$R_{\wedge, i}$, $i = 1, 2$}
 \UnaryInfC{$\wedge_{E, i}(\wedge I(T_1, T_2)) \equiv T_i$}
 \AxiomC{}
 \RightLabel{$R_{\vee, i}$, $i = 1, 2$}
 \UnaryInfC{$\vee E \xi^{\alpha_2} \xi^{\alpha_2}(\vee I [\alpha_i \rhd \alpha_1 \vee \alpha_2)(T), U_1(\xi^{\alpha_1}), U_2(\xi^{\alpha_2})) \equiv U_i(T)$}
 \noLine
 \BinaryInfC{}
\end{prooftree}
\medskip

\begin{prooftree}
 \AxiomC{}
 \RightLabel{$R_\rightarrow$}
 \UnaryInfC{$\rightarrow E(\rightarrow I \xi^\alpha(T(\xi^\alpha)), U) \equiv T(U)$}
 \AxiomC{}
 \RightLabel{$R_\forall$}
 \UnaryInfC{$\forall E[\forall x \alpha(x) \rhd \alpha(t/x)](\forall I x(T(x))) \equiv T(t)$}
 \noLine
 \BinaryInfC{}
\end{prooftree}
\medskip

\begin{prooftree}
 \AxiomC{}
 \RightLabel{$R_\exists$}
 \UnaryInfC{$\exists E \ x \ \xi^{\alpha(x)}(\exists I[\alpha(t/x) \rhd \exists x \alpha(x)](T), U(x, \xi^{\alpha(x)})) = U(t, T)$}
\end{prooftree}

\subsubsection{Logic}

The logical rules amount to standard Gentzen natural deduction rules for first-order intuitionistic logic. We just give those for the quantifiers and the explosion principle, since the others are straightforward - observe that, in particular, we have an assumption rule, saying that any \emph{formula} of $\texttt{Gen}^\varepsilon$ can be assumed.
\medskip

\begin{prooftree}
 \AxiomC{$A(\nu)$}
 \RightLabel{($\Pi_I$)}
 \UnaryInfC{$\Pi \mu A(\mu/\nu)$}
 \AxiomC{$\Pi \nu A(\nu)$}
 \RightLabel{($\Pi_E$)}
 \UnaryInfC{$A(\tau/\nu)$}
 \AxiomC{$A(\tau/\nu)$}
 \RightLabel{($\mathfrak{E}_I$)}
 \UnaryInfC{$\mathfrak{E} \nu A(\nu)$}
 \AxiomC{$\mathfrak{E} \nu A(\mu/\nu)$}
 \AxiomC{$[A(\mu)]$}
 \noLine
 \UnaryInfC{$\vdots$}
 \noLine
 \UnaryInfC{$B$}
 \RightLabel{($\mathfrak{E}_E$)}
 \BinaryInfC{$B$}
 \AxiomC{$\bot^G$}
 \RightLabel{($\bot^G$)}
 \UnaryInfC{$A$}
 \noLine
 \QuinaryInfC{}
\end{prooftree}
\medskip

\noindent In these rules, $\nu, \mu$ are individual, typed of functional variables, while $\tau$ stands for a term of $L$ or for a term of $\texttt{Gen}$. We have standard restrictions on \emph{proper variables}, plus the following restriction for ($\Pi_E$) and ($\mathfrak{E}_I$): if $\nu$ is $\xi^\alpha$ or $\texttt{f}^\alpha$, $\tau$ is of type $\alpha$; if $\nu$ is $\texttt{h}^\beta_x$, $\tau$ is of type $\beta$ and such that $x \notin FV(\alpha)$ for $\xi^\alpha \in FV^T(\tau)$. 

\subsubsection{Derivations and remarks}

We call our system $\texttt{GG}$, i.e. \emph{Gentzen-grounding}. Derivations in $\texttt{GG}$ can be defined in a usual inductive way. The base of the inductive definition states that the single node consisting of an application of the assumption rule is a derivation, and that the application of a no-premises rule (usually understood as an axiom in the Gentzen format) is a derivation. The way in which typing rules are combined with logical rules is shown in examples below. We indicate with $\mathfrak{G} \vdash_{\texttt{GG}} A$ the fact that $A$ is derivable in $\texttt{GG}$ from a set of assumptions $\mathfrak{G}$. We sometimes indicate with $\texttt{DER}_{\texttt{GG}}$ the set of the derivations of $\texttt{GG}$. Observe that type elimination rules are inter-derivable with elimination rules of a more straightforward shape, e.g.

\begin{prooftree}
 \AxiomC{$T : \exists x \alpha(x)$}
 \UnaryInfC{$\mathfrak{E} x \mathfrak{E} \xi^{\alpha(x)}(T \equiv \exists I[\alpha(x) \rhd \exists x \alpha(x)](\xi^{\alpha(x)}) \times \xi^{\alpha(x)} : \alpha(x))$}
\end{prooftree}
We chose to adopt generalized elimination rules because they render somewhat easier the proof of some results in Section 6.

One may wonder why we do not have the inverses of the equivalence rules for equivalent applications of non-primitive operational symbols. These rules are indeed not plausible because, if we had them, we could derive for example the following

\begin{prooftree}
 \AxiomC{}
 \RightLabel{$R_{\wedge, 1}$}
 \UnaryInfC{$\wedge_{E, 1}(\wedge I(\xi^\alpha, \xi^\beta)) \equiv \xi^\alpha$}\AxiomC{}
 \RightLabel{$R_{\wedge, 1}$}
 \UnaryInfC{$\wedge_{E, 1}(\wedge I(\xi^\alpha, \xi^\gamma)) \equiv \xi^\alpha$}
 \RightLabel{$\equiv_S$}
 \UnaryInfC{$\xi^\alpha \equiv \wedge_{E, 1}(\wedge I(\xi^\alpha, \xi^\gamma))$}
 \RightLabel{$\equiv_T$}
 \BinaryInfC{$\wedge_{E, 1}(\wedge I(\xi^\alpha, \xi^\beta)) \equiv \wedge_{E, 1}(\wedge I(\xi^\alpha, \xi^\gamma))$}
 \UnaryInfC{$\wedge I(\xi^\alpha, \xi^\beta) \equiv \wedge I (\xi^\alpha, \xi^\gamma)$}
\end{prooftree}
which is clearly undesirable when $\beta \neq \gamma$ - and of course violates the restriction for the well-forming of formulas with $\equiv$ as main sign.

\subsection{Some theorems within the system}

We list some theorems that can be proved within the system.

\subsubsection{Theorems expressing the ground-clauses}

We indicate bi-conditional with $\Leftrightarrow$, and as usual we put $A \Leftrightarrow B \stackrel{def}{=} (A \supset B) \times (B \supset A)$. The following results hold:

\begin{itemize}
    \item[($\wedge$)*] $\vdash_{\texttt{GG}} \Pi \xi^{\alpha \wedge \beta} (\xi^{\alpha \wedge \beta} : \alpha \wedge \beta \Leftrightarrow \mathfrak{E} \xi^\alpha \mathfrak{E} \xi^\beta (\xi^{\alpha \wedge \beta} \equiv \wedge I (\xi^\alpha, \xi^\beta) \times (\xi^\alpha : \alpha \times \xi^\beta : \beta)))$
    \item[($\vee$)*] $\vdash_{\texttt{GG}} \Pi \xi^{\alpha \vee \beta} (\xi^{\alpha \vee \beta} : \alpha \vee \beta \Leftrightarrow \mathfrak{E} \xi^\alpha \mathfrak{\xi^\beta}((\xi^{\alpha \vee \beta} \equiv \vee I (\xi^\alpha) \times \xi^\alpha : \alpha) + (\xi^{\alpha \vee \beta} \equiv \vee I (\xi^\beta) \times \xi^\beta : \beta)))$
    \item[($\rightarrow$)*] $\vdash_{\texttt{GG}} \Pi \xi^{\alpha \rightarrow \beta} (\xi^{\alpha \rightarrow \beta} : \alpha \rightarrow \beta \Leftrightarrow \mathfrak{E} \texttt{f}^\beta (\xi^{\alpha \rightarrow \beta} \equiv \ \rightarrow I \xi^\alpha (\texttt{f}^\beta(\xi^\alpha)) \times \Pi \xi^\alpha (\xi^\alpha : \alpha \supset \texttt{f}^\beta (\xi^\alpha) : \beta)))$
    \item[($\forall$)*] $\vdash_{\texttt{GG}} \Pi \xi^{\forall x \alpha(x)} (\xi^{\forall x \alpha(x)} : \forall x \alpha(x) \Leftrightarrow \mathfrak{E} \texttt{h}^{\alpha(x)} (\xi^{\forall x \alpha(x)} \equiv \forall I x (\texttt{h}^{\alpha(x)}(x)) \times \Pi x (\texttt{h}^{\alpha(x)}(x) : \alpha (x))))$
    \item[($\exists$)*] $\vdash_{\texttt{GG}} \Pi \xi^{\exists x \alpha(x)} (\xi^{\exists x \alpha(x)} : \exists x \alpha(x) \Leftrightarrow \mathfrak{E} x \mathfrak{E} \xi^{\alpha(x)} (\xi^{\exists x \alpha(x)} \equiv \exists I (\xi^{\alpha(x)}) \times \xi^{\alpha(x)} : \alpha(x)))$
\end{itemize}
As can be seen, ($k$)* can be considered as a translation in $\texttt{Gen}^\varepsilon$ of the corresponding ground-clause ($k$) of Section 2. Given our type introduction rules, the fact that the clauses are provable is of course not surprising.

We now exemplify one of the derivations in $\texttt{GG}$ establishing the theorems above. It is the derivation of ($\exists$)*. The left-to-right implication can be proved by applying ($\supset_I$) to the following derivation:

\begin{prooftree}
\AxiomC{$\xi^{\exists x \alpha (x)} : \exists x \alpha(x)$}
\AxiomC{$1$}
\noLine
\UnaryInfC{$[\xi^{\exists x \alpha(x)} \equiv \exists I (\xi^{\alpha(x)})]$}
\AxiomC{$2$}
\noLine
\UnaryInfC{$[\xi^{\alpha(x)} : \alpha(x)]$}
\RightLabel{($\times_I$)}
\BinaryInfC{$\xi^{\exists x \alpha(x)} \equiv \exists I (\xi^{\alpha(x)}) \times \xi^{\alpha(x)} : \alpha(x)$}
\RightLabel{($\mathfrak{E}_I$)}
\UnaryInfC{$\mathfrak{E} \xi^{\alpha(x)} (\xi^{\exists x \alpha(x)} \equiv \exists I (\xi^{\alpha(x)}) \times \xi^{\alpha(x)} : \alpha(x))$}
\RightLabel{($\mathfrak{E}_I$)}
\UnaryInfC{$\mathfrak{E} x \mathfrak{E} \xi^{\alpha(x)} (\xi^{\exists x \alpha(x)} \equiv \exists I (\xi^{\alpha(x)}) \times \xi^{\alpha(x)} : \alpha(x))$}
\RightLabel{$D_\exists$, $1, 2$}
\BinaryInfC{$\mathfrak{E} x \mathfrak{E} \xi^{\alpha(x)} (\xi^{\exists x \alpha(x)} \equiv \exists I (\xi^{\alpha(x)}) \times \xi^{\alpha(x)} : \alpha(x))$}
\end{prooftree}
\medskip 

\noindent The right-to-left implication can be proved by applying twice ($\mathfrak{E}_E$) and then ($\supset_I$) to the following derivation from the assumption $\mathfrak{E} x \mathfrak{E} \xi^{\alpha(x)} (\xi^{\exists x \alpha(x)} \equiv \exists I (\xi^{\alpha(x)}) \times \xi^{\alpha(x)} : \alpha(x))$

\begin{prooftree}
\AxiomC{$\xi^{\exists x \alpha(x)} \equiv \exists I (\xi^{\alpha(x)}) \times \xi^{\alpha(x)} : \alpha(x)$}
\RightLabel{($\times_{E, 2}$)}
\UnaryInfC{$\xi^{\alpha(x)} : \alpha(x)$}
\RightLabel{$\exists I$}
\UnaryInfC{$\exists I (\xi^{\alpha(x)}) : \exists x \alpha(x)$}
\AxiomC{$\xi^{\exists x \alpha(x)} \equiv \exists I (\xi^{\alpha(x)}) \times \xi^{\alpha(x)} : \alpha(x)$}
\RightLabel{($\times_{E, 1}$)}
\UnaryInfC{$\xi^{\exists x \alpha(x)} \equiv \exists I (\xi^{\alpha(x)})$}
\RightLabel{$\equiv_P$}
\BinaryInfC{$\xi^{\exists x \alpha(x)} : \exists x \alpha(x)$}
\end{prooftree}
\medskip

\noindent Then one takes the two derivations above and introduces $\times$. Finally, one universally quantifies over $\xi^{\exists x \alpha(x)}$ by ($\Pi_I$).

The ground-clauses of Section 2 are usually accompanied by a clause saying that every closed derivation in the atomic system $\texttt{S}$ of an atomic base $\mathfrak{B}$ with conclusion $\alpha$ is a ground for $\alpha$ on $\mathfrak{B}$, and by a clause saying that every $\mathfrak{B}$ is consistent, i.e. that $\nvdash_{\texttt{S}} \bot$ - or, as one usually puts it, that there is no ground on $\mathfrak{B}$ for $\bot$. The first of such clauses is clearly provable in $\texttt{GG}$ because we have assumed as an axiom that, for every individual constant $\delta$ naming a closed $\Delta \in \texttt{DER}_{\texttt{S}}$ with conclusion $\alpha$, $\delta : \alpha$ holds. As regards the clause for $\bot$, it is provable too, i.e. we have $\vdash_{\texttt{GG}} \neg \mathfrak{E} \xi^\bot (\xi^\bot : \bot)$. Finally, observe that we can easily prove under what conditions we are allowed to apply the operational symbol $\bot_\alpha$, i.e. $\vdash_{\texttt{GG}} \Pi \xi^\bot (\xi^\bot : \bot \supset \bot_\alpha(\xi^\bot) : \alpha)$.

\subsubsection{Checking the equations}

The type of an operational symbol is intended, not only to provide a ‘‘syntactic" typing, but also to indicate that the symbol itself stands for an operation on grounds of that type. So for example, if $T$ stands for a ground for $\alpha_1 \wedge \alpha_2$, we expect $\wedge_{E, i}(T)$ to stand for a ground for $\alpha_i$ ($i = 1, 2$) [the expression ‘‘stand for" is of course ambiguous here, but a precise development of it can be given in terms of the denotation functions introduced in d'Aragona 2021a]. That the operational symbols of $\texttt{Gen}$ actually behave in this way should be granted by the equations associated to them in $\texttt{GG}$. We could therefore expect that $\texttt{GG}$ proves that the operational symbols are well-defined, and this is indeed the case. In other words, the following results hold:

\begin{itemize}
    \item[($\wedge^{\text{w}}$)] $\vdash_{\texttt{GG}} \Pi \xi^{\alpha_1 \wedge \alpha_2}(\xi^{\alpha_1 \wedge \alpha_2} : \alpha_1 \wedge \alpha_2 \supset \wedge_{E, i}(\xi^{\alpha_1 \wedge \alpha_2}) : \alpha_i)$
    \item[($\vee^{\text{w}}$)] $\vdash_{\texttt{GG}} \Pi \xi^{\alpha_1 \vee \alpha_2} \Pi \texttt{f}^\beta_1 \Pi \texttt{f}^\beta_2 (((\xi^{\alpha_1 \vee \alpha_2} : \alpha_1 \vee \alpha_2 \times \Pi \xi^{\alpha_1} (\xi^{\alpha_1} : \alpha_1 \supset \texttt{f}^\beta_1(\xi^{\alpha_1}) : \beta)) \times$
    
    $\Pi \xi^{\alpha_2} (\xi^{\alpha_2} : \alpha_2 \supset \texttt{f}^\beta_2(\xi^{\alpha_2}) : \beta)) \supset \vee E \xi^{\alpha_1} \xi^{\alpha_2} (\xi^{\alpha_1 \vee \alpha_2}, \texttt{f}^\beta_1(\xi^{\alpha_1}), \texttt{f}^\beta_2(\xi^{\alpha_2})) : \beta)$
    \item[($\rightarrow^{\text{w}}$)] $\vdash_{\texttt{GG}} \Pi \xi^{\alpha \rightarrow \beta} \Pi \xi^\alpha (\xi^{\alpha \rightarrow \beta} : \alpha \rightarrow \beta \times \xi^\alpha : \alpha \supset \ \rightarrow E (\xi^{\alpha \rightarrow \beta}, \xi^\alpha) : \beta)$
    \item[($\forall^{\text{w}}$)] $\vdash_{\texttt{GG}} \Pi \xi^{\forall x \alpha(x)} (\xi^{\forall x \alpha(x)} : \forall x \alpha(x) \supset \forall E (\xi^{\forall x \alpha(x)}) : \alpha(t))$
    \item[($\exists^{\text{w}}$)] $\vdash_{\texttt{GG}} \Pi \xi^{\exists x \alpha(x)} \Pi \texttt{f}^\beta (\xi^{\exists x \alpha(x)} : \exists x \alpha(x) \times \Pi x \Pi \xi^{\alpha(x)} (\xi^{\alpha(x)} : \alpha(x) \supset \texttt{f}^\beta(\xi^{\alpha(x)}) : \beta)$
    
    $\supset \exists E \ x \ \xi^{\alpha(x)} (\xi^{\exists x \alpha(x)}, \texttt{f}^\beta(\xi^{\alpha(x)})) : \beta)$
\end{itemize}
It is easily seen that ($k^{\text{w}}$) proves that $kE$ yields something which is a ground for the co-domain of $kE$, under the assumption that the arguments to which $kE$ is applied are grounds for the domain of $kE$ - remember that $T : \alpha$ is just an abbreviation for $Gr(T, \alpha)$.

As an example, we prove ($\exists^{\text{w}}$), that we abbreviate with $\Pi \xi^{\exists x \alpha(x)} \Pi \texttt{f}^\beta Th_\exists$. First, we eliminate conjunction on an assumption that corresponds to the antecedent of $Th_\exists$ - we indicate this derivation with $\Delta_0$ -

\begin{prooftree}
\AxiomC{$1$}
\noLine
\UnaryInfC{$[\xi^{\exists x \alpha(x)} : \exists x \alpha(x) \times \Pi x \Pi \xi^{\alpha(x)} (\xi^{\alpha(x)} : \alpha(x) \supset \texttt{f}^\beta(\xi^{\alpha(x)}) : \beta)]$}
\RightLabel{($\times_{E, 1}$)}
\UnaryInfC{$\xi^{\exists x \alpha(x)} : \exists x \alpha(x)$}
\end{prooftree}
\medskip

\noindent We now want to apply type elimination to the conclusion of $\Delta_0$, so to obtain the consequent of $Th_\exists$. So we prove that the application of $\exists E$ to $\xi^{\exists x \alpha(x)}$ and to an arbitrary term in the free individual variable $x$ and in the free typed-variable $\xi^{\alpha(x)}$ is equal to an application of $\exists E$ to an arbitrary canonical term and to the same arbitrary term as above - we call this derivation $\Delta^1_1$ -

\begin{prooftree}
\AxiomC{$2$}
\noLine
\UnaryInfC{$[\xi^{\exists x \alpha(x)} \equiv \exists I (\xi^{\alpha(x)})]$}
\AxiomC{}
\RightLabel{$\equiv_R$}
\UnaryInfC{$\texttt{f}^\beta(\xi^{\alpha(x)}) \equiv \texttt{f}^\beta(\xi^{\alpha(x)})$}
\RightLabel{$\equiv^\exists_3$}
\BinaryInfC{$\exists E \ x \ \xi^{\alpha(x)} (\xi^{\exists x \alpha(x)}, \texttt{f}^\beta(\xi^{\alpha(x)})) \equiv \exists E \ x \ \xi^{\alpha(x)} (\exists I (\xi^{\alpha(x)}), \texttt{f}^\beta(\xi^{\alpha(x)}))$}
\end{prooftree}
\medskip

\noindent Then we apply the equation for $\exists E$ to prove via transitivity that the application of $\exists E$ to $\xi^{\exists x \alpha(x)}$ and to an arbitrary term in the free individual variable $x$ and in the free typed-variable $\xi^{\alpha(x)}$ is equivalent to something of type $\beta$, i.e. we have the derivation $\Delta^2_1$

\begin{prooftree}
\AxiomC{$\Delta^1_1$}
\AxiomC{}
\RightLabel{$R_\exists$}
\UnaryInfC{$\exists E \ x \ \xi^{\alpha(x)} (\exists I (\xi^{\alpha(x)}), \texttt{f}^\beta(\xi^{\alpha(x)}))) \equiv \texttt{f}^\beta(\xi^{\alpha(x)})$}
\RightLabel{$\equiv_T$}
\BinaryInfC{$\exists E \ x \ \xi^{\alpha(x)} (\xi^{\exists x \alpha(x)}, \texttt{f}^\beta(\xi^{\alpha(x)}) \equiv \texttt{f}^\beta(\xi^{\alpha(x)})$}
\end{prooftree}
\medskip

\noindent and the derivation $\Delta_2$

\begin{prooftree}
\AxiomC{$1$}
\noLine
\UnaryInfC{$[\xi^{\exists x \alpha(x)} : \exists x \alpha(x) \times \Pi x \Pi \xi^{\alpha(x)} (\xi^{\alpha(x)} : \alpha(x) \supset \texttt{f}^\beta(\xi^{\alpha(x)}) : \beta)]$}
\RightLabel{($\times_{E, 2}$)}
\UnaryInfC{$\Pi x \Pi \xi^{\alpha(x)} (\xi^{\alpha(x)} : \alpha(x) \supset \texttt{f}^\beta(\xi^{\alpha(x)}) : \beta)$}
\RightLabel{($\Pi_E$)}
\UnaryInfC{$\Pi \xi^{\alpha(x)} (\xi^{\alpha(x)} : \alpha(x) \supset \texttt{f}^\beta(\xi^{\alpha(x)}) : \beta)$}
\RightLabel{($\Pi_E$)}
\UnaryInfC{$\xi^{\alpha(x)} : \alpha(x) \supset \texttt{f}^\beta(\xi^{\alpha(x)}) : \beta$}
\AxiomC{$3$}
\noLine
\UnaryInfC{$[\xi^{\alpha(x)} : \alpha(x)]$}
\RightLabel{($\supset_E$)}
\BinaryInfC{$\texttt{f}^\beta(\xi^{\alpha(x)}) : \beta$}
\end{prooftree}
\medskip

\noindent Finally, via preservation of denotation, we can apply type elimination

\begin{prooftree}
\AxiomC{$\Delta_0$}
\AxiomC{$\Delta^2_1$}
\AxiomC{$\Delta_2$}
\RightLabel{$\equiv_P$}
\BinaryInfC{$\exists E \ x \ \xi^{\alpha(x)} (\xi^{\exists x \alpha(x)}, \texttt{f}^\beta(\xi^{\alpha(x)})) : \beta$}
\RightLabel{$D_\exists$, $2, 3$}
\BinaryInfC{$\exists E \ x \ \xi^{\alpha(x)} (\xi^{\exists x \alpha(x)}, \texttt{f}^\beta(\xi^{\alpha(x)})) : \beta$}
\RightLabel{($\supset_I$), $1$}
\UnaryInfC{$Th_\exists$}
\end{prooftree}
Then one applies ($\Pi_I$) twice on $\texttt{f}^\beta$ and on $\xi^{\exists x \alpha(x)}$ to obtain $\Pi \xi^{\exists x \alpha(x)} \Pi \texttt{f}^\beta Th_\exists$.\footnote{Observe that, strictly speaking, we should require $\texttt{f}^\beta$ to range over terms that stand for operations on grounds of operational type $\Gamma \rhd \beta$ \emph{with $x$ not free in any element of $\Gamma$ other than $\alpha(x)$, and $x$ not free in $\beta$}. Thus, we should have introduced a \emph{third} kind of functional variables, say $\texttt{f}^\beta_x\{\alpha(x)\}$, which in a way merge the two kinds we already have. The substitution of terms $U$ of type $\beta$ for such functional variables in formulas has then to respect the principle that $U$ is of type $\beta$, $x \notin FV^I(\beta)$, $x \notin FV^I(\delta)$ for $\xi^\delta \in FV^T(U)$ and $\delta \neq \alpha(x)$. We have omitted this third kind of functional variables not to overburden the notation, and also because the restriction on $\exists E$ in definition 1 is sufficient for obtaining the same.}

\subsubsection{Derived types}

Let us build an expansion $\texttt{Gen}^{++}$ of $\texttt{Gen}^\varepsilon$ by adding a non-primitive operational symbol $F[\alpha \vee \beta, \neg \alpha \rhd \beta]$. Let us then build an expansion $\texttt{GG}^+$ of $\texttt{GG}$ by adding the following equivalence rules

\begin{prooftree}
\AxiomC{$T \equiv U$}
\AxiomC{$V \equiv W$}
\RightLabel{$\equiv^F$}
\BinaryInfC{$F(T, V) \equiv F(U, W)$}
\AxiomC{}
\RightLabel{$R_F$}
\UnaryInfC{$F(\vee I[\beta \rhd \alpha \vee \beta](T), U) \equiv T$}
\noLine
\BinaryInfC{}
\end{prooftree}
\medskip

\noindent In $\texttt{GG}^+$ we can prove that every term of $\texttt{Gen}^{++}$ can be rewritten as a term of $\texttt{Gen}^\varepsilon$, i.e.:

\begin{itemize}
    \item[] $\vdash_{\texttt{GG}^+} \Pi \xi^{\alpha \vee \beta} \ \Pi \xi^{\neg \alpha} \ (\xi^{\alpha \vee \beta} : \alpha \vee \beta \times \xi^{\neg \alpha} : \neg \alpha \supset F(\xi^{\alpha \vee \beta}, \xi^{\neg \alpha}) \equiv \vee E \xi^{\alpha} \xi^{\beta} (\xi^{\alpha \vee \beta}, \bot_\beta(\rightarrow E (\xi^{\alpha}, \xi^{\neg \alpha})), \xi^\beta))$.
\end{itemize}
[this is a ‘‘deductive" version of a conservativity property between languages of grounds defined in d'Aragona 2021a]. The derivation is as follows (we use the theorems discussed in the previous section, and abbreviate the theorem to be proved with $\Pi \xi^{\alpha \vee \beta} \Pi \xi^{\neg \alpha} Th_F$). Let $\Delta^0_1$ be the following derivation

\begin{prooftree}
\AxiomC{$2$}
\noLine
\UnaryInfC{$[\xi^{\alpha \vee \beta} : \alpha \vee \beta \times \xi^{\neg \alpha} : \alpha]$}
\RightLabel{($\times_{E, 2}$)}
\UnaryInfC{$\xi^{\neg \alpha} : \neg \alpha$}
\AxiomC{$3$}
\noLine
\UnaryInfC{$[\xi^{\alpha} : \alpha]$}
\RightLabel{($\times_I$)}
\BinaryInfC{$\xi^{\neg \alpha} : \neg \alpha \times \xi^{\alpha} : \alpha$}
\end{prooftree}
The conclusion of $\Delta^0_1$ can be used as minor premise in an application of ($\supset_E$), whose major premise is obtained from ($\rightarrow^\text{w}$), in such a way as to obtain a term of type $\bot$, so that we can apply the explosion principle to derive $Th_F$ - i.e. we have the following derivation $\Delta_1$

\begin{prooftree}
\AxiomC{$\Delta^0_1$}
\AxiomC{theorem of Section 4.2.2}
\UnaryInfC{$\Pi \xi^{\neg \alpha} \Pi \xi^{\alpha} (\xi^{\neg \alpha} : \neg \alpha \times \xi^{\alpha} : \alpha \supset \ \rightarrow E (\xi^{\neg \alpha}, \xi^{\alpha}) : \bot)$}
\RightLabel{($\Pi_E$)}
\UnaryInfC{$\Pi \xi^{\alpha} (\xi^{\neg \alpha} : \neg \alpha \times \xi^{\alpha} : \alpha \supset \ \rightarrow E (\xi^{\neg \alpha}, \xi^{\alpha}) : \bot)$}
\RightLabel{($\Pi_E$)}
\UnaryInfC{$\xi^{\neg \alpha} : \neg \alpha \times \xi^{\alpha} : \alpha \supset \ \rightarrow E (\xi^{\neg \alpha}, \xi^{\alpha}) : \bot$}
\RightLabel{($\supset_E$)}
\BinaryInfC{$\rightarrow E (\xi^{\neg \alpha}, \xi^{\alpha}) : \bot$}
\RightLabel{$\bot$}
\UnaryInfC{$\bot^G$}
\RightLabel{($\bot^G$)}
\UnaryInfC{$F(\xi^{\alpha \vee \beta}, \xi^{\neg \alpha}) \equiv \vee E \xi^{\alpha} \xi^{\beta} (\xi^{\alpha \vee \beta}, \bot_\beta(\rightarrow E (\xi^{\alpha}, \xi^{\neg \alpha})), \xi^\beta)$}
\end{prooftree}
This constitutes a derivation from the left side-assumptions of an application of type elimination to $\xi^{\alpha \vee \beta} : \alpha \vee \beta$, which is the final goal of our proof. The derivation from the right side-assumptions is instead obtained by applying the equations for $F$ and for $\vee E$ on an arbitrary canonical term of type $\alpha \vee \beta$ and on an arbitrary term of type $\beta$, and then by showing that what one obtains is transitively equivalent to an application of $\vee E$ to $\xi^{\alpha \vee \beta}$, $\rightarrow E (\xi^{\neg \alpha}, \xi^\alpha)$ and $\xi^\beta$, which is in turn provably equivalent to $F(\xi^{\alpha \vee \beta}, \xi^{\neg \alpha})$. This reasoning amounts to four sub-derivations. First, we have the derivation  $\Delta^0_2$

\begin{prooftree}
\AxiomC{$4$}
\noLine
\UnaryInfC{$[\xi^{\alpha \vee \beta} \equiv \vee I (\xi^\beta)]$}
\AxiomC{}
\RightLabel{$\equiv_R$}
\UnaryInfC{$\xi^{\neg \alpha} \equiv \xi^{\neg \alpha}$}
\RightLabel{$\equiv^F$}
\BinaryInfC{$F(\xi^{\alpha \vee \beta}, \xi^{\neg \alpha}) \equiv F(\vee I (\xi^\beta), \xi^{\neg \alpha})$}
\AxiomC{}
\RightLabel{$R_F$}
\UnaryInfC{$F(\vee I (\xi^\beta), \xi^{\neg \alpha}) \equiv \xi^\beta$}
\RightLabel{$\equiv_T$}
\BinaryInfC{$F(\xi^{\alpha \vee \beta}, \xi^{\neg \alpha}) \equiv \xi^\beta$}
\end{prooftree}
Then we have the derivation $\Delta^{1, 1}_2$

\begin{prooftree}
\AxiomC{$5$}
\noLine
\UnaryInfC{$[\xi^{\alpha \vee \beta} \equiv \vee I (\xi^\beta)]$}
\AxiomC{}
\RightLabel{$\equiv_R$}
\UnaryInfC{$\bot_\beta(\rightarrow E (\xi^{\neg \alpha}, \xi^{\alpha})) \equiv \bot_\beta(\rightarrow E (\xi^{\neg \alpha}, \xi^{\alpha}))$}
\AxiomC{}
\RightLabel{$\equiv_R$}
\UnaryInfC{$\xi^\beta \equiv \xi^\beta$}
\RightLabel{$\equiv^\vee_3$}
\TrinaryInfC{$\vee E \xi^\alpha \xi^\beta (\xi^{\alpha \vee \beta}, \bot_\beta(\rightarrow E (\xi^{\neg \alpha}, \xi^{\alpha})), \xi^\beta) \equiv \vee E \xi^\alpha \xi^\beta (\vee I (\xi^\beta), \bot_\beta(\rightarrow E (\xi^{\neg \alpha}, \xi^{\alpha})), \xi^\beta)$}
\end{prooftree}
Then the derivation $\Delta^1_2$

\begin{prooftree}
\AxiomC{$\Delta^{1, 1}_2$}
\AxiomC{}
\RightLabel{$R_{\vee, 2}$}
\UnaryInfC{$\vee E \xi^\alpha \xi^\beta (\vee I (\xi^\beta), \bot_\beta(\rightarrow E (\xi^{\neg \alpha}, \xi^{\alpha})), \xi^\beta) \equiv \xi^\beta$}
\RightLabel{$\equiv_T$}
\BinaryInfC{$\vee E \xi^\alpha \xi^\beta (\xi^{\alpha \vee \beta}, \bot_\beta(\rightarrow E (\xi^{\neg \alpha}, \xi^{\alpha})), \xi^\beta) \equiv \xi^\beta$}
\RightLabel{$\equiv_S$}
\UnaryInfC{$\xi^\beta \equiv \vee E \xi^\alpha \xi^\beta (\xi^{\alpha \vee \beta}, \bot_\beta(\rightarrow E (\xi^{\neg \alpha}, \xi^{\alpha})), \xi^\beta)$}
\end{prooftree}
And finally we have the derivation $\Delta_2$

\begin{prooftree}
\AxiomC{$\Delta^0_2$}
\AxiomC{$\Delta^1_2$}
\RightLabel{$\equiv_T$}
\BinaryInfC{$F(\xi^{\alpha \vee \beta}, \xi^{\neg \alpha}) \equiv \vee E \xi^\alpha \xi^\beta (\xi^{\alpha \vee \beta}, \bot_\beta(\rightarrow E (\xi^{\neg \alpha}, \xi^{\alpha})), \xi^\beta)$}
\end{prooftree}
We can now apply type elimination to $\xi^{\alpha \vee \beta} : \alpha \vee \beta$ to obtain the consequent of $Th_F$, to which we can apply ($\supset_I$) to discharge the assumption $\xi^{\alpha \vee \beta} \times \xi^{\neg \alpha} : \neg \alpha$, and then we can apply twice ($\Pi_I$), respectively binding $\xi^{\neg \alpha}$ and $\xi^{\alpha \vee \beta}$, i.e.

\begin{prooftree}
\AxiomC{$1$}
\noLine
\UnaryInfC{$[\xi^{\alpha \vee \beta} : \alpha \vee \beta \times \xi^{\neg \alpha} : \neg \alpha]$}
\RightLabel{($\times_{E, 1}$)}
\UnaryInfC{$\xi^{\alpha \vee \beta} : \alpha \vee \beta$}
\AxiomC{$\Delta_1$}
\AxiomC{$\Delta_2$}
\RightLabel{$D_\vee$, $3, 4, 5$}
\TrinaryInfC{$F(\xi^{\alpha \vee \beta}, \xi^{\neg \alpha}) \equiv \vee E \xi^\alpha \xi^\beta (\xi^{\alpha \vee \beta}, \bot_\beta(\rightarrow E (\xi^{\neg \alpha}, \xi^{\alpha})), \xi^\beta)$}
\RightLabel{($\supset_I$), $1, 2$}
\UnaryInfC{$Th_F$}
\RightLabel{($\Pi_I$)}
\UnaryInfC{$\Pi \xi^{\neg \alpha} \ Th_F$}
\RightLabel{($\Pi_I$)}
\UnaryInfC{$\Pi \xi^{\alpha \vee \beta} \ \Pi \xi^{\neg \alpha} \ Th_F$}
\end{prooftree}
\medskip

\noindent This property also relates to a completeness issue, that we discuss in the Conclusions of this paper.

\subsection{A theorem about the system}

In this section, we prove a meta-theorem that, as said, can be considered as a kind of \emph{correctness} result for $\texttt{GG}$ - restricted to terms of $\texttt{Gen}$ [it is a ‘‘deductive" version of an analogous theorem proved in d'Aragona 2021a for denotation functions]. The theorem could be given in a more fine-grained form, by taking into account typed-variables with equal type but different indexes. However, we shall omit these details here.

\begin{thm}
 Let $U : \beta \in \texttt{TERM}_\texttt{Gen}$ with $FV^T(U) = \{\xi^{\alpha_1}, ..., \xi^{\alpha_n}\}$. Then $\xi^{\alpha_1} : \alpha_1, ..., \xi^{\alpha_n} : \alpha_n \vdash_{\texttt{GG}} U : \beta$.
\end{thm}

\begin{proof}
 By induction on the complexity of $U$, using the results proved in Section 4.2.2. The induction basis is trivial, i.e. for every $\xi^\alpha$, $\xi^\alpha : \alpha \vdash_{\texttt{GG}} \xi^\alpha : \alpha$, and for every $\delta$ naming a derivation in the atomic system with conclusion $\alpha$, $\vdash_\texttt{GG} \delta : \alpha$ - observe that $\delta$ is a closed term. For the induction step we prove just some more relevant examples.
 
 \begin{itemize}
\item $\vee I (Z) : \alpha_1 \vee \alpha_2$ with $Z : \alpha_i$ ($i = 1, 2$) and $FV^T(\vee I (Z)) = FV^T(Z) \Rightarrow$ by induction hypothesis, there is $\Delta \in \texttt{DER}_\texttt{GG}$ that satisfies the required properties, so that

\begin{prooftree}
\AxiomC{$\Delta$}
\noLine
\UnaryInfC{$Z : \alpha_i$}
\RightLabel{$\vee I$}
\UnaryInfC{$\vee I (Z) : \alpha_1 \vee \alpha_2$}
\end{prooftree}
\medskip

\noindent satisfies the required properties;
      
\item $\vee E \xi^{\beta_1} \xi^{\beta_2} (Z_1, Z_2, Z_3) : \beta_3$ with $Z_1 : \beta_1 \vee \beta_2$, $Z_2 : \beta_3$, $Z_3 : \beta_3$ and

\begin{center}
$FV^T(\vee E \xi^{\beta_1} \xi^{\beta_2} (Z_1, Z_2, Z_3)) = FV^T(Z_1) \cup (FV^T(Z_2) - \{\xi^{\beta_1}\}) \cup (FV^T(Z_3) - \{\xi^{\beta_2}\})$
\end{center}
$\Rightarrow$ by induction hypothesis, there are $\Delta_1, \Delta_2, \Delta_3 \in \texttt{DER}_\texttt{GG}$ that satisfy the required properties, so that, called $\Delta$ the following derivation

\begin{prooftree}
\AxiomC{$\Delta_1$}
\noLine
\UnaryInfC{$Z_1 : \beta_1 \vee \beta_2$}
\AxiomC{$[\xi^{\beta_1} : \beta_1]$}
\noLine
\UnaryInfC{$\Delta_2$}
\noLine
\UnaryInfC{$Z_2 : \beta_3$}
\RightLabel{($\supset_I$)}
\UnaryInfC{$\xi^{\beta_1} : \beta_1 \supset Z_2 : \beta_3$}
\RightLabel{($\Pi_I$)}
\UnaryInfC{$\Pi \xi^{\beta_1} (\xi^{\beta_1} : \beta_1 \supset Z_2 : \beta_3)$}
\RightLabel{($\times_I$)}
\BinaryInfC{$Z : \beta_1 \vee \beta_2 \times \Pi \xi^{\beta_1} (\xi^{\beta_1} : \beta_1 \supset Z_2 : \beta_3)$}
\AxiomC{$[\xi^{\beta_2} : \beta_2]$}
\noLine
\UnaryInfC{$\Delta_3$}
\noLine
\UnaryInfC{$Z : \beta_3$}
\RightLabel{($\supset_I$)}
\UnaryInfC{$\xi^{\beta_2} : \beta_2 \supset Z_3 : \beta_3$}
\RightLabel{($\Pi_I$)}
\UnaryInfC{$\Pi \xi^{\beta_2} (\xi^{\beta_2} : \beta_2 \supset Z_3 : \beta_3)$}
\RightLabel{($\times_I$)}
\BinaryInfC{$(Z : \beta_1 \vee \beta_2 \times\Pi \xi^{\beta_1} (\xi^{\beta_1} : \beta_1 \supset Z_2 : \beta_3)) \times \Pi \xi^{\beta_2} (\xi^{\beta_2} : \beta_2 \supset Z_3 : \beta_3)$}
\end{prooftree}
\medskip

\noindent we have that

\begin{prooftree}
\AxiomC{$\Delta$}
\AxiomC{theorem of Section 4.2.2}
\UnaryInfC{$\Pi \xi^{\beta_1 \vee \beta_2} \Pi \xi^{\beta_3}_1 \Pi \xi^{\beta_3}_2 \ Th_\vee$}
\RightLabel{($\Pi_E$)}
\UnaryInfC{$\Pi \xi^{\beta_3}_1 \Pi \xi^{\beta_3}_2 \ Th_\vee [Z_1/\xi^{\beta_1 \vee \beta_2}]$}
\RightLabel{($\Pi_E$)}
\UnaryInfC{$\Pi \xi^{\beta_3}_2 \ (Th_\vee [Z_1/\xi^{\beta_1 \vee \beta_2}]) [Z_2/\xi^{\beta_3}_1]$}
\RightLabel{($\Pi_E$)}
\UnaryInfC{$((Th_\vee [Z_1/\xi^{\beta_1 \vee \beta_2}]) [Z_2/\xi^\beta_1]) [Z_3/\xi^{\beta_3}_2]$}
\RightLabel{($\supset_E$)}
\BinaryInfC{$\vee E \xi^{\beta_1} \xi^{\beta_2} (Z_1, Z_2, Z_3) : \beta_3$}
\end{prooftree}
\medskip

\noindent satisfies the required properties;

\item $\forall I x (Z) : \forall x \beta(x)$ with $Z : \beta(x)$ and $FV^T(\forall I x (Z)) = FV^T(Z) \Rightarrow$ by induction hypothesis, there is $\Delta \in \texttt{DER}_{\texttt{GG}}$ that satisfies the required properties (recall that, for the hypotheses of the theorem, $\forall I x (Z) \in \texttt{TERM}_{\texttt{Gen}}$, therefore $x$ cannot occur free in $\gamma$ for $\xi^\gamma \in FV^T(Z)$ and hence, again by induction hypothesis, $x$ does not occur free in any undischarged assumption of $\Delta$), so that

\begin{prooftree}
\AxiomC{$\Delta$}
\noLine
\UnaryInfC{$Z : \beta(x)$}
\RightLabel{$\forall I$}
\UnaryInfC{$\forall x I (Z) : \forall x \beta(x)$}
\end{prooftree}
\medskip

\noindent satisfies the required properties;
\item $\forall E (Z) : \beta(t)$ with $Z : \forall x \beta(x)$ and $FV^T(\forall E (Z)) = FV^T(Z) \Rightarrow$ by induction hypothesis, there is $\Delta \in \texttt{DER}_{\texttt{GG}}$ that satisfies the required properties, so that

\begin{prooftree}
\AxiomC{$\Delta$}
\noLine
\UnaryInfC{$Z : \forall x \beta(x)$}
\AxiomC{theorem of Section 4.2.2}
\UnaryInfC{$\Pi \xi^{\forall x \beta(x)} (\xi^{\forall x \beta(x)} : \forall x \beta(x) \supset \forall E (\xi^{\forall x \beta(x)}) : \beta(t))$}
\RightLabel{($\Pi_E$)}
\UnaryInfC{$Z : \forall x \beta(x) \supset \forall E (Z) : \beta(t)$}
\RightLabel{($\supset_E$)}
\BinaryInfC{$\forall E (Z) : \beta(t)$}
\end{prooftree}
\medskip

\noindent satisfies the required properties.
 \end{itemize}
The other cases are in all ways analogous.
\end{proof}

\noindent We remark that the converse of theorem 3 clearly does not hold - e.g.

\begin{prooftree}
\AxiomC{}
\RightLabel{$R_{\wedge, 1}$}
\UnaryInfC{$\wedge_{E, 1}(\wedge I(\delta, \xi^\beta)) \equiv \delta$}
\AxiomC{}
\RightLabel{$\texttt{C}$}
\UnaryInfC{$\delta : \alpha$}
\RightLabel{$\equiv_P$}
\BinaryInfC{$\wedge_{E, 1}(\wedge I(\delta, \xi^\beta)) : \alpha$}
\end{prooftree}
but $FV^T(\wedge_{E, 1}(\wedge I(\delta, \xi^\beta))) = \{\xi^\beta\} \neq \emptyset$. This depends on the fact that the equations ruling the non-primitive operational symbols provide ‘‘syntactic" means of transformations between terms. It is no coincidence that Prawitz defines denotation functions of open terms by appealing to \emph{closed instances} of the operation they denote. One usually says in these cases that terms denote, not their \emph{normal} form, but their \emph{full-evaluated} form, i.e. their form reduced up to open sub-terms [see e.g. Tranchini 2019].

\section{A class of systems}

When Prawitz speaks of a general notion of language of grounding $\Lambda$ over an atomic base $\mathfrak{B}$ with background language $L$, he seems to understand $\Lambda$ as coming with the same structure as $\texttt{Gen}$ over $\mathfrak{B}$, i.e. $\Lambda$ contains:

\begin{itemize}
    \item an alphabet $\texttt{Al}_{\Lambda}$ with individual constants for closed atomic derivations in the atomic system of $\mathfrak{B}$, variables typed on $L$, and a recursive set of primitive and non-primitive operational symbols typed on $L$;
    \item a set of typed terms $\texttt{TERM}_{\Lambda}$, defined inductively according to the types of the operational symbols.
\end{itemize}
Based on this, one can easily provide a general notion of expansion of $\Lambda$, that is, a language of grounding $\Lambda^*$ on $\mathfrak{B}^+$ such that $\texttt{Al}_{\Lambda} \subseteq \texttt{Al}_{\Lambda^*}$, and with the following properties: $\mathfrak{B} \subseteq \mathfrak{B}^+$ or $\Omega_\Lambda \subseteq \Omega_{\Lambda^*}$, where $\Omega_\Lambda$ and $\Omega_{\Lambda^*}$ are the sets of the operational symbols of $\Lambda$ and $\Lambda^*$ respectively - we write $\Lambda \subseteq \Lambda^*$ [see d'Aragona 2021a for a development of this idea].

General notions of system of grounding and of expansion of a system of grounding are now easily obtained. Given a language of grounding $\Lambda$, we build an enriched version $\Lambda^+$ of it as we did for $\texttt{Gen}$ with $\texttt{Gen}^\varepsilon$ - observe that this enriched version \emph{is not} an expansion of $\Lambda$ in the sense defined above - and then a formal system $\Sigma$ over $\Lambda^+$ as we did for $\texttt{GG}$ over $\texttt{Gen}^\varepsilon$. Therefore $\Sigma$ contains:

\begin{itemize}
    \item[($\tau$)] type introductions and type eliminations;
    \item[($\epsilon$)] equivalence rules, i.e.
    \begin{itemize}
        \item[($\epsilon_1$)] standard rules for reflexivity, simmetry and transitivity, plus preservation of denotation;
        \item[($\epsilon_2$)] \textbf{rules permitting to prove that operational symbols applied to equivalent terms yield equivalent terms};
        \item[($\epsilon_3$)] \textbf{equations for computing non-canonical terms};
    \end{itemize}
    \item[($\lambda$)] first-order logical rules.
\end{itemize}
The bold character used for the groups of rules ($\epsilon_2$) and ($\epsilon_3$) is no coincidence. While ($\tau$), ($\epsilon_1$) and ($\lambda$) contain the \emph{same} schematic rules for whatever system of grounding, ($\epsilon_2$) and ($\epsilon_3$) depend on the set of non-primitive operational symbols of $\Lambda$, and can be thus called \emph{characteristic} rules of $\Sigma$. This also means that an expansion $\Sigma^+$ of $\Sigma$ differs from $\Sigma$ only with respect to ($\epsilon_2$) and ($\epsilon_3$). $\Sigma^+$ is a system of grounding over $\Lambda^{*+}$, where $\Lambda^{*+}$ is the enriched version of $\Lambda^*$, for $\Lambda \subseteq \Lambda^*$ - we write $\Sigma \subseteq \Sigma^+$ (observe however that $\Sigma^+$ differs from $\Sigma$ also with respect to the axioms for the individual constants, when the base of $\Lambda^*$ properly contains the base of $\Lambda$).

To make somewhat more precise the class of systems of grounding just introduced, we have to say something about the general form of the characteristic rules of a system of grounding $\Sigma$ over $\Lambda$. As for the rules in the group ($\epsilon_2$), since the \emph{primitive} operational symbols are the same in any language of grounding, the rules concerned are only those saying that, whenever a \emph{non-primitive} operational symbol is applied to equivalent terms, it yields equivalent terms. Therefore, the general form of the characteristic rule is in this case very straightforward. Given $F[\alpha_1, ..., \alpha_n \rhd \beta] \in \texttt{Al}_\Lambda$, binding (possibly empty) sequences of individual and typed-variables $\underline{x}, \underline{\xi}$, and given $T_i, U_i : \alpha_i \in \texttt{TERM}_{\Lambda^+}$ ($1 \leq i \leq n$),

\begin{prooftree}
 \AxiomC{$T_1 \equiv U_1$}
 \AxiomC{$\dots$}
 \AxiomC{$T_n \equiv U_n$}
 \RightLabel{$\equiv^F$}
 \TrinaryInfC{$F \ \underline{x} \ \underline{\xi} \ (T_1, ..., T_n) \equiv F \ \underline{x} \ \underline{\xi} \ (U_1, ..., U_n)$}
\end{prooftree}
As regards the equations for the computation of the non-primitive operational symbols, they have to satisfy some properties, similar to those required by [Prawitz 1973] for the justification procedures on arguments. Given $F[\alpha_1, ..., \alpha_n \rhd \beta] \in \texttt{Al}_\Lambda$, binding (possibly empty) sequences of individual and typed-variables $\underline{x}, \underline{\xi}$, and given $U_i : \alpha_i \in \texttt{TERM}_{\Lambda^+}$ ($1 \leq i \leq n$): (1) each rule for computing $F$ has the form

\begin{prooftree}
 \AxiomC{}
 \RightLabel{$R_F$}
 \UnaryInfC{$F \ \underline{x} \ \underline{\xi} \ (U_1, ..., U_n) \equiv Z$}
\end{prooftree}
for some $Z : \beta \in \texttt{TERM}_{\Lambda^+}$; (2) for each $R_F$ of the kind indicated in (1), it holds that $\star(F \ \underline{x} \ \underline{\xi} \ (U_1, ..., U_n)) \subseteq \star(Z)$, for $\star = FV^I, FV^T, FV^F$; (3) for each $R_F$ of the kind indicated in (1), for every substitution $[\ast/\circ]$ of variables with terms, there is a rule for computing $F$ of the form

\begin{prooftree}
 \AxiomC{}
 \RightLabel{$R'_F$}
 \UnaryInfC{$(F \ \underline{x} \ \underline{\xi} \ (U_1, ..., U_n))[\ast/\circ] \equiv W$}
\end{prooftree}
such that $\vdash_\Sigma W \equiv Z[\ast/\circ]$. (3) is here a linearity condition for $F$ under substitutions of variables. If we indicate a substitution with $\texttt{sub}$, the application of $F$ to given arguments with $F(a)$, and a computation of $F(a)$ with $\texttt{comp}(F(a))$, we have $\vdash_\Sigma \texttt{sub}(\texttt{comp}(F(a))) \equiv \texttt{comp}(\texttt{sub}(F(a))) = \texttt{comp}(F(\texttt{sub}(a)))$. So for example take $R_{\wedge, 1}$. We have the two instances

\begin{prooftree}
 \AxiomC{}
 \RightLabel{$R_{\wedge, 1}$}
 \UnaryInfC{$\wedge_{E, 1}(\wedge I (\xi^\alpha, \xi^\beta)) \equiv \xi^\alpha$}
 \AxiomC{}
 \RightLabel{$R_{\wedge, 1}$}
 \UnaryInfC{$\wedge_{E, 1}(\wedge I (T, \xi^\beta)) \equiv T$}
 \noLine
 \BinaryInfC{}
\end{prooftree}
and - with a little but intuitive abuse of notation -

\begin{center}
    $(\wedge_{E, i}(\wedge I (\xi^\alpha, \xi^\beta)))[T/\xi^\alpha] \equiv \xi^\alpha[T/\xi^\alpha] = T \equiv \wedge_{E, 1}(\wedge I(T, \xi^\beta)) = \wedge_{E, i}((\wedge I(\xi^\alpha, \xi^\beta)[T/\xi^\alpha])$.
\end{center}

We may assume that all the systems of grounding described in this section allow to prove results analogous to those proved within $\texttt{GG}$ in Sections 4.2.1, 4.2.2 and 4.2.3. This means in particular that we may assume that the equations for computing the non-primitive operational symbols are always well-given, and require that a meta-theorem analogous to the one proved about $\texttt{GG}$ in Section 4.3 be provable for each system of our class.

\section{Normalization}

We now prove a normalization result that abstracts from characteristic rules, and hence holds for each of the systems of grounding outlined in Section 5.

Normalization is based, as usual, on the possibility of transforming a derivation with maximal formulas into one without maximal formulas. The idea of the normalization proof is of course that of progressively lowering the complexity of the maximal formulas through appropriate reduction and permutation functions. We deal with reduction functions first, then we outline permutation functions, and finally we introduce some definitions as well as some conventions on the derivations that make the normalization proof easier. Finally, we provide a measure for maximal formulas, and we prove the theorem itself.

\subsection{Reduction functions}

At variance with a standard logical system where only \emph{logically} maximal formulas (LMF) are at issue, in systems of grounding we also have \emph{typing} maximal formulas (TMF), namely, formulas which occur both as conclusions of a type introduction, and as major premises of a type elimination - when there is no need to distinguish, we shall refer to both LMF and TMF as \emph{maximal points}.

Both LMFs and TMFs require reduction functions that eliminate redundant passages. Standard reduction functions for LMFs can be found in [Prawitz 2006]. Reduction functions for TMFs are instead defined below, but basically follow the same intuition as reduction functions for LMFs. It is redundant to prove by type-introduction (possibly under assumptions) that a term $T$ denotes a ground for $\alpha$, and then to apply type elimination to $T : \alpha$ to obtain $A$, since the conditions required by the assumptions of the minor premises $A$ in the type elimination - namely, $T$ is equivalent to some canonical form whose immediate sub-terms denote (operations on) grounds for the immediate sub-formulas of $\alpha$ - are already satisfied by $T$ itself and by the derivations (possibly under assumptions) of the premises of the type-introduction.

However, while standard reduction functions for LMFs remove \emph{logically complex} formulas of the enriched language of grounding $\Lambda^+$, those associated to the TMFs remove \emph{atomic formulas} of $\Lambda^+$. The redundant passage in a TMF, therefore, do not concern the introduction-elimination of a logical constant of $\Lambda^+$, but the introduction-elimination of a primitive operational symbol $s I$ in a term, and of a logical constant $s$ in the formula of the background language for which that term is said to denote a ground. In a sense, the reduction functions for the TMFs can be considered as an ‘‘embedding" into a system $\Sigma$ for $\Lambda^+$ of the so-called \emph{inversion principle} [Prawitz 2006] relative to a proof-system, say $\texttt{S}$, \emph{for the background language} of $\Lambda^+$. The principle says that, given a derivation in $\texttt{S}$ of $\alpha$ from $\Gamma$ whose last inference is an elimination rule the major premise of which is inferred by an introduction, the immediate sub-derivations already ‘‘contain" a derivation of $\alpha$ from $\Gamma$. Now, if we take the derivations of $\texttt{S}$ to be Curry-Howard translated by the terms of $\Lambda$, what the ‘‘embedding" shows is that the ‘‘denotational" properties of the derivations of $\texttt{S}$, expressed by the formulas of $\Lambda^+$, remain stable modulo normalization in $\texttt{S}$. To make just an example, that the ‘‘denotational" properties of a non-normal derivation $\Delta$ in $\texttt{S}$, say

\begin{prooftree}
\AxiomC{$\Delta_1$}
\noLine
\UnaryInfC{$\alpha_1$}
\AxiomC{$\Delta_2$}
\noLine
\UnaryInfC{$\alpha_2$}
\RightLabel{($\wedge_I$)}
\BinaryInfC{$\alpha_1 \wedge \alpha_2$}
\RightLabel{($\wedge_I$)}
\UnaryInfC{$\alpha_i$}
\end{prooftree}
are the same as those enjoyed by

\begin{prooftree}
\AxiomC{$\Delta_i$}
\noLine
\UnaryInfC{$\alpha_i$}
\end{prooftree}
can be seen in $\Sigma$ by showing that the consequences we can draw by applying a type elimination after a type-introduction on a term $T : \alpha_1 \wedge \alpha_2$ that translates the immediate sub-derivation of $\Delta$, and whose immediate sub-terms translate $\Delta_1$ and $\Delta_2$ respectively, are the same as those we would obtain from the premises of the type-introduction alone - plus an application of reflexivity to $T$. This is indeed what the TMF reduction for $\wedge$ does, i.e.

\begin{prooftree}
\AxiomC{$\Delta_1$}
\noLine
\UnaryInfC{$T : \alpha$}
\AxiomC{$\Delta_2$}
\noLine
\UnaryInfC{$U : \beta$}
\RightLabel{$\wedge I$}
\BinaryInfC{$\wedge I (T, U) : \alpha \wedge \beta$}
\AxiomC{$[\wedge I (T, U) \equiv \wedge I (\xi^\alpha, \xi^\beta)]$}
\AxiomC{$[\xi^\alpha : \alpha]$}
\AxiomC{$[\xi^\beta : \beta]$}
\noLine
\TrinaryInfC{}
\noLine
\UnaryInfC{$\Delta_3(\xi^\alpha, \xi^\beta)$}
\noLine
\UnaryInfC{}
\noLine
\UnaryInfC{$A$}
\RightLabel{$D_\wedge$}
\BinaryInfC{$A$}
\end{prooftree}
\medskip

\noindent reduces to

\begin{prooftree}
\AxiomC{}
\RightLabel{$\equiv_R$}
\UnaryInfC{$\wedge I (T, U) \equiv \wedge I (T, U)$}
\AxiomC{$\Delta_1$}
\noLine
\UnaryInfC{$T : \alpha$}
\AxiomC{$\Delta_2$}
\noLine
\UnaryInfC{$U : \beta$}
\noLine
\TrinaryInfC{}
\noLine
\UnaryInfC{$\Delta_3(T/\xi^\alpha, U/\xi^\beta)$}
\noLine
\UnaryInfC{}
\noLine
\UnaryInfC{$A$}
\end{prooftree}
\medskip

\noindent Spelled out in detail, we have replaced the passage giving rise to the TMF $\wedge I(T, U) : \alpha \wedge \beta$ with a derivation obtained from the derivation $\Delta_3(\xi^\alpha, \xi^\beta)$ of the minor premise $A$ of $D_\wedge$ where: (1) the assumption $\wedge I(T, U) \equiv \wedge I(\xi^\alpha, \xi^\beta)$ is replaced by an application to $\wedge I(T, U)$ of the reflexivity axiom for equivalence; (2) the assumption $\xi^\alpha : \alpha$ is replaced by the derivation $\Delta_1$ of $T : \alpha$, and $\xi^\alpha$ is replaced with $T$ throughout $\Delta_3(\xi^\alpha, \xi^\beta)$; (3) the assumption $\xi^\beta : \beta$ is replaced by the derivation $\Delta_2$ of $U : \beta$, and $\xi^\beta$ is replaced with $U$ throughout $\Delta_3(\xi^\alpha, \xi^\beta)$ of the minor premise $A$ of $D_\wedge$.

The TMF reductions for $\vee$ and $\exists$ basically run in the same way, so we shall dispense ourselves with spelling out in detail how they work. As for $\vee$,

\begin{prooftree}
\AxiomC{$\Delta_1$}
\noLine
\UnaryInfC{$T : \alpha_i$}
\RightLabel{$\vee I$}
\UnaryInfC{$\vee I [\alpha_i \rhd \alpha_1 \vee \alpha_2] (T) : \alpha_1 \vee \alpha_2$}
\AxiomC{$[\vee I (T) \equiv \vee I (\xi^{\alpha_1})]$}
\AxiomC{$[\xi^{\alpha_1} : \alpha_1]$}
\noLine
\BinaryInfC{}
\noLine
\UnaryInfC{$\Delta_2(\xi^{\alpha_1})$}
\noLine
\UnaryInfC{}
\noLine
\UnaryInfC{$A$}
\AxiomC{$[\vee I (T) \equiv \vee I (\xi^{\alpha_2})]$}
\AxiomC{$[\xi^{\alpha_2} : \alpha_2]$}
\noLine
\BinaryInfC{}
\noLine
\UnaryInfC{$\Delta_3(\xi^{\alpha_2})$}
\noLine
\UnaryInfC{}
\noLine
\UnaryInfC{$A$}
\RightLabel{$D_\vee$}
\TrinaryInfC{$A$}
\end{prooftree}
\medskip

\noindent reduces to

\begin{prooftree}
\AxiomC{}
\RightLabel{$\equiv_R$}
\UnaryInfC{$\vee I (T) \equiv \vee I (T)$}
\AxiomC{$\Delta_1$}
\noLine
\UnaryInfC{$T : \alpha_i$}
\noLine
\BinaryInfC{}
\noLine
\UnaryInfC{$\Delta_{i + 1}(T/\xi^{\alpha_i})$}
\noLine
\UnaryInfC{}
\noLine
\UnaryInfC{$A$}
\end{prooftree}
\medskip

\noindent As for $\exists$,

\begin{prooftree}
\AxiomC{$\Delta_1$}
\noLine
\UnaryInfC{$T : \alpha(t/x)$}
\RightLabel{$\exists I$}
\UnaryInfC{$\exists I [\alpha(t/x) \rhd \exists x \alpha(x)] (T) : \exists x \alpha(x)$}
\AxiomC{$[\exists I (T) \equiv \exists I (\xi^{\alpha(x)})]$}
\AxiomC{$[\xi^{\alpha(x)} : \alpha(x)]$}
\noLine
\BinaryInfC{}
\noLine
\UnaryInfC{$\Delta_2(x, \xi^{\alpha(x)})$}
\noLine
\UnaryInfC{}
\noLine
\UnaryInfC{$A$}
\RightLabel{$D_\exists$}
\BinaryInfC{$A$}
\end{prooftree}
\medskip

\noindent reduces to

\begin{prooftree}
\AxiomC{}
\RightLabel{$\equiv_R$}
\UnaryInfC{$\exists I (T) \equiv \exists I (T)$}
\AxiomC{$\Delta_1$}
\noLine
\UnaryInfC{$T : \alpha(t/x)$}
\noLine
\BinaryInfC{}
\noLine
\UnaryInfC{$\Delta_2(t/x, T/\xi^{\alpha(x)})$}
\noLine
\UnaryInfC{}
\noLine
\UnaryInfC{$A$}
\end{prooftree}
\medskip

\noindent The TMF reduction for $\rightarrow$ and $\forall$ are more complex, since the derivation produced by the reduction involves some new inferences - but no new assumptions. As for $\rightarrow$,

\begin{prooftree}
\AxiomC{$[\xi^\alpha : \alpha]$}
\noLine
\UnaryInfC{$\Delta_1(\xi^\alpha)$}
\noLine
\UnaryInfC{$T(\xi^\alpha) : \beta$}
\RightLabel{$\rightarrow I$}
\UnaryInfC{$\rightarrow I \xi^\alpha (T(\xi^\alpha)) : \alpha \rightarrow \beta$}
\AxiomC{$[\rightarrow I \xi^\alpha (T(\xi^\alpha)) \equiv \ \rightarrow I \xi^\alpha (\texttt{f}^\beta(\xi^\alpha))]$}
\AxiomC{$[\Pi \xi^\alpha (\xi^\alpha : \alpha \supset \texttt{f}^\beta(\xi^\alpha) : \beta)]$}
\noLine
\BinaryInfC{}
\noLine
\UnaryInfC{$\Delta_2(\texttt{f}^\beta(\xi^\alpha))$}
\noLine
\UnaryInfC{}
\noLine
\UnaryInfC{$A$}
\RightLabel{$D_\rightarrow$}
\BinaryInfC{$A$}
\end{prooftree}
\medskip

\noindent reduces to

\begin{prooftree}
\AxiomC{}
\RightLabel{$\equiv_R$}
\UnaryInfC{$\rightarrow I \xi^\alpha (T(\xi^\alpha)) \equiv \ \rightarrow I \xi^\alpha (T(\xi^\alpha))$}
\AxiomC{$[\xi^\alpha : \alpha]$}
\noLine
\UnaryInfC{$\Delta_1(\xi^\alpha)$}
\noLine
\UnaryInfC{$T(\xi^\alpha) : \beta$}
\RightLabel{($\supset_I$)}
\UnaryInfC{$\xi^\alpha : \alpha \supset T(\xi^\alpha) : \beta$}
\RightLabel{($\Pi_I$)}
\UnaryInfC{$\Pi \xi^\alpha (\xi^\alpha : \alpha \supset T(\xi^\alpha) : \beta)$}
\noLine
\BinaryInfC{}
\noLine
\UnaryInfC{$\Delta_2(T(\xi^\alpha)/\texttt{f}^\beta(\xi^\alpha))$}
\noLine
\UnaryInfC{}
\noLine
\UnaryInfC{$A$}
\end{prooftree}
\medskip

\noindent The passage giving rise to the TMF $\rightarrow I \xi^\alpha (T(\xi^\alpha)) : \alpha \rightarrow \beta$ is replaced with a derivation obtained from the derivation $\Delta_2(\texttt{f}^\beta(\xi^\alpha))$ of the minor premise $A$ of $D_\rightarrow$ where: (1) the assumption $\rightarrow I \xi^\alpha (T(\xi^\alpha)) \equiv \ \rightarrow I \xi^\alpha (\texttt{f}^\beta(\xi^\alpha))$ is replaced by an application to $\rightarrow I \xi^\alpha (T(\xi^\alpha))$ of the reflexivity axiom for equivalence; (2) the assumption $\Pi \xi^\alpha (\xi^\alpha : \alpha \supset \texttt{f}^\beta(\xi^\alpha) : \beta)$ is replaced with the derivation of $\Pi \xi^\alpha (\xi^\alpha : \alpha \supset T(\xi^\alpha) : \beta)$ obtained by applying, first ($\supset_I$) to the immediate sub-derivation $\Delta_1(\xi^\alpha)$ of the premise of the TMF discharging $\xi^\alpha : \alpha$, and thereby obtaining $\xi^\alpha : \alpha \supset T(\xi^\alpha) : \beta$, and then ($\Pi_I$) to $\xi^\alpha : \alpha \supset T(\xi^\alpha) : \beta$, binding $\xi^\alpha$. The term $\texttt{f}^\beta(\xi^\alpha)$ is replaced with $T(\xi^\alpha)$ throughout $\Delta_2(\texttt{f}^\beta(\xi^\alpha))$. Observe that the application of ($\Pi_I$) in the new derivation is licit, thanks to the restriction on $\rightarrow I$. As for $\forall$,

\begin{prooftree}
\AxiomC{$\Delta_1(x)$}
\noLine
\UnaryInfC{$T(x) : \alpha(x)$}
\RightLabel{$\forall I$}
\UnaryInfC{$\forall I x (T(x)) : \forall x \alpha(x)$}
\AxiomC{$[\forall I x (T(x)) \equiv \forall I x (\texttt{h}^{\alpha(x)}(x))]$}
\AxiomC{$[\Pi x (\texttt{h}^{\alpha(x)}(x) : \alpha(x))]$}
\noLine
\BinaryInfC{}
\noLine
\UnaryInfC{$\Delta_2(\texttt{h}^{\alpha(x)}(x))$}
\noLine
\UnaryInfC{}
\noLine
\UnaryInfC{$A$}
\RightLabel{$D_\forall$}
\BinaryInfC{$A$}
\end{prooftree}
\medskip

\noindent reduces to

\begin{prooftree}
\AxiomC{}
\RightLabel{$\equiv_R$}
\UnaryInfC{$\forall I x (T(x)) \equiv \forall I x (T(x))$}
\AxiomC{$\Delta_1(x)$}
\noLine
\UnaryInfC{$T(x) : \alpha(x)$}
\RightLabel{($\Pi_I$)}
\UnaryInfC{$\Pi x (T(x) : \alpha(x))$}
\noLine
\BinaryInfC{}
\noLine
\UnaryInfC{$\Delta_2(T/\texttt{h}^{\alpha(x)}(x))$}
\noLine
\UnaryInfC{}
\noLine
\UnaryInfC{$A$}
\end{prooftree}
\medskip

\noindent In this case, the passage giving rise to the TMF $\forall I x (T(x)) : \forall x \alpha(x)$ is obtained from the derivation $\Delta_2(\texttt{h}^{\alpha(x)}(x))$ of the minor premise $A$ of $D_\forall$ where: (1) the assumption $\forall I x (T(x)) \equiv \forall I x (\texttt{h}^{\alpha(x)}(x))$ is replaced by an application to $\forall I x (T(x)) : \forall x \alpha(x)$ of the reflexivity axiom for equivalence; (2) the assumption $\Pi x(\texttt{h}^{\alpha(x)}(x) : \alpha(x))$ is replaced with the derivation of $\Pi x (T(x) : \alpha(x))$ obtained by applying ($\Pi_I$) to the immediate sub-derivation $\Delta_1(x)$ of the premise of the TMF, binding $x$. The term $\texttt{h}^{\alpha(x)}(x)$ is finally replaced with $T(x)$ throughout $\Delta_2(\texttt{h}^{\alpha(x)}(x))$. Observe that the application of ($\Pi_I$) is licit, since the type-introduction $\forall I$ for $\forall$ and $\Pi$-introduction have to respect the same restriction on $x$.

\subsection{Permutation functions}

Permutation functions are defined like in [Prawitz 2006], along the following pattern. Suppose we have

\begin{prooftree}
\AxiomC{$\Delta_n$}
\noLine
\UnaryInfC{$A_n$}
\AxiomC{$\Delta_2$}
\noLine
\UnaryInfC{$A_2$}
\AxiomC{$\Delta_1$}
\noLine
\UnaryInfC{$A_1$}
\AxiomC{$\Delta^a$}
\noLine
\UnaryInfC{$\vdots$}
\RightLabel{$\texttt{intro}$}
\UnaryInfC{$S$}
\AxiomC{$\Delta^b$}
\noLine
\UnaryInfC{$\vdots$}
\RightLabel{$\star_1$}
\TrinaryInfC{$S$}
\AxiomC{$\Delta^c$}
\noLine
\UnaryInfC{$\vdots$}
\RightLabel{$\star_2$}
\TrinaryInfC{$S$}
\noLine
\UnaryInfC{$\vdots$}
\noLine
\UnaryInfC{$S$}
\AxiomC{$\Delta^d$}
\noLine
\UnaryInfC{$\vdots$}
\RightLabel{$\star_n$}
\TrinaryInfC{$S$}
\AxiomC{$\Delta^e$}
\noLine
\UnaryInfC{$\vdots$}
\AxiomC{$\Delta^f$}
\noLine
\UnaryInfC{$\vdots$}
\RightLabel{$\texttt{elim}$}
\TrinaryInfC{$B$}
\end{prooftree}
where $\texttt{intro}$ is either a type introduction or an introduction for a logical constant $\kappa$ of $\texttt{Gen}^\varepsilon$, and where $\texttt{elim}$ is accordingly either the corresponding type elimination or the corresponding elimination for $\kappa$. Also, $\Delta^b$, $\Delta^c$ and $\Delta^d$ could in turn end with an introduction. $\Delta^b$, $\Delta^c$, $\Delta^d$, $\Delta^e$ and $\Delta^f$ could be ‘‘empty", that is, the rules to which they correspond could have a number of premises lower than that indicated in our general representation. Finally, $\star_i$ ($1 \leq i \leq n$) is ($+_E$), or ($\mathfrak{E}_E$), or a type elimination rule. $S$ can be considered as a maximal point, although its occurrence as conclusion of an introduction is very far from its occurrence as major premise of an elimination. In order ‘‘to shorten" this distance, we ‘‘bring up" the last elimination rule within the derivations of the minor premises of the various $\star_i$ ($1 \leq i \leq n$). More specifically, the procedure is applied $n$ times, and the first passage yields

\begin{small}
\begin{prooftree}
\AxiomC{$\Delta_n$}
\noLine
\UnaryInfC{$A_n$}
\AxiomC{$\Delta_2$}
\noLine
\UnaryInfC{$A_2$}
\AxiomC{$\Delta_1$}
\noLine
\UnaryInfC{$A_1$}
\AxiomC{$\Delta^a$}
\noLine
\UnaryInfC{$\vdots$}
\RightLabel{$\texttt{intro}$}
\UnaryInfC{$S$}
\AxiomC{$\Delta^b$}
\noLine
\UnaryInfC{$\vdots$}
\RightLabel{$\star_1$}
\TrinaryInfC{$S$}
\AxiomC{$\Delta^c$}
\noLine
\UnaryInfC{$\vdots$}
\RightLabel{$\star_2$}
\TrinaryInfC{$S$}
\noLine
\UnaryInfC{$\vdots$}
\noLine
\UnaryInfC{$S$}
\AxiomC{$\Delta^e$}
\noLine
\UnaryInfC{$\vdots$}
\AxiomC{$\Delta^f$}
\noLine
\UnaryInfC{$\vdots$}
\RightLabel{$\texttt{elim}$}
\TrinaryInfC{$B$}
\AxiomC{$\Delta^{d^*}$}
\noLine
\UnaryInfC{$\vdots$}
\RightLabel{$\star_n$}
\TrinaryInfC{$B$}
\end{prooftree}
\end{small}
where $\Delta^{d^*}$ is obtained applying

\begin{prooftree}
\AxiomC{$S$}
\AxiomC{$\Delta^e$}
\noLine
\UnaryInfC{$\vdots$}
\AxiomC{$\Delta^f$}
\noLine
\UnaryInfC{$\vdots$}
\RightLabel{$\texttt{elim}$}
\TrinaryInfC{$B$}
\end{prooftree}
to the conclusion of $\Delta^d$. After $n - 1$ passages we have

\begin{small}
\begin{prooftree}
\AxiomC{$\Delta_n$}
\noLine
\UnaryInfC{$A_n$}
\AxiomC{$\Delta_2$}
\noLine
\UnaryInfC{$A_2$}
\AxiomC{$\Delta_1$}
\noLine
\UnaryInfC{$A_1$}
\AxiomC{$\Delta^a$}
\noLine
\UnaryInfC{$\vdots$}
\RightLabel{$\texttt{intro}$}
\UnaryInfC{$S$}
\AxiomC{$\Delta^b$}
\noLine
\UnaryInfC{$\vdots$}
\RightLabel{$\star_1$}
\TrinaryInfC{$S$}
\AxiomC{$\Delta^e$}
\noLine
\UnaryInfC{$\vdots$}
\AxiomC{$\Delta^f$}
\noLine
\UnaryInfC{$\vdots$}
\RightLabel{$\texttt{elim}$}
\TrinaryInfC{$B$}
\AxiomC{$\Delta^{c^*}$}
\noLine
\UnaryInfC{$\vdots$}
\RightLabel{$\star_2$}
\TrinaryInfC{$B$}
\noLine
\UnaryInfC{$\vdots$}
\noLine
\UnaryInfC{$B$}
\AxiomC{$\Delta^{d^*}$}
\noLine
\UnaryInfC{$\vdots$}
\RightLabel{$\star_n$}
\TrinaryInfC{$B$}
\end{prooftree}
\end{small}
where $\Delta^{c^*}$ is obtained from $\Delta^c$ as indicated for $\Delta^{d^*}$. The last passage finally returns

\begin{small}
\begin{prooftree}
\AxiomC{$\Delta_n$}
\noLine
\UnaryInfC{$A_n$}
\AxiomC{$\Delta_2$}
\noLine
\UnaryInfC{$A_2$}
\AxiomC{$\Delta_1$}
\noLine
\UnaryInfC{$A_1$}
\AxiomC{$\Delta^a$}
\noLine
\UnaryInfC{$\vdots$}
\RightLabel{$\texttt{intro}$}
\UnaryInfC{$S$}
\AxiomC{$\Delta^e$}
\noLine
\UnaryInfC{$\vdots$}
\AxiomC{$\Delta^f$}
\noLine
\UnaryInfC{$\vdots$}
\RightLabel{$\texttt{elim}$}
\TrinaryInfC{$B$}
\AxiomC{$\Delta^{b^*}$}
\noLine
\UnaryInfC{$\vdots$}
\RightLabel{$\star_1$}
\TrinaryInfC{$B$}
\AxiomC{$\Delta^{c^*}$}
\noLine
\UnaryInfC{$\vdots$}
\RightLabel{$\star_2$}
\TrinaryInfC{$B$}
\noLine
\UnaryInfC{$\vdots$}
\noLine
\UnaryInfC{$B$}
\AxiomC{$\Delta^{d^*}$}
\noLine
\UnaryInfC{$\vdots$}
\RightLabel{$\star_n$}
\TrinaryInfC{$B$}
\end{prooftree}
\end{small}
where $\Delta^{b^*}$ is obtained in the same way as $\Delta^{d^*}$ and $\Delta^{c^*}$. As a concrete example (with typing rules),

\begin{prooftree}
\AxiomC{$\Delta_1$}
\noLine
\UnaryInfC{$T : \alpha_1 \vee \alpha_2$}
\AxiomC{$[T \equiv \vee I (\xi^{\alpha_1})]$}
\AxiomC{$[\xi^{\alpha_1} : \alpha_1]$}
\noLine
\BinaryInfC{$\Delta_2$}
\noLine
\UnaryInfC{$A$}
\AxiomC{$[T \equiv \vee I (\xi^{\alpha_2})]$}
\AxiomC{$[\xi^{\alpha_2} : \alpha_2]$}
\noLine
\BinaryInfC{$\Delta_3$}
\noLine
\UnaryInfC{$A$}
\RightLabel{$D_\vee$}
\TrinaryInfC{$A$}
\AxiomC{$\Delta_4$}
\AxiomC{$\Delta_5$}
\TrinaryInfC{$B$}
\noLine
\UnaryInfC{$\Delta_6$}
\end{prooftree}
\medskip

\noindent permutes into

\begin{prooftree}
\AxiomC{$\Delta_1$}
\noLine
\UnaryInfC{$T : \alpha_1 \vee \alpha_2$}
\AxiomC{$\Delta^*_{2}$}
\AxiomC{$\Delta^*_{3}$}
\RightLabel{$D_\vee$}
\TrinaryInfC{$B$}
\noLine
\UnaryInfC{$\Delta_6$}
\end{prooftree}
\medskip

\noindent where $\Delta^*_{i + 1}$ ($i = 1, 2$) is

\begin{prooftree}
\AxiomC{$[T \equiv \vee I (\xi^{\alpha_i})]$}
\AxiomC{$[\xi^{\alpha_i} : \alpha_i]$}
\noLine
\BinaryInfC{$\Delta_{i + 1}$}
\noLine
\UnaryInfC{$A$}
\AxiomC{$\Delta_4$}
\AxiomC{$\Delta_5$}
\TrinaryInfC{$B$}
\end{prooftree}
\medskip

\noindent On the other typing rules the definition is analogous, while it is standard for ($+_E$) and ($\mathfrak{E}_E$).

\subsection{Cut-segments and conventions on derivations}

\begin{defn}
 A \emph{cut-segment} in $\Delta$ is a sequence $\sigma$ of $n$ occurrences of a formula $A$ in $\Delta$ such that:
 
 \begin{itemize}
     \item the first occurrence of $A$ is the conclusion of an application of a type introduction rule $s I$ for a logical constant $s$ of the background language, or it is the conclusion of an introduction rule ($\kappa_I$) for a logical constant $\kappa$ of the (enriched) language of grounding;
     \item for every $1 \leq i \leq n$, the $i + 1$-th occurrence of $A$ is the conclusion of an application of ($+_E$), or of ($\mathfrak{E}_E$), or of a type elimination rule;
     \item the $n$-th occurrence of $A$ is the major premise of an application of $D_s$ if the fist occurrence of $A$ in $\sigma$ was a conclusion of an application of $s I$, and the major premise of an application of ($\kappa_E$) if the first occurrence of $\sigma$ was the conclusion of an application of ($\kappa_I$).
 \end{itemize}
\end{defn}
\noindent Given a cut-segment $\sigma = A_1, ..., A_n$, we say that $\sigma$ has \emph{length} $n$. We also say that $A$ is \emph{the formula} of $\sigma$, and that $A_i$ ($1 \leq i \leq n$) is \emph{an occurrence} of $A$ in $\sigma$ - of course, an occurrence in $\sigma$ whose formula is $A$ is an occurrence of $A$ in $\Delta$. Clearly, if $n = 1$, $\sigma$ is a maximal point. Given two cut-segments $\sigma_1, \sigma_2$, we say that $\sigma_1$ is \emph{disjoint} from $\sigma_2$ if, and only if, none of the occurrences in $\sigma_1$ is also an occurrence in $\sigma_2$.

The following conventions cause no loss of generality thanks to well-know results [Prawitz 2006, Cellucci 1978].

\begin{cnv}
 For every application of ($\bot^G$) in $\Delta$ with conclusion $A$, $A \in \texttt{ATOM}_{\Lambda^+}$.
\end{cnv}

\begin{cnv}
 (1) Free and bound individual and typed-variables are all distinct in $\Delta$; (2) proper and non-proper variables are all distinct in $\Delta$, and each proper variable occurs at most once.
\end{cnv}

\subsection{Measure for maximal points}

Maximal points undergo a double measure. We now give this measure, and then we explain why, at variance with what happens in normalization for simply logical system, the measure has to be double. In what follows, we use $\texttt{max} \colon \mathbb{N}^2 \to \mathbb{N}$, which is the standard function

\begin{center}
    $\texttt{max}(n, m) = \begin{cases} n & \text{if} \ m < n \\ m & \text{if} \ n \leq m\end{cases}$
\end{center}
and, put $\mathbb{F} = \{S \ | \ S \subset \mathbb{N} \ \text{and} \ S \ \text{finite}\}$, we use the function $\texttt{MAX} \colon \mathbb{F} \to \mathbb{N}$ such that

\begin{center}
    $\texttt{MAX}(S) = \begin{cases} 0 & \text{if} \ S = \emptyset \\ n \in S & \text{such that, for every} \ m \in S, m \leq n \end{cases}$
\end{center}
$\texttt{MAX}$ can be extended in a natural way to set of pairs of integers ordered according to an alphabetic order. The \emph{measure} of $\alpha \in \texttt{FORM}_L$ is a function $k^1 \colon \texttt{FORM}_L \to \mathbb{N}$ defined in a standard way. Given $A \in \texttt{FORM}_{\Lambda^+}$, we put

\begin{center}
$T_A = \{k^1(\alpha) \ | \ U : \alpha \ \text{sub-formula of} \ A\}$
\end{center}

\begin{defn}
The \emph{type-measure} of $A \in \texttt{FORM}_{\Lambda^+}$ is a function $\tau \colon \texttt{FORM}_{\Lambda^+} \to \mathbb{N}$ such that

\begin{center}
$\tau(A) = \texttt{MAX}(T_A)$.
\end{center}
\end{defn}

\begin{defn}
The \emph{logical measure} of $A \in \texttt{FORM}_{\Lambda^+}$ is a function $k^2 \colon \texttt{FORM}_{\Lambda^+} \to \mathbb{N}$ inductively defined as follows:

\begin{itemize}
\item $A \in \texttt{ATOM}_{\Lambda^+} \Rightarrow k^2(A) = 0$
\item $A = B \square C \rightarrow k^2(A) = \texttt{max}(k^2(B), k^2(C)) + 1$ ($\square = \times, +, \supset$)
\item $A = \square \ \nu \ B \Rightarrow k^2(A) = k^2(B) + 1$ ($\square = \Pi, \mathfrak{E}$)
\end{itemize}
\end{defn}

\begin{defn}
The \emph{measure} of $A \in \texttt{FORM}_{\Lambda^+}$ is a function $\mu \colon \texttt{FORM}_\Lambda \to \mathbb{N}^2$ such that

\begin{center}
$\mu(A) = (\tau(A), k^2(A))$.
\end{center}
\end{defn}
\noindent Thus, the measure of $A \in \texttt{FORM}_{\Lambda^+}$ is a pair whose first element is the measure of the most complex $\alpha \in \texttt{FORM}_L$ occurring in sub-formulas of $A$ of the form $... : ---$, and whose second element is the standard measure of the logical complexity of $A$. As a concrete example, let us compute the complexity of

\begin{center}
   $A = \Pi \xi^{\alpha \vee \beta} (\xi^{\alpha \vee \beta} : \alpha \vee \beta \supset \mathfrak{E} \xi^\alpha \mathfrak{\xi^\beta}((\xi^{\alpha \vee \beta} \equiv \vee I (\xi^\alpha) \times \xi^\alpha : \alpha) + (\xi^{\alpha \vee \beta} \equiv \vee I (\xi^\beta) \times \xi^\beta : \beta)))$.
\end{center}
Assume $k^1(\alpha) = 3$ and $k^1(\beta) = 5$. Thus, $k^1(\alpha \vee \beta) = \texttt{max}(k^1(\alpha), k^1(\beta)) + 1 = \texttt{max}(3, 5) + 1 = 5 + 1 = 6$. Our $A$ contains three atomic sub-formulas of the form $... : ---$, i.e.

\begin{center}
    $\xi^{\alpha \vee \beta} : \alpha \vee \beta$, $\xi^\alpha : \alpha$ and $\xi^\beta : \beta$.
\end{center}
Thus, $T_A = \{3, 5, 6\}$, so that $\tau(A) = \texttt{MAX}(T_A) = 6$. We have now to compute the logical complexity of $A$, that is - not all passages are displayed -

\begin{itemize}
    \item[] $k^2(\Pi \xi^{\alpha \vee \beta} (\xi^{\alpha \vee \beta} : \alpha \vee \beta \supset \mathfrak{E} \xi^\alpha \mathfrak{\xi^\beta}((\xi^{\alpha \vee \beta} \equiv \vee I (\xi^\alpha) \times \xi^\alpha : \alpha) + (\xi^{\alpha \vee \beta} \equiv \vee I (\xi^\beta) \times \xi^\beta : \beta)))) =$
    \item[] $= \texttt{max}(k^2(\xi^{\alpha \vee \beta} : \alpha \vee \beta), k^2(\mathfrak{E} \xi^\alpha \mathfrak{\xi^\beta}((\xi^{\alpha \vee \beta} \equiv \vee I (\xi^\alpha) \times \xi^\alpha : \alpha) + (\xi^{\alpha \vee \beta} \equiv \vee I (\xi^\beta) \times \xi^\beta : \beta)))) + 1 =$
    \item[] $= \texttt{max}(0, k^2((\xi^{\alpha \vee \beta} \equiv \vee I (\xi^\alpha) \times \xi^\alpha : \alpha) + (\xi^{\alpha \vee \beta} \equiv \vee I (\xi^\beta) \times \xi^\beta : \beta)) + 2) + 1 =$
    \item[] $= \texttt{max}(0, \texttt{max}(1, 1) + 1) + 2) + 1 = 5$.
\end{itemize}
So in conclusion we have $\mu(A) = (\tau(A), k^2(A)) = (6, 5)$. Once the measure has been set, on the formulas of $\Lambda^+$ it is possible to define a strict order relation. Since the measure is double, the order is of an ‘‘alphabetic" kind, that is, for every $A, B \in \texttt{FORM}_{\Lambda^+}$,

\begin{center}
$\mu(A) < \mu(B) \Leftrightarrow \begin{cases} \tau(A) < \tau(B) & \text{or} \\ k^2(A) < k^2(B)  & \text{for} \ \tau(A) = \tau(B) \end{cases}$
\end{center}
The reasons for adopting the double measure are essentially two. First, observe that the measure of an atomic formula of the form $T : \alpha$ \emph{is not} $0$, but $(k^1(\alpha), 0)$. This is important for, as said, our strategy for proving normalization will be based on the idea of progressively lowering the complexity of the maximal points. Since maximal points may well be TMFs, i.e. they may have the form $T : \alpha$, we must hence be able to set a non-null measure of them for the lowering process to go through - observe that $\alpha$ cannot be atomic, so the measure is greater than $(0, 0)$ also in these cases. Second, reductions for typing rules may generate new TMFs or LMFs. In the specific case of new LMFs, we must ensure that they have a lower measure than the eliminated TMFs - that the measure is lower in the case of new TMFs depends on the fact that the new TMFs have as type formulas which are sub-formulas of the type of the starting TMFs. Suppose we are in the following situation:

\begin{prooftree}
\AxiomC{$[\xi^\alpha : \alpha]$}
\noLine
\UnaryInfC{$\Delta_1$}
\noLine
\UnaryInfC{$T(\xi^\alpha) : \beta$}
\RightLabel{$\rightarrow I$}
\UnaryInfC{$\rightarrow I \xi^\alpha (T(\xi^\alpha)) : \alpha \rightarrow \beta$}
\AxiomC{$[\rightarrow I \xi^\alpha (T(\xi^\alpha)) \equiv \ \rightarrow I \xi^\alpha (\texttt{f}^\beta(\xi^\alpha))]$}
\AxiomC{$[\Pi \xi^\alpha (\xi^\alpha : \alpha \supset \texttt{f}^\beta(\xi^\alpha) : \beta)]$}
\RightLabel{($\Pi_E$)}
\UnaryInfC{$U : \alpha \supset \texttt{f}^\beta(U) : \beta$}
\noLine
\BinaryInfC{$\Delta_2(\texttt{f}^\beta)$}
\noLine
\UnaryInfC{$A$}
\RightLabel{$D_\rightarrow$}
\BinaryInfC{$A$}
\noLine
\UnaryInfC{$\Delta_3$}
\end{prooftree}
\medskip

\noindent The reduction for $\rightarrow$ of Section 6.1 gives

\begin{prooftree}
\AxiomC{}
\RightLabel{$\equiv_R$}
\UnaryInfC{$\rightarrow I \xi^\alpha (T(\xi^\alpha)) \equiv \ \rightarrow I \xi^\alpha (T(\xi^\alpha))$}
\AxiomC{$[\xi^\alpha : \alpha]$}
\noLine
\UnaryInfC{$\Delta_1$}
\noLine
\UnaryInfC{$T(\xi^\alpha) : \beta$}
\RightLabel{($\supset_I$)}
\UnaryInfC{$\xi^\alpha : \alpha \supset T(\xi^\alpha) : \beta$}
\RightLabel{($\Pi_I$)}
\UnaryInfC{$\Pi \xi^\alpha (\xi^\alpha : \alpha \supset T(\xi^\alpha) : \beta)$}
\RightLabel{($\Pi_E$)}
\UnaryInfC{$U : \alpha \supset T(U) : \beta$}
\noLine
\BinaryInfC{$\Delta_2(T/\texttt{f}^\beta)$}
\noLine
\UnaryInfC{$A$}
\noLine
\UnaryInfC{$\Delta_3$}
\end{prooftree}
with a new LMF, namely, $\Pi \xi^\alpha (\xi^\alpha : \alpha \supset T(\xi^\alpha) : \beta)$. Now note that $\mu(\rightarrow I \xi^\alpha (T(\xi^\alpha)) : \alpha \rightarrow \beta) = (k^1(\alpha \rightarrow \beta), 0)$, and $\mu(\Pi \xi^\alpha (\xi^\alpha : \alpha \supset T(\xi^\alpha) : \beta)) = (\texttt{max}(k^1(\alpha), k^1(\beta)), 2)$. Since $k^1(\alpha \rightarrow \beta) = \texttt{max}(k^1(\alpha), k^1(\beta)) + 1$, we have $\texttt{max}(k^1(\alpha), k^1(\beta)) < k^1(\alpha \rightarrow \beta)$, and hence $\mu(\Pi \xi^\alpha (\xi^\alpha : \alpha \supset T(\xi^\alpha) : \beta)) < \mu(\rightarrow I \xi^\alpha (T(\xi^\alpha)) : \alpha \rightarrow \beta)$ by the alphabetic order.

As happens in Prawitz’s standard normalization theory, also the reductions associated with logical rules could produce new TMFs or new LMFs. In these cases, however, the following applies. Let $A$ be an LMF, and let $B$ be a new TMF or LMF produced by the reduction on $A$; since $B$ is a sub-formula of $A$, we have that $\tau(B) \leq \tau(A)$, and $k^2(B) < k^2(A)$. Therefore, $\mu(B) < \mu(A)$.

\subsection{Normalization theorem and the form of normal derivations}

\begin{defn}
Given $\Delta$, let $\sigma$ be a cut-segment in $\Delta$ the formula of which is $A$. The \emph{measure} of $\sigma$ - written $\mu(\sigma)$ - is $\mu(A)$.
\end{defn}
\noindent To every $\Delta$ we associate the set

\begin{center}
$M_\Delta = \{\mu(\sigma) \ | \ \sigma \ \text{cut-segment of} \ \Delta\}$.
\end{center}

\begin{defn}
The \emph{degree} of $\Delta$ is the value of a function $\delta \colon \texttt{DER}_{\Sigma} \to \mathbb{N}^2 \times \mathbb{N}$ such that

\begin{center}
$\delta(\Delta) = (\texttt{MAX}(M_\Delta), n)$
\end{center}
where $n$ is the sum of the lengths of the cut-segments $\sigma$ of $\Delta$ such that $\mu(\sigma) = \texttt{MAX}(M_\Delta)$.
\end{defn}
\noindent Also on derivations it is possible to define a strict order relation. It is again of an ‘‘alphabetic" kind, that is, for every $\Delta_1, \Delta_2$ such that $\delta(\Delta_i) = (\texttt{MAX}(M_{\Delta_i}), n_i)$ ($i = 1, 2$),

\begin{center}
$\delta(\Delta_1) < \delta(\Delta_2) \Leftrightarrow \begin{cases} \texttt{MAX}(M_{\Delta_1}) < \texttt{MAX}(M_{\Delta_2}) & \text{or} \\ n_1 < n_2 & \text{for} \ \texttt{MAX}(M_{\Delta_1}) = \texttt{MAX}(M_{\Delta_2}) \end{cases}$
\end{center}

\begin{defn}
 $\Delta$ is \emph{normal} iff $\delta(\Delta) = ((0,0), 0)$. Otherwise, $\Delta$ is \emph{non-normal}.
\end{defn}

\begin{defn}
$\Delta_a$ \emph{immediately reduces} to $\Delta_b$ - written $\Delta_a \succeq \Delta_b$ - iff $\Delta_a = \Delta_b$, or $\Delta_b$ can be obtained from $\Delta_a$ by applying one of the reduction functions for TMFs, or one of the reduction functions for LMFs, or a permutation function. Moreover, $\Delta_a$ \emph{reduces} to $\Delta_b$ - written $\Delta_a \succ \Delta_b$ - iff there is a sequence $\Delta_1, ..., \Delta_n$ with $\Delta_1 = \Delta_a$, $\Delta_n = \Delta_b$, and $\Delta_i \succeq \Delta_{i + 1}$ for every $1 \leq i \leq n$.
\end{defn}

\noindent We take for granted the notion of a formula occurrence \emph{standing} (\emph{immediately}) \emph{above} (respectively, \emph{below}) another formula occurrence, and the notion of \emph{side-connected} occurrences of formulas - basically, formulas that occur as premises of the same inference. The \emph{upper edge} of a cut-segment $\sigma$ whose formula is $A$ is the occurrence of $A$ in $\sigma$ that stands above all the other occurrences of $A$ in $\sigma$ - the \emph{lower edge} of $\sigma$ is accordingly the occurrence of $A$ in $\sigma$ that stands below all the other occurrences of $A$ in $\sigma$.

As said, the normalization strategy is based on the idea of progressively lowering the degree of a non-normal derivation. A standard way for doing this - that we shall also follow here - is that of finding a topmost rightmost maximal cut-segment $\sigma$ in the derivation. This means that $\sigma$: (1) has measure equal to the first element of the degree of the derivation; (2) there is no cut-segment of equal measure (2.1) whose lower edge stands above the upper edge of $\sigma$ or (2.2) in which there is an occurrence side-connected to the lower edge of $\sigma$ or (2.3) whose lower edge stands above a formula side-connected with the lower edge of $\sigma$. We can thereby grant that the application of a reduction or permutation function on the lower edge of $\sigma$ actually decreases the degree of the derivation, i.e. either all the cut-segments of the resulting derivation are of a lower measure, or the overall length of the cut-segments of highest measure in the resulting derivation is smaller. The existence of a topmost cut-segment in a non-normal derivation can be proved by easy induction on the length of the derivation.

\begin{prp}
If $\Delta$ is a non-normal derivation, there is a cut-segment $\sigma$ in $\Delta$ such that $\mu(\sigma) = \texttt{MAX}(M_\Delta)$ and such that, for no cut-segment $\sigma^*$, the lower edge of $\sigma^*$ stands above the upper edge of $\sigma$ and $\mu(\sigma^*) = \mu(\sigma)$.
\end{prp}

\noindent Then we show that, if a topmost maximal cut-segment $\sigma$ is not also a rightmost maximal cut-segment, then there is a topmost maximal cut-segment disjoint from $\sigma$, and that, given a sequence $\sigma_1, ..., \sigma_n$ of topmost maximal cut-segments which are not also rightmost maximal cut-segments, then there is a topmost maximal cut-segment disjoint from each $\sigma_i$ ($i \leq n$) - we abbreviate ‘‘topmost maximal cut-segment" with tmcs.

\begin{prp}
Given a non-normal derivation $\Delta$, and given a tmcs $\sigma^1$ of $\Delta$, suppose there is a cut-segment $\sigma^*$ in $\Delta$ such that $\mu(\sigma^*) = \mu(\sigma^1)$ and such that the lower edge of $\sigma^1$ is side-connected either with an occurrence in $\sigma^*$ or with a formula occurrence in $\Delta$ that stands below the lower edge of $\sigma^*$. Then, there is a tmcs $\sigma^2$ of $\Delta$ disjoint from $\sigma^1$. Moreover, let $\sigma^1, \sigma^2, \sigma^3, ..., \sigma^n$ be a sequence of pairwise disjoint tmcs of $\Delta$ such that, for every $i < n$, $\sigma^{i + 1}$ stands with $\sigma^i$ in the same relation as the one occurring between $\sigma^*$ and $\sigma^1$ above, and suppose there is a cut-segment $\sigma^{**}$ in $\Delta$ that stands with $\sigma^n$ in the same relation as the one occurring between $\sigma^*$ and $\sigma^1$ above. Then, there is a tmcs $\sigma^{n + 1}$ of $\Delta$ such that, for every $\sigma^i$ ($1 \leq i \leq n$), $\sigma^{n + 1}$ is disjoint from $\sigma^i$.
\end{prp}

\begin{proof}
If $\sigma^*$ is a tmcs of $\Delta$, we can put $\sigma^2 = \sigma^*$ since $\sigma^*$ is disjoint from $\sigma^1$. If $\sigma^*$ is not a tmcs of $\Delta$, then the sub-derivation $\Delta^*$ of $\Delta$ whose conclusion is the upper edge of $\sigma^*$ is non-normal, and we can apply proposition 14 to find a tmcs $\sigma^2$ of $\Delta^*$. Since $\mu(\sigma^2) = \mu(\sigma^*) = \mu(\sigma^1)$, $\sigma^2$ is also a tmcs of $\Delta$. Finally, $\sigma^2$ is disjoint from $\sigma^1$. This proves the first part of the proposition. As for the second, since $\sigma^{**}$ stands with $\sigma^n$ in the same relation as the one one occurring between $\sigma^*$ and $\sigma^1$ above, with the same reasoning as above we can prove that there is a tmcs $\sigma^{n + 1}$ of $\Delta$ disjoint from $\sigma^n$. But now observe that also $\sigma^{n + 1}$ stands with $\sigma^n$ in the same relation as $\sigma^*$ and $\sigma^1$ above, and since we have assumed that this holds for all the pairs $\sigma^i, \sigma^{i + 1}$ of our sequence of pairwise disjoint tmcs of $\Delta$, we can conclude that

\begin{itemize}
    \item[($*$)] for every $i < n + 1$, the lower edge of $\sigma^i$ is side-connected either with an occurrence in $\sigma^{i + 1}$ or with a formula occurrence in $\Delta$ that stands below the lower edge of $\sigma^{i + 1}$ - in other words, $\sigma^{i + 1}$ is to the right of $\sigma^i$.
\end{itemize}
Suppose now that there is $\sigma^i$ not disjoint from $\sigma^{n + 1}$ for some $i < n - 1$. Then, the lower edge of $\sigma^i$ and the lower edge of $\sigma^{n + 1}$ are the same formula occurrence in $\Delta$. So the lower edge of $\sigma^{n + 1}$ is side-connected either with an occurrence in $\sigma^{i + 1}$ or with a formula occurrence in $\Delta$ that stands below the lower edge of $\sigma^{i + 1}$. But this is clearly impossible because of ($*$).
\end{proof}

\noindent Finally, we prove the existence of a topmost rightmost cut-segment in a non-normal derivation.

\begin{prp}
If $\Delta$ is a non-normal derivation, then there is a tmcs $\sigma$ of $\Delta$ such that, for no cut-segment $\sigma^*$ in $\Delta$, $\mu(\sigma^*) = \mu(\sigma)$ and the lower edge of $\sigma$ is side-connected with an occurrence in $\sigma^*$ or with a formula occurrence in $\Delta$ that stands below the lower edge of $\sigma^*$.
\end{prp}

\begin{proof}
The existence of a tmcs in $\Delta$ is granted by proposition 14. Moreover, if no tcms satisfied the required property, by proposition 15 $\Delta$ would contain an infinite number of tmcs, which is clearly impossible.
\end{proof}

\noindent Now, by choosing a topmost rightmost maximal cut-segment - abbreviated with trmcs - we can prove the following fundamental theorem.

\begin{thm}
Let $\Delta$ be a non-normal derivation. Then, there is a derivation $\Delta^*$ such that $\Delta \succeq \Delta^*$ and $\delta(\Delta^*) < \delta(\Delta)$. 
\end{thm}

\begin{proof}
Let us choose a trmcs $\sigma$ of $\Delta$ whose formula is $A$. Suppose that $\sigma$ has length greater than $1$. Then the lower edge of $\sigma$ is the conclusion either of an elimination rule for a logical constant of the (enriched) language of grounding, or of a type elimination rule. To give just an example of each case, we may have either

\begin{prooftree}
\AxiomC{$\Delta_1$}
\noLine
\UnaryInfC{$B + C$}
\AxiomC{$\Delta_2$}
\RightLabel{$\texttt{intro}$}
\UnaryInfC{$\textbf{A}$}
\noLine
\UnaryInfC{$\Delta_3$}
\noLine
\UnaryInfC{$\textbf{A}$}
\AxiomC{$\Delta_4$}
\noLine
\UnaryInfC{$A$}
\RightLabel{($+_E$)}
\TrinaryInfC{$\textbf{A}$}
\AxiomC{$\Delta_5$}
\AxiomC{$\Delta_6$}
\RightLabel{$\texttt{elim}$}
\TrinaryInfC{$D$}
\end{prooftree}
- where the bold character indicates the chosen trmcs, and where $\Delta_5$ or $\Delta_6$ may be empty - or we may have

\begin{prooftree}
\AxiomC{$\Delta_1$}
\noLine
\UnaryInfC{$T : \exists x \alpha(x)$}
\AxiomC{$\Delta_2$}
\RightLabel{$\texttt{intro}$}
\UnaryInfC{$A$}
\noLine
\UnaryInfC{$\Delta_3$}
\noLine
\UnaryInfC{$A$}
\RightLabel{$D_\exists$}
\BinaryInfC{$A$}
\AxiomC{$\Delta_4$}
\AxiomC{$\Delta_5$}
\RightLabel{$\texttt{elim}$}
\TrinaryInfC{$B$}
\end{prooftree}
- where $\Delta_4$ or $\Delta_5$ may be empty. In all cases, by applying a permutation, we obtain a derivation $\Delta^*$ with $\delta(\Delta^*) = (\texttt{MAX}(M_{\Delta^*}), n^*)$ for $\texttt{MAX}(M_{\Delta^*}) = \texttt{MAX}(M_\Delta)$ and $n^* < n$, whence $\delta(\Delta^*) < \delta(\Delta)$. In the first exemplified case, for example, we have

\begin{prooftree}
\AxiomC{$\Delta_1$}
\noLine
\UnaryInfC{$B + C$}
\AxiomC{$\Delta_2$}
\RightLabel{$\texttt{intro}$}
\UnaryInfC{$\textbf{A}$}
\noLine
\UnaryInfC{$\Delta_3$}
\noLine
\UnaryInfC{$\textbf{A}$}
\AxiomC{$\Delta_5$}
\AxiomC{$\Delta_6$}
\RightLabel{$\texttt{elim}$}
\TrinaryInfC{$D$}
\AxiomC{$\Delta_4$}
\noLine
\UnaryInfC{$A$}
\AxiomC{$\Delta_5$}
\AxiomC{$\Delta_6$}
\RightLabel{$\texttt{elim}$}
\TrinaryInfC{$D$}
\RightLabel{($+_E$)}
\TrinaryInfC{$D$}
\end{prooftree}
- observe that the conclusion of $\Delta_4$ may again be the lower edge of a trmcs, whose length is decreased too - while in the second exemplified case we have

\begin{prooftree}
\AxiomC{$\Delta_1$}
\noLine
\UnaryInfC{$T : \exists x \alpha(x)$}
\AxiomC{$\Delta_2$}
\RightLabel{$\texttt{intro}$}
\UnaryInfC{$A$}
\noLine
\UnaryInfC{$\Delta_3$}
\noLine
\UnaryInfC{$A$}
\AxiomC{$\Delta_4$}
\AxiomC{$\Delta_5$}
\RightLabel{$\texttt{elim}$}
\TrinaryInfC{$B$}
\RightLabel{$D_\exists$}
\BinaryInfC{$B$}
\end{prooftree}

\noindent If $\sigma$ has length $1$, then $\sigma$ is either an LMF or a TMF. To give just an example of each case, we may have

\begin{prooftree}
\AxiomC{$[A]$}
\noLine
\UnaryInfC{$\Delta_1$}
\noLine
\UnaryInfC{$B$}
\RightLabel{($\supset_I$)}
\UnaryInfC{$A \supset B$}
\AxiomC{$\Delta_2$}
\noLine
\UnaryInfC{$A$}
\RightLabel{($\supset_E$)}
\BinaryInfC{$B$}
\end{prooftree}
or

\begin{prooftree}
\AxiomC{$\Delta_1$}
\noLine
\UnaryInfC{$T(x) : \alpha(x)$}
\RightLabel{$\forall I$}
\UnaryInfC{$\forall I x (T(x)) : \forall x \alpha(x)$}
\AxiomC{$\Delta_2$}
\noLine
\UnaryInfC{$A$}
\RightLabel{$D_\forall$}
\BinaryInfC{$A$}
\end{prooftree}
In all cases, by applying a reduction, we obtain a derivation $\Delta^*$ for which one of the following circumstances applies:

\begin{itemize}
\item $\Delta^*$ does not contain cut-segments not already occurring in $\Delta$. In this case: either all the cut-segments of $\Delta^*$ have lower measure than the cut-segments of $\Delta$ of maximal measure, so that $\texttt{MAX}(M_{\Delta^*}) < \texttt{MAX}(M_\Delta)$, and hence $\delta(\Delta^*) < \delta(\Delta)$; or $\Delta^*$ contains cut-segments of the same measure of the cut-segments of $\Delta$ of maximal measure, so that $\delta(\Delta^*) = (\texttt{MAX}(M_{\Delta^*}), n^*)$ for $\texttt{MAX}(M_{\Delta^*}) = \texttt{MAX}(M_\Delta)$ and $n^* < n$, whence $\delta(\Delta^*) < \delta(\Delta)$;
\item $\Delta^*$ contains cut-segments which were not already in $\Delta$. If $B$ is the formula $A$ of a new cut-segment, then $\mu(B) < \mu(A)$. This means that, either $\texttt{MAX}(M_{\Delta^*}) < \texttt{MAX}(M_\Delta)$, whence $\delta(\Delta^*) < \delta(\Delta)$, or $\delta(\Delta^*) = (\texttt{MAX}(M_{\Delta^*}), n^*)$ for $\texttt{MAX}(M_{\Delta^*}) = \texttt{MAX}(M_\Delta)$ and $n^* < n$, whence again $\delta(\Delta^*) < \delta(\Delta)$.
\end{itemize}
In the first exemplified case, the reduction yields

\begin{prooftree}
\AxiomC{$\Delta_2$}
\noLine
\UnaryInfC{$A$}
\noLine
\UnaryInfC{$\Delta_1$}
\noLine
\UnaryInfC{$B$}
\end{prooftree}
The only new cut-segment $\sigma^*$ produced by the reduction can have $A$ as formula. Independently of whether $A$ is atomic or logically complex, we have that $\tau(A) \leq \tau(A \supset B)$ and $k^2(A) < k^2(A \supset B)$, whence $\mu(A) < \mu(A \supset B)$, and hence $\mu(\sigma^*) < \mu(\sigma)$. In the second exemplified case, if $\Delta_2$ does not depend on the assumptions $\forall I x (T(x)) \equiv \forall I x (\texttt{h}^{\alpha(x)}(x))$ and $\Pi x (\texttt{h}^{\alpha(x)}(x) : \alpha(x))$, the reduction yields simply

\begin{prooftree}
\AxiomC{$\Delta_2$}
\noLine
\UnaryInfC{$A$}
\end{prooftree}
and no new cut-segment can be produced. If instead $\Delta_2$ depends on the assumption $\forall I x (T(x)) \equiv \forall I x (\texttt{h}^{\alpha(x)}(x))$ but not on the assumption $\Pi x (\texttt{h}^{\alpha(x)}(x) : \alpha(x))$, the reduction produces the derivation

\begin{prooftree}
\AxiomC{$\Delta_3$}
\AxiomC{}
\RightLabel{$\equiv_R$}
\UnaryInfC{$\forall I x (T(x)) \equiv \forall I x (T(x))$}
\AxiomC{$\Delta_4$}
\noLine
\TrinaryInfC{$\Delta_5$}
\noLine
\UnaryInfC{$A$}
\end{prooftree}
which is obtained from $\Delta_2$ by replacing the assumption $\forall I x (T(x)) \equiv \forall I x (\texttt{h}^{\alpha(x)}(x))$ with an application of the reflexivity axiom for equivalence to $\forall I x (T(x))$, and then replacing everywhere $\texttt{h}^{\alpha(x)}(x)$ with $T(x)$. Once again, no new cut-segments can be produced. If finally $\Delta_2$ depends on the assumption $\Pi x (\texttt{h}^{\alpha(x)}(x) : \alpha(x))$ but not on the assumption $\forall I x (T(x)) \equiv \forall I x (\texttt{h}^{\alpha(x)}(x))$, or if it depends on both assumptions, the reduction produces the derivation

\begin{prooftree}
\AxiomC{$\Delta_3$}
\AxiomC{$\Delta_1$}
\noLine
\UnaryInfC{$T(x) : \alpha(x)$}
\RightLabel{($\Pi_I$)}
\UnaryInfC{$\Pi x (T(x) : \alpha(x))$}
\AxiomC{$\Delta_4$}
\noLine
\TrinaryInfC{$\Delta_5$}
\noLine
\UnaryInfC{$A$}
\end{prooftree}
which is obtained from $\Delta_2$ by replacing the assumption $\Pi x (\texttt{h}^{\alpha(x)}(x) : \alpha(x))$ with the derivation

\begin{prooftree}
\AxiomC{$\Delta_1$}
\noLine
\UnaryInfC{$T(x) : \alpha(x)$}
\RightLabel{($\Pi_I$)}
\UnaryInfC{$\Pi x (T(x) : \alpha(x))$}
\end{prooftree}
and then replacing everywhere $\texttt{h}^{\alpha(x)}(x)$ with $T(x)$ - and replacing the assumption $\forall I x (T(x)) \equiv \forall I x (\texttt{h}^{\alpha(x)}(x))$ with an application of the reflexivity axiom for equivalence to $\forall I x (T(x))$. The only new cut-segment $\sigma^*$ produced by the reduction can have $\Pi x (T(x) : \alpha(x))$ as formula. Clearly, $\tau(\Pi x (T(x) : \alpha(x))) < \tau(\forall I x (T(x)) : \forall x \alpha(x))$, so $\mu(\Pi x (T(x) : \alpha(x))) < \mu(\forall I x (T(x)) : \forall x \alpha(x)))$, and hence $\mu(\sigma^*) < \mu(\sigma)$.
\end{proof}

\begin{crl}
For every $\Delta$, there is a normal $\Delta^*$ such that $\Delta \succ \Delta^*$.
\end{crl}

\noindent The normalization theorem is thus a corollary of theorem 17, obtained by repeatedly applying reduction and permutation functions. It should be added that the normal derivation $\Delta^*$ to which $\Delta$ reduces depends on a set of assumptions $\Gamma^*$ and involves unbound individual variables $x_1, ..., x_n$ such that, called $\Gamma$ the set of assumptions on which $\Delta$ depends, and called $y_1, ..., y_m$ the unbound individual variables of $\Delta$, $\Gamma^* \subseteq \Gamma$ and $\{x_1, ..., x_n\} \subseteq \{y_1, ..., y_m\}$. That this is the case is easily seen from how reduction and permutation functions are defined.

Normalization results are usually accompanied by derived theorems concerning the structure of normal derivations, the sub-formula property, and so on. We conclude this section by briefly discussing some of these topics, although we remark that our discussion is not meant to be exhaustive.

As known, in Gentzen's natural deduction systems for, say, intuitionistic or classical - without $\vee$ and $\exists$ - first-order logic (respectively, $\texttt{IL}$ and $\texttt{CL}$) one can prove that paths of normal derivations enjoy some most interesting structural properties [the results we are going to mention are due to Prawitz 2006]. A path in a derivation $\Delta$ in $\texttt{IL}$ or $\texttt{CL}$ is a sequence $a_1, ..., a_n$ of formula occurrences in $\Delta$ such that: (1) $a_1$ is a top-formula in $\Delta$ not discharged by a disjunction or existential elimination; (2) for every $i < n$, $a_i$ is not the minor premise of an implication elimination and, either $a_i$ is not the major premise of a disjunction or existential elimination and $a_{i + 1}$ is the formula occurrence immediately below $a_i$, or $a_i$ is the major premise of a disjunction or existential elimination and $a_{i + 1}$ is an assumption discharged by this elimination; (3) $a_n$ is the minor premise of an implication elimination, or the end formula of $\Delta$, or the major premise of a disjunction or existential elimination that does not discharge any assumption. When $\Delta$ is normal, each path in $\Delta$ can be divided into a sequence $s_1, ..., s_n$ of segments (i.e. of sequences of occurrences of one and the same formula in $\Delta$, like for cut-segments but without requiring that the upper edge is the conclusion of an introduction and the lower edge the major premise of an elimination) where, for some $i < n$, a minimal segment $s_i$ splits the path into two parts: the so-called E-part where, given $j < i$, the formula of each $s_{j + 1}$ is an immediate sub-formula of the formula of $s_j$, and the so-called I-part where, given $i < j < n$, the formula of $s_j$ is an immediate sub-formula of the formula of $s_{j + 1}$. So morally, the minimal segment splits the path into a part where only elimination rules are applied, and a part where only introduction rules are applied; the minimal segment itself only consist of the application of a $\bot$-rule. From this one can finally prove that, if $\Delta$ has assumptions $\Gamma$ and conclusion $\alpha$, then every formula occurring in $\Delta$ is either a sub-formula of $\alpha$ or a sub-formula of some elements of $\Gamma$.

The question is now whether similar results also hold for normal derivations in the systems of grounding we have been discussing so far. These systems can be understood as Gentzen's first-order natural deduction systems with underlying systems given by the equivalence and the typing rules. Although the latter provide means for proving atomic formulas and for drawing consequences from them, the underlying systems \emph{are not} atomic. Clearly, this depends on the type elimination rules, which either mention non-atomic formulas - $D_\rightarrow$ and $D_\forall$ - or allow for non-atomic derivations from discharged assumptions to the minor premise(s). On the other hand, we know that it may be possible to prove results about the form of normal derivations similar to those mentioned above, when $\texttt{IL}$ or $\texttt{CL}$ have a suitable underlying atomic system [see Cellucci 1978 for a survey]. The E-part and the I-part are now separated by a minimal segment where only atomic rules are applied, so that every formula occurrence in the E-part contains its immediate successor as sub-formula, and every formula in the I-part contains its immediate predecessor as sub-formula; it also follows a sub-formula principle limited to the E-part and to the I-part, i.e. every formula occurrence either in the E-part or in the I-part is a sub-formula of the conclusion or of one of the undischarged assumptions. The idea could therefore be that of having something similar in our systems of grounding, i.e. that each path in a normal derivation in such systems splits into an E-part where only elimination rules for the logical constants of the language of grounding are applied, a minimal segment where only equivalence or typing rules are applied, and finally an I-part where only introduction rules for the logical constants of the language of grounding are applied.

Whether this holds or not depends on how we define paths. If we want to stick to a notion of path as the one introduced above, we should say what follows. Given a derivation $\Delta$ in a system of grounding, a path in $\Delta$ is a sequence $a_1, ..., a_n$ of formula occurrences in $\Delta$ such that: (1) $a_1$ is a top-formula in $\Delta$ not discharged by ($+_E$), ($\mathfrak{E}_E$), or by a type elimination; (2) for every $i < n$, $a_i$ is not the minor premise of ($\supset_E$) and, either $a_i$ is not the major premise of ($+_E$), ($\mathfrak{E}_E$) or of a type elimination rule and $a_{i + 1}$ is the formula occurrence immediately below $a_i$, or $a_i$ is the major premise of ($+_E$), ($\mathfrak{E}_E$) or of a type elimination rule and $a_{i + 1}$ is an assumption discharged by this elimination; (3) $a_n$ is the minor premise of ($\supset_E$), or the end formula of $\Delta$, or the major premise of ($+_E$), ($\mathfrak{E}_E$) or of a type elimination that does not discharge any assumption. But now consider the following normal derivation:

\begin{prooftree}
\AxiomC{$1$}
\noLine
\UnaryInfC{$\xi^{\alpha \rightarrow \beta} : \alpha \rightarrow \beta$}
\AxiomC{$2$}
\noLine
\UnaryInfC{$\xi^\alpha : \alpha$}
\AxiomC{$3$}
\noLine
\UnaryInfC{$[\Pi \xi^\alpha (\xi^\alpha : \alpha \supset \texttt{f}^\beta(\xi^\alpha) : \beta)]$}
\RightLabel{($\Pi_E$)}
\UnaryInfC{$\xi^\alpha : \alpha \supset \texttt{f}^\beta(\xi^\alpha) : \beta$}
\RightLabel{($\supset_E$)}
\BinaryInfC{$\texttt{f}^\beta(\xi^\alpha) : \beta$}
\RightLabel{($\mathfrak{E}_I$)}
\UnaryInfC{$\mathfrak{E} \texttt{f}^\beta (\texttt{f}^\beta(\xi^\alpha) : \beta))$}
\RightLabel{$D_\rightarrow$, $3$}
\BinaryInfC{$\mathfrak{E} \texttt{f}^\beta (\texttt{f}^\beta(\xi^\alpha) : \beta))$}
\end{prooftree}
Here, we have a path where the major premise of $D_\rightarrow$ is \emph{followed} by the major premise of ($\Pi_E$). Of course, this depends on the fact that $D_\rightarrow$ - like $D_\forall$ - has a side-assumption which, at variance with its atomic major premise, is logically complex.

One may however expect that, when the normal derivation \emph{does not} involve applications of $D_\rightarrow$ or of $D_\forall$, situations similar to the one exemplified above do not obtain. In such cases, each path in a normal derivation may split into an E-part with only logical eliminations, a central part with only equivalence and typing rules, and an I-part with only logical introductions. One may thus ask whether, similarly to what happens in the E-part and the I-part, in the central part we never have a type-introduction followed by a type elimination. The answer would be in the negative, though, and this mainly depends on the rule of preservation of denotation. Consider for example the following normal derivation:

\begin{tiny}
\begin{prooftree}
\AxiomC{$1$}
\noLine
\UnaryInfC{$\xi^\alpha : \alpha$}
\AxiomC{$2$}
\noLine
\UnaryInfC{$\xi^\beta : \beta$}
\RightLabel{$\wedge I$}
\BinaryInfC{$\wedge I (\xi^\alpha, \xi^\beta) : \alpha \wedge \beta$}
\AxiomC{$3$}
\noLine
\UnaryInfC{$\wedge I (\xi^\alpha, \xi^\beta) \equiv \xi^{\alpha \wedge \beta}$}
\RightLabel{$\equiv_P$}
\BinaryInfC{$\xi^{\alpha \wedge \beta} : \alpha \wedge \beta$}
\AxiomC{$4$}
\noLine
\UnaryInfC{$[\xi^\alpha_1 : \alpha]$}
\AxiomC{}
\RightLabel{$\equiv_\wedge$}
\UnaryInfC{$\wedge_{E, 1}(\wedge I (\xi^\alpha_1, \xi^\beta_1)) \equiv \xi^{\alpha}_1$}
\RightLabel{$\equiv_P$}
\BinaryInfC{$\wedge_{E, 1}(\wedge I (\xi^{\alpha}_1, \xi^\beta_1)) : \alpha$}
\AxiomC{$5$}
\noLine
\UnaryInfC{$[\xi^{\alpha \wedge \beta} \equiv \wedge I (\xi^\alpha_1, \xi^\beta_1)]$}
\RightLabel{$\equiv^\wedge_{3, 1}$}
\UnaryInfC{$\wedge_{E, 1}(\xi^{\alpha \wedge \beta}) \equiv \wedge_{E, 1}(\wedge I(\xi^\alpha_1, \xi^\beta_1))$}
\RightLabel{$\equiv_P$}
\BinaryInfC{$\wedge_{E, 1}(\xi^{\alpha \wedge \beta}) : \alpha$}
\RightLabel{$D_\wedge$, $4, 5$}
\BinaryInfC{$\wedge_{E, 1}(\xi^{\alpha \wedge \beta}) : \alpha$}
\end{prooftree}
\end{tiny}
Here, we have two paths where a premise of $\wedge I$ is \emph{followed} by the major premise of $D_\wedge$, without this creating a TMF because of the intermediate leftmost application of $\equiv_P$.\footnote{One may also expect that a similar situation may occur because of an application of ($\bot^G$), but in a normal derivation this cannot happen, for it generally holds what follows. Let $\pi$ be a path in a normal derivation, and let $s_1, ..., s_n$ be the segments of $\pi$ whose formulas are respectively $A_1, ..., A_n$. If, for some $i < j \leq n$ there are $A_i, A_j$ such that $A_i = T : \alpha$ and $A_j = U : \beta$ with $k^1(\alpha) > k^2(\beta)$, then there is $i < h < j$ such that $A_h$ is the major premise of a type elimination rule. This suggests in turn that an order in normal derivations may be found not with respects to formulas, but with respect to \emph{types}, i.e. the type of each term occurring in a normal derivation is a sub-formula of the type of some term which occurs either in some sub-formula of the conclusion, or in some sub-formula of some element of the set of undischarged assumptions. In particular, given a path $\pi$ in a normal derivation without applications of $D_\rightarrow$ or $D_\forall$ (or defined in such a way that each path splits into an E-part, an I-part and a central part with only equivalence and typing rules), and given the sequence of segments $s_1, ..., s_n$ in its central part whose formulas are respectively $A_1, ..., A_n$, for every $i \leq n$ there is $j \leq n$ such that either the type of the term(s) occurring in $A_i$ is an immediate sub-formula of the type of the term(s) occurring in $A_j$, or vice versa. We may then call \emph{maximum type-point} any $A_i$ such that the type of the term(s) occurring in $A_i$ is not the immediate sub-formula of the type of the term(s) of any $A_j$. The problem of finding an order in the central part of paths of normal derivations becomes then the problem of whether these paths can have a \emph{single} maximum type-point.}

One may envisage different ways to deal with the problems just discussed, and seek whether it is possible to solve them in such a way as to obtain results like those mentioned above about the form and the sub-formulas of normal derivations. For example, one may require that a path can pass through an assumption discharged by a type elimination only when this assumption is atomic - and hence allow a path to start from a non-atomic assumption discharged by a type elimination. As for $\equiv_P$, one may require that passages where a type-introduction ends in one of the premises of $\equiv_P$ \emph{do} give rise to maximal points, and then try to define suitable reduction functions for removing the application of $\equiv_P$. A discussion of whether and how these strategies are feasible would however lead us too far - but may constitute the topic of future works.

\section{Conclusions}

In these concluding remarks, we want to raise a quick point concerning Prawitz's conjecture [Prawitz 1973], namely, about whether first-order intuitionistic logic is complete with respect to Prawitz's semantics.

Strictly speaking, Prawitz's original conjecture concerns the semantics of valid arguments. It states that, if an inference rule is logically valid in this semantics, then it is derivable in first-order intuitionistic logic [Prawitz 1973]. A ground-theoretic reformulation of the conjecture requires a notion of universal (operation on) ground(s). An (operation on) ground(s) $g$ is universal iff, for every atomic base $\mathfrak{B}$, $g$ is an (operation on) ground(s) over $\mathfrak{B}$. An inference from $\Gamma$ to $\alpha$ is then logically valid iff there is a universal operation on grounds of type $\Gamma \rhd \alpha$ [see d'Aragona 2021a for further details].

Prawitz's conjecture for valid arguments - in a weaker form - has been refuted by Piecha and Schroeder-Heister [Piecha \& Schroeder-Heister 2018, see also Piecha, de Campos Sanz \& Schroeder-Heister 2015], and [as shown for example in d'Aragona 2019b] Piecha and Schroeder-Heister's proof can be adapted to the theory of grounds. From a ground-theoretic point of view, the weaker form relies upon the following definition of logical validity: an inference from $\Gamma$ to $\alpha$ is logically valid iff, for every atomic base $\mathfrak{B}$, there is an operation on grounds over $\mathfrak{B}$ of type $\Gamma \rhd \alpha$.

One can also put forward a strong completeness conjecture, stating that \emph{every} universal operation on grounds $\phi$ of type $\Gamma \rhd \alpha$ corresponds to a term $T$ of $\texttt{Gen}$ such that the denotation of $T$ is equivalent to $\phi$ [see d'Aragona 2019b for details]. The strong completeness conjecture may not be far-fetched, even if Prawitz's conjecture turned out to be false in the stronger form. Indeed, Piecha and Schroeder-Heister [Piecha \& Schroeder-Heister 2018] suggest that completeness may hold for some intermediate logic, and so they reformulate Prawitz's conjecture in terms of existence of an appropriate intermediate logic. The strong completeness conjecture may then apply to such intermediate logic.

The strong completeness conjecture [as formulated e.g. in d'Aragona 2019b], has a seemingly interesting consequence on our systems of grounding. It is equivalent to the following proposition.

\begin{prp}
Let $\texttt{Gen}^*$ be an expansion of $\texttt{Gen}$ obtained by adding non-primitive operational symbols $F_1[\omega_1], ..., F_n[\omega_n]$, and let $\Sigma^+$ be a system of grounding over an enriched version of $\texttt{Gen}^*$. If the equations of $\Sigma^+$ ruling the non-primitive operational symbols of $\texttt{Gen}^*$ are such that $F_i[\omega_i]$ stands for a universal operation on grounds ($1 \leq i \leq n$), then there is a system of grounding $\Sigma$ over an enriched version of $\texttt{Gen}$ such that, for every $T \in \texttt{TERM}_{\texttt{Gen}^*}$, there is $U \in \texttt{TERM}_{\texttt{Gen}}$ such that $\vdash_{\Sigma \cup \Sigma^+} T \equiv U$.
\end{prp}
\noindent Since strong completeness is formulated in terms of denotation functions, a proof of the equivalence requires to connect denotation functions with equations ruling the non-primitive operational symbols. This can be done by first of all introducing a notion of language of grounding $\Lambda$ \emph{self-contained} with respect to a denotation function $den$, namely, for every operational symbol $F$ of $\Lambda$, if $den(F)$ is defined by an equation of the form $den(F) = f_1 \circ ... \circ f_n(x_1, ..., x_n)$ then, for every $f_i$ ($1 \leq i \leq n$), there is an operational symbol $H$ of $\Lambda$ such that $den(H) = f_i$. Secondly, one says that a system of grounding $\Sigma$ over a language of grounding $\Lambda$ \emph{totally interprets} $\Lambda$ with respect to a denotation function $den$ iff the equations of $\Sigma$ ruling the non-primitive operational symbols $F$ of $\Lambda$ ‘‘internalize" the equations that define $den(F)$. Then, given a language of grounding $\Lambda$, one builds a self-contained expansion $\Lambda^*$ of it, and defines a system of grounding $\Sigma$ that totally interprets $\Lambda^*$. With this done, it is easily seen that strong completeness implies ‘‘re-writability" in $\Sigma$ plus an appropriate system of grounding over $\texttt{Gen}$ (along the lines of what shown in Section 4.2), and vice versa [for more precise technical details, see again d'Aragona 2019b].

\paragraph{Acknowledgments} I am grateful to Cesare Cozzo, Gabriella Crocco, Enrico Moriconi and, above all, Dag Prawitz for helpful suggestions. I am also grateful to the anonymous reviewers, whose comments helped me to improve an earlier draft of this paper.

\begin{footnotesize}
\paragraph{References}
\begin{itemize}
\item C. Cellucci (1978), \emph{Teoria delle dimostrazione. Normalizzazioni e assegnazioni di numeri ordinali}, Bollati Boringhieri, Torino.
\item C. Cozzo (1994), \emph{Meaning and argument. A theory of meaning centred on immediate argumental role}, Almqvist \& Wiksell, Uppsala.
\item A. P. d'Aragona (2018), \emph{A partial calculus for Dag Prawitz's theory of grounds and a decidability issue}, in in A. Christian, D. Hommen, N. Retzlaff, G. Schurz (eds), \emph{Philosophy of Science. European Studies in Philosophy of Science}, vol 9. Springer, Berlin Heidelberg New York.
\item --- (2019a), \emph{Dag Prawitz on proofs, operations and grounding}, in G. Crocco \& A. P. d'Aragona, \emph{Inferences and proofs}, special issue of \emph{Topoi}, DOI: https://doi.org/10.1007/s11245-017-9473-9.
\item --- (2019b), \emph{Dag Prawitz's \emph{theory of grounds}}, PhD dissertation, Aix-Marseille University, ‘‘Sapienza" University of Rome, HAL Id: tel-02482320, version 1.
\item --- (2021a), \emph{Denotational semantics for languages of epistemic grounding based on Prawitz's theory of grounds}, in \emph{Studia logica}, DOI: https://doi.org/10.1007/s11225-021-09969-8.
\item --- (2021b), \emph{Proofs, grounds and empty functions: epistemic compulsion in Prawitz's semantics}, in \emph{Journal of philosophical logic}, DOI: https://doi.org/10.1007/s10992-021-09621-9.
\item M. Dummett (1991), \emph{The logical basis of metaphysics}, Harvard University Press, Cambridge.
\item N. Francez (2015), \emph{Proof-theoretic semantics}, College Publications, London.
\item G. Gentzen (1934 - 1935), \emph{Untersuchungen \"{u}ber das logische Schlie\ss{}en}, in \emph{Matematische Zeitschrift}, XXXIX, DOI: BF01201353.
\item A. Heyting (1956), \emph{Intuitionism. An introduction}, North-Holland Publishing Company, Amsterdam.
\item P. Martin-L\"{o}f (1984), \emph{Intuitionistic type theory}, Bibliopolis, Napoli.
\item T. Piecha, W. de Campos Sanz \& P. Schroeder-Heister (2015), \emph{Failure of completeness in proof-theoretic semantics}, in \emph{Journal of philosophical logic}.
\item T. Piecha \& P. Schroeder-Heister (2018), \emph{Incompleteness of intuitionistic propositional logic with respect to proof-theoretic semantics}, in \emph{Studia Logica}.
\item D. Prawitz (1971), \emph{Ideas and results in proof theory}, in \emph{Proceedings of the Second Scandinavian Logic Symposium}.
\item --- (1973), \emph{Towards a foundation of a general proof-theory}, in P. Suppes, (ed) \emph{Logic methodology and philosophy of science IV}, North-Holland Publishing Company, Amsterdam, DOI: S0049-237X(09)70361-1.
\item --- (1977), \emph{Meaning and proofs: on the conflict between classical and intuitionistic logic}, in \emph{Theoria}, DOI: j.1755-2567.1977.tb00776.x.
\item --- (2006), \emph{Natural deduction. A proof-theoretical study}, Dover, New York.
\item --- (2015), \emph{Explaining deductive inference}, in H. Wansing (ed), \emph{Dag Prawitz on proofs and meaning}, Springer, Berlin Heidelberg New York, DOI: 978-3-319-11041-7\_3.
\item --- (2019), \emph{The seeming interdependence between the concepts of valid inference and proof}, in G. Crocco \& A. Piccolomini d'Aragona (eds), \emph{Inferences and proofs}, special issue of \emph{Topoi}, DOI: s11245-017-9506-4.
\item --- (2021), \emph{The validity of inference and argument}, forthcoming.
\item P. Schroeder-Heister (1984a), \emph{A natural extension for natural deduction}, in \emph{Journal of symbolic logic}.
\item --- (1984b), \emph{Generalized rules for quantifiers and the completeness of the intuitionistic operators $\wedge$, $\vee$, $\rightarrow$, $\forall$, $\exists$}, in M. M. Richter, E. B\"{o}rger, W. Oberschelp, B. Schinzel \& W. Thomas (eds.), \emph{Computation and proof theory. Proceedings of the Logic Colloquium held in Aachen, July 18-23, 1983, Part II}, Springer, Berlin Heidelberg New York Tokyo.
\item --- (1991), \emph{Uniform proof-theoretic semantics for logical constants. Abstract}, in \emph{Journal of Symbolic Logic}.
\item --- (2018), \emph{Proof-theoretic semantics}, in Edward N. Zalta (ed.), \emph{The Stanford Encyclopedia of Philosophy \emph{(}Spring 2018 Edition\emph{)}}.
\item L. Tranchini (2019), \emph{Proof, meaning and paradox: some remarks}, in G. Crocco \& A. Piccolomini d'Aragona (eds), \emph{Inferences and proofs}, special issue of \emph{Topoi}, DOI: s11245-018-9552-6.
\item G. Usberti (2015), \emph{A notion of $C$-justification for empirical statements}, in H. Wansing (ed), \emph{Dag Prawitz on proofs and meaning}, Springer, Berlin Heidelberg New York, DOI: 978-3-319-11041-7\_18.
\item --- (2019), \emph{Inference and epistemic transparency}, in G. Crocco \& A. Piccolomini d'Aragona (eds), \emph{Inferences and proofs}, special issue of \emph{Topoi}, DOI: s11245-017-9497-1.
\end{itemize}
\end{footnotesize}

\end{document}